\documentclass[11pt,twoside]{article}
\usepackage{fullpage}

\usepackage{epsf}
\usepackage{fancyhdr}
\usepackage{microtype}

\usepackage{algorithmic}
\usepackage{color,xcolor}
\usepackage[linesnumbered,ruled]{algorithm2e}
\DontPrintSemicolon
\usepackage{color}

\usepackage{amsthm}
\usepackage{amsfonts}
\usepackage{amsmath}
\usepackage{amssymb,bbm}
\usepackage{comment}

\newcommand*\Uspace[0]{\mathbb{U}}

\def\ball{\mathbb{B}}

\newcommand{\real}{\ensuremath{\mathbb{R}}}

\newcommand{\Xspace}{\ensuremath{\mathbb{X}}}

\newcommand{\abss}[1]{\left| #1 \right |}

\newcommand{\law}{\ensuremath{\mathcal{L}}}

\newcommand{\mydefn}{\ensuremath{:=}}

\newcommand{\defn}{:=}

\newcommand{\matsnorm}[2]{|\!|\!| #1 | \! | \!|_{{#2}}}
\newcommand{\vecnorm}[2]{\| #1\|_{#2}}

\newcommand{\opnorm}[1]{\ensuremath{\matsnorm{#1}{\tiny{\mbox{op}}}}}

\newcommand{\inprod}[2]{\ensuremath{\langle #1 , \, #2 \rangle}}

\newcommand{\kull}[2]{\ensuremath{D_{\text{KL}}\left(#1\; \| \; #2 \right)}}

\newcommand{\Exs}{\ensuremath{{\mathbb{E}}}}
\newcommand{\Prob}{\ensuremath{{\mathbb{P}}}}

\DeclareMathOperator{\diag}{diag}
\DeclareMathOperator{\var}{var}

\DeclareMathOperator{\trace}{trace}

\newtheoremstyle{named}{}{}{\itshape}{}{\bfseries}{.}{.5em}{\thmnote{#3's }#1}
\theoremstyle{named}

\theoremstyle{plain}

\newtheorem{theorem}{Theorem}
\newtheorem{proposition}{Proposition}

\newtheorem{lemma}{Lemma}

\newtheorem{corollary}{Corollary}
\newtheorem{definition}{Definition}

\newlength{\widebarargwidth}
\newlength{\widebarargheight}
\newlength{\widebarargdepth}
\DeclareRobustCommand{\widebar}[1]{%
  \settowidth{\widebarargwidth}{\ensuremath{#1}}%
  \settoheight{\widebarargheight}{\ensuremath{#1}}%
  \settodepth{\widebarargdepth}{\ensuremath{#1}}%
  \addtolength{\widebarargwidth}{-0.3\widebarargheight}%
  \addtolength{\widebarargwidth}{-0.3\widebarargdepth}%
  \makebox[0pt][l]{\hspace{0.3\widebarargheight}%
    \hspace{0.3\widebarargdepth}%
    \addtolength{\widebarargheight}{0.3ex}%
    \rule[\widebarargheight]{0.95\widebarargwidth}{0.1ex}}%
  {#1}}

\makeatletter
\long\def\@makecaption#1#2{
        \vskip 0.8ex
        \setbox\@tempboxa\hbox{\small {\bf #1:} #2}
        \parindent 1.5em  
        \dimen0=\hsize
        \advance\dimen0 by -3em
        \ifdim \wd\@tempboxa >\dimen0
                \hbox to \hsize{
                        \parindent 0em
                        \hfil
                        \parbox{\dimen0}{\def\baselinestretch{0.96}\small
                                {\bf #1.} #2
                                }
                        \hfil}
        \else \hbox to \hsize{\hfil \box\@tempboxa \hfil}
        \fi
        }
\makeatother

\long\def\comment#1{}
\definecolor{battleshipgrey}{rgb}{0.52, 0.52, 0.51}
\definecolor{darkgray}{rgb}{0.66, 0.66, 0.66}
\definecolor{darkgreen}{rgb}{0.0, 0.2, 0.13}
\definecolor{darkspringgreen}{rgb}{0.09, 0.45, 0.27}
\definecolor{dukeblue}{rgb}{0.0, 0.0, 0.61}
\definecolor{olivedrab7}{rgb}{0.24, 0.2, 0.12}
\definecolor{darkblue}{rgb}{0.0, 0.0, 0.55}
\definecolor{darkscarlet}{rgb}{0.34, 0.01, 0.1}
\definecolor{candyapplered}{rgb}{1.0, 0.03, 0.0}
\definecolor{ao(english)}{rgb}{0.0, 0.5, 0.0}
\definecolor{applegreen}{rgb}{0.55, 0.71, 0.0}

\newcommand{\widgraph}[2]{\includegraphics[keepaspectratio,width=#1]{#2}}

\newcommand{\avgset}{\widehat{S}_{\xzero, \actionzero}}

\newcommand{\vbar}{\widebar{v}}
\newcommand{\rkhs}{\mathbb{H}}

\usepackage{url}
\usepackage[colorlinks,linkcolor=magenta,citecolor=blue, backref=true]{hyperref}

\usepackage{nicefrac}
\usepackage{comment}

\newcommand{\myunder}[2]{\underbrace{{#1}}_{{\mbox{\small{#2}}}}}
\newcommand{\myundermath}[2]{\underbrace{{#1}}_{#2}}

\usepackage{scalerel}
\newcommand{\smallsuper}[2]{#1^{\scaleto{#2}{4pt}}}
\newcommand{\smallsub}[2]{#1_{\scaleto{#2}{4pt}}}
\newcommand{\medsuper}[2]{#1^{\scaleto{#2}{6pt}}}
\newcommand{\proplogistic}{\propscore^{\mathrm{L}}}
\newcommand{\propnormal}{\propscore^{\mathrm{N}}}
\newcommand{\probxlogistic}{\probx^{\mathrm{L}}}
\newcommand{\probxnormal}{\probx^{\mathrm{N}}}
\newcommand{\probxcauchy}{\probx^{\mathrm{C}}}

\usepackage{chngpage}

\usepackage[mathscr]{euscript}
 \usepackage{tabularx}%

\usepackage{enumitem}
\setlist[enumerate,0]{leftmargin=5em}

\usepackage{booktabs}

\usepackage{cleveref}

\newcommand{\ltwolfour}[1]{\ensuremath{\|#1\|_{2 \rightarrow 4}}}

\usepackage{bm,bbm}
\usepackage{mathtools}

\newcommand{\ctwofour}{c_{2 \rightarrow 4}}

\newcommand{\regu}{\rho}

\newcommand{\datablock}{\mathcal{B}}

\newcommand{\actionspace}{\mathbb{A}}

\newcommand{\varbound}{\bar{\sigma}}

\newcommand{\dudley}{\mathcal{J}}

\newcommand{\statespace}{\mathbb{X}}

\newcommand{\actionzero}{\action_0}
\newcommand{\minimaxRisk}{\mathscr{M}_\numobs}

\newcommand{\tauhat}{\widehat{\avgtreat}}
\newcommand{\Ltwospace}{\mathbb{L}^2}

\newcommand{\lowersigma}{\underline{\sigma}}

\newcommand{\kernelfunc}{\mathcal{K}}

\newtheorem{assumption}{Assumption}
\newcommand{\numobs}{\ensuremath{n}}

\newcommand{\usedim}{\ensuremath{d}}

\newcommand{\basisfunc}{\phi}

\newcommand{\smallexponent}{\nu}

\newcommand{\sigbar}{\overline{\sigma}}

\newcommand{\effdim}{D}

\newcommand{\eigen}{\lambda}
\newcommand{\eigval}{\eigen}

\newcommand{\totalvariation}{d_{\mathrm{TV}}}
\newcommand{\radiusnumobs}{r_\numobs}

\newcommand{\plainvec}{z}

\newcommand{\abstrEst}{\mathcal{A}}

\newcommand{\rkhsinprod}[2]{\inprod{#1}{#2}_{\rkhs}}
\newcommand{\rkhsnorm}[1]{\vecnorm{#1}{\rkhs}}

\newcommand{\rkhsprime}{\rkhs_\sigma}

\newcommand{\rkhsprimenorm}[1]{\vecnorm{#1}{\rkhsprime}}
\newcommand{\kernelfuncprime}{\kernelfunc_\sigma}

\newcommand{\probx}{\xi}

\newcommand{\highorder}{\mathcal{H}}

\newcommand{\avgtreat}{\tau}

\newcommand{\pitarget}{T}
\newcommand{\bigcov}{\mymatrix{\Gamma}}

\newcommand{\Ball}{\ensuremath{\mathbb{B}}}

\newcommand{\maux}{m}

\newcommand{\myiid}{\mathrm{i.i.d.}}

\newcommand{\featweighbar}{\widebar{\featweigh}}

\newcommand{\Event}{\mathscr{E}}

\newcommand{\chisqdiv}[2]{\chi^2 \left( #1 ~||~ #2 \right)}
\long\def\comment#1{}

\long\def\comment#1{}

\newcommand{\funcClass}{\mathcal{F}}
\newcommand{\state}{s}

\newcommand{\probInstance}{\mathcal{I}}

\newcommand{\weightfunc}{\omega}

\newcommand{\propscore}{\pi}
\newcommand{\treateff}{\plaintreateff^*}

\newcommand{\plaintreateff}{\mu}

\newcommand{\outcome}{Y}

\newcommand{\localradius}{r_0}

\newcommand{\dsone}[1]{\medsuper{#1}{$(I)$}}
\newcommand{\dstwo}[1]{\medsuper{#1}{$(II)$}}
\newcommand{\dsthree}[1]{\medsuper{#1}{$(III)$}}
\newcommand{\dsfour}[1]{\medsuper{#1}{$(IV)$}}

\newcommand{\torus}{\mathbb{T}}

\newcommand{\Rbound}{\ensuremath{R}}

\newcommand{\smoothorder}{s}

\newcommand{\IdMat}{\ensuremath{\mymatrix{I}}}

\newcommand{\Nplain}[1]{\neighborhood_{#1}}

\newcommand{\Nval}[1]{\smallsuper{\neighborhood}{val}_{#1}}

\newcommand{\noise}{\varepsilon}
\newcommand{\smallbias}{b}

\newcommand{\offpar}{\ensuremath{s}}

\newcommand{\neighborhood}{\mathscr{N}}

\newcommand{\taustar}{\tau^*}

\newcommand{\sighat}{\widehat{\sigma}}

\newcommand{\xzero}{\ensuremath{{x_0}}}

\newcommand{\Delhat}{\widehat{\Delta}}

\newcommand{\dimx}{\usedim_x}
\newcommand{\dima}{\usedim_a}

\newcommand{\mutil}{\widetilde{\mu}}

\newcommand{\naturals}{\mathbb{N}}

\newcommand{\plainmu}{\plaintreateff}

\newcommand{\muhat}{\widehat{\mu}}

\newcommand{\vecspace}{\mathbb{V}}
\newcommand{\phimax}{\phi_{\max}}

\newcommand{\featweigh}{u}
\newcommand{\projparam}{\beta}
\newcommand{\featweighhat}{\widehat{\featweigh}}
\newcommand{\projparamhat}{\widehat{\projparam}}

\newcommand{\subgaussian}{\gamma}
\renewcommand{\state}{\ensuremath{x}}
\newcommand{\action}{\ensuremath{a}}
\newcommand{\State}{\ensuremath{X}}
\newcommand{\Action}{\ensuremath{A}}
\newcommand{\StateSpace}{\Xspace}
\newcommand{\ActionSpace}{\ensuremath{\mathbb{A}}}

\newcommand{\Outcome}{\ensuremath{Y}}

\newcommand{\probxstar}{\probx^*}

\newcommand{\bigcovsig}{\ensuremath{\bigcov_\sigma}}

\newcommand{\mymatrix}[1]{\ensuremath{\mathbf{#1}}}

\newcommand{\EigenMat}{\mymatrix{\Lambda}}

\usepackage[linewidth=1pt]{mdframed}

\usepackage[firstinits=true,
  maxnames = 99,
  backref=true,
  backend=bibtex,
  style=alphabetic,
  sortlocale=de_DE,
  natbib=true,
  doi=false,
  isbn=false,
  url=false]{biblatex}

\newenvironment{carlist}
 {\begin{list}{$\bullet$}
 {\setlength{\topsep}{0in} \setlength{\partopsep}{0in}
  \setlength{\parsep}{0in} \setlength{\itemsep}{\parskip}
  \setlength{\leftmargin}{0.15in} \setlength{\rightmargin}{0.08in}
  \setlength{\listparindent}{0in} \setlength{\labelwidth}{0.08in}
  \setlength{\labelsep}{0.1in} \setlength{\itemindent}{0in}}}
 {\end{list}}

\newcommand{\bcar}{\begin{carlist}}
\newcommand{\ecar}{\end{carlist}}

\newenvironment{narrowpara}
  {\par\addvspace{\smallskipamount}\narrower\noindent\ignorespaces}
  {\par\addvspace{\smallskipamount}}

  \newcommand{\ALGSTEPS}[1]{\begin{mdframed}
        \begin{narrowpara} #1
        \end{narrowpara}
    \end{mdframed}
  }
  
\newcommand{\figdir}{figs}

\newcommand{\Fclass}{\ensuremath{\mathcal{F}}}

\newcommand{\overexp}{\ensuremath{\alpha}}

\newcommand{\Vbase}{V}
\newcommand{\vsemi}{\ensuremath{\smallsub{\Vbase}{semi}}}
\newcommand{\probInstanceStar}{\ensuremath{\probInstance^*}}

\newcommand{\vprobx}{\ensuremath{\Vbase_{\probxstar}}}
\newcommand{\vclassical}{\ensuremath{\Vbase_\sigma}}
\newcommand{\impratio}{\ensuremath{\frac{d \weightfunc}{d \propscore}}}
\newcommand{\smallimpratio}{\ensuremath{\tfrac{d \weightfunc}{d \propscore}}}

\newcommand{\vwenplain}{\ensuremath{\Vbase_{\sigma, \numobs}}}
\newcommand{\vwenplaintil}{\ensuremath{\tilde{\Vbase}_{\sigma, \numobs}}}

\newcommand{\vdouble}[2]{\ensuremath{\Vbase_{\sigma,
      \numobs}(#1, \mypair; #2)}}

\newcommand{\vdoublesq}[2]{\ensuremath{\Vbase^2_{\sigma,
      \numobs}(#1, \mypair; #2)}}

\newcommand{\vsuper}[1]{\ensuremath{\Vbase_{\sigma,
      \numobs}(#1, \mypair; \funcClass)}}

\newcommand{\vsupersq}[1]{\ensuremath{\Vbase^2_{\sigma, \numobs}(#1,
    \mypair; \funcClass)}}

    \newcommand{\vsuperhil}[1]{\ensuremath{\Vbase_{\sigma,
      \numobs}(#1, \mypair; \Ball_\Hil(R))}}

\newcommand{\mypair}{\ensuremath{\propscore, \weightfunc}}

\newcommand{\vwenlong}{\ensuremath{\vwenplain(\probxstar, \mypair; \funcClass)}}
\newcommand{\vwenlongsq}{\ensuremath{\vwenplain^2(\probxstar, \mypair;
    \funcClass)}}

\newcommand{\vwenlongxzero}{\ensuremath{\vsuper{\delta_\xzero}}}

\newcommand{\vwenlongsqxzero}{\ensuremath{\vsupersq{\delta_\xzero}}}

\newcommand{\vwenlonghil}{\ensuremath{\vwenplain(\probxstar, \mypair; \Ball_\Hil(R))}}

\newcommand{\vwenlonghilbar}{\ensuremath{\Vbase_{\bar{\sigma},
      \numobs} (\probxstar, \mypair; \Ball_\Hil(R))}}

\newcommand{\hackhilsq}{\ensuremath{\vwenplaintil^2(\mypair;
    \Ball_\Hil(R))}}

\newcommand{\Hil}{\ensuremath{\mathscr{H}}}

\newcommand{\wenprop}[2]{\ensuremath{\propscore(#1 \mid #2)}}

\newcommand{\plainweight}{\ensuremath{\omega}}
\newcommand{\wenweight}[2]{\ensuremath{\plainweight(#1 \mid #2)}}

\newcommand{\LinFunc}{\ensuremath{\mathscr{L}_{\plainweight}}}
\newcommand{\LinFuncZero}{\LinFunc (\treateff, \delta_{\xzero})}

\newcommand{\Joint}[2]{\ensuremath{#1 \cdot  #2}}
\newcommand{\optfun}{\ensuremath{f}}

\newcommand{\newxdist}{\ensuremath{\nu}}

\newcommand{\taustarzero}{\ensuremath{\taustar_\xzero}}

\newcommand{\myoptfun}{\ensuremath{f}}

\newcommand{\pstar}{\ensuremath{\mathbb{P}^*}}
\newcommand{\noisevar}{\ensuremath{W}}

\newcommand{\wenbias}{\ensuremath{t}}
\newcommand{\Wensamp}{\ensuremath{M}}

\newcommand{\Set}{\ensuremath{\mathcal{S}}}

\newcommand{\eigfun}{\ensuremath{\phi}}

\renewcommand{\Uspace}{\ensuremath{\mathbb{U}}}
\newcommand{\uvar}{\ensuremath{u}}

\newcommand{\Hball}{\ensuremath{\Ball_\Hil(\Rbound)}}

\newcommand{\newAmat}{\ensuremath{\mymatrix{A}}}
\newcommand{\Bmat}{\ensuremath{\mymatrix{B}}}
\newcommand{\Mmat}{\ensuremath{\mymatrix{M}}}

\newcommand{\reguhomo}{\ensuremath{\regu_\numobs}}
\newcommand{\regustageone}{\ensuremath{\medsuper{\regu_{\numobs}}{$(I)$}}}
\newcommand{\regustagethree}{\ensuremath{\medsuper{\regu_{\numobs}}{$(III)$}}}

\begin{document}

\begin{center}
  {\bf{\Large{Kernel-based off-policy estimation without overlap: \\
        Instance optimality beyond semiparametric efficiency}}}

\vspace*{.2in} {\large{
 \begin{tabular}{cccc}
 Wenlong Mou$^{ \diamond}$ & 
  Peng Ding$^{\dagger}$  &
  Martin J. Wainwright$^{\diamond, \dagger,
    \ddagger}$ &Peter L.  Bartlett$^{\diamond, \dagger, \star}$
 \end{tabular}
}

\vspace*{.2in}

 \begin{tabular}{c}
 Department of Electrical Engineering and Computer
 Sciences$^\diamond$\\ Department of Statistics$^\dagger$ \\ UC
 Berkeley
 \end{tabular}

 \medskip
 
 \begin{tabular}{c}
   Laboratory for Information and Decision Systems$^\ddagger$
   \\
   Statistics and Data Science Center$^\ddagger$ \\
   EECS and Mathematics \\
   Massachusetts Institute of Technology
 \end{tabular}
  \medskip
 
 \begin{tabular}{c}
    Google Research$^\star$
 \end{tabular}
}

\vspace*{.2in}

\begin{abstract}
We study optimal procedures for estimating a linear functional based
 on observational data. In many problems of this kind, a widely used assumption is strict overlap, i.e., uniform boundedness of the importance ratio, which measures how well the observational data covers the directions of interest.  When it is
 violated, the classical semi-parametric efficiency bound can easily become infinite,
 so that the instance-optimal risk depends on the function class used
 to model the regression function.  For any convex and symmetric
 function class $\mathcal{F}$, we derive a non-asymptotic local
 minimax bound on the mean-squared error in estimating a broad class
 of linear functionals.  This lower bound refines the classical
 semi-parametric one, and makes connections to moduli of continuity in
 functional estimation.  When $\mathcal{F}$ is a reproducing kernel
 Hilbert space, we prove that this lower bound can be achieved up to a
 constant factor by analyzing a computationally simple regression
 estimator.  We apply our general results to various families of
 examples, thereby uncovering a spectrum of rates that interpolate
 between the classical theories of semi-parametric efficiency (with
 $\sqrt{n}$-consistency) and the slower minimax rates associated with
 non-parametric function estimation.
\end{abstract}

\end{center}


\section{Introduction}
\label{SecIntro}

Estimation and inference problems based on observational data arise in
various applications, and are studied in the fields of causal
inference, econometrics, and reinforcement learning.  An interesting
subclass of such problems are semi-parametric in nature: they involve
estimating the value of a linear functional in the presence of one or
more unknown non-parametric ``nuisance'' functions.

More concretely, suppose that we observe $\numobs$ i.i.d. triples of
the form $(\State_i, \Action_i, \Outcome_i)$, where each triple is drawn
according to the following procedure
\bcar
\item the state variable $\State_i$ is drawn from some distribution
  $\probxstar$ over the state space $\StateSpace$.
\item the action variable $\Action_i$ is drawn with conditional
  distribution $\Action_i \mid \State_i \sim
  \wenprop{\cdot}{\State_i}$, where $\propscore$ is a \emph{behavioral
  policy}, also known as the propensity score in the causal inference
  literature.
\item the response or outcome variable $\Outcome_i$ has conditional
  expectation $\Exs \big[ \Outcome_i \mid \State_i, \Action_i \big] =
  \treateff(\State_i, \Action_i)$, where $\treateff$ is the
  \emph{regression function}, also known as the treatment effect.
\ecar

\medskip

The distribution $\probxstar$ and regression function $\treateff$ are
both unknown, and we would like to estimate a known functional
that depends on both of them.  More precisely, given a known family of
functions $ \{ \wenweight{\cdot}{\state} \mid \state \in \Xspace \}$,
where each $\wenweight{\cdot}{\state}$ is a signed Radon measure over
the action space $\actionspace$, consider the functional
\begin{align}
\label{EqnDefnFunctional}  
(\plainmu, \probx)\mapsto \LinFunc(\plainmu, \probx) \mydefn
\Exs_{\State \sim \probx} \left[ \int_\actionspace \treateff (\State,
  \action) d \wenweight{\action}{\State} \right] \; = \;
\int_\Xspace \int_\actionspace \treateff (\state, \action) \; d
\wenweight{\action}{\state} \, d \probx(\state).
\end{align}
Our goal is to estimate $\taustar \defn \LinFunc(\treateff,
\probxstar)$---the value of this functional at the unknown pair
$(\treateff, \probxstar)$.  The behaviorial policy $\propscore$ is
also unknown, and it plays the role of another non-parametric
nuisance, since it affects the joint distribution of the samples
$(\State_i, \Action_i)$ that we observe.  Special cases of this set-up
include estimating the average treatment effect (ATE), and off-policy
evaluation for contextual bandit problems.  We also consider a variant
in which, instead of taking the expectation over $\State \sim \probx$,
we evaluate at a fixed state $\state_0$.  This latter set-up is
appropriate for the conditional average treatment effect (CATE).

There are a variety of settings---involving particular assumptions on
the regression function and behavioral policy---under which estimates
of $\taustar$ based on $\numobs$ samples are consistent at the
classical $\sqrt{\numobs}$-rate.  Moreover, via the classical notion
of semi-parametric efficiency~\cite{levit1978infinite}, we have a
refined understanding of the optimal instance-dependent constants that
should accompany this $\sqrt{\numobs}$-rate~\cite{hahn1998role}.
However, there are also various settings---of interest in
practice---in which the efficiency bound is infinite, and the
classical $\sqrt{\numobs}$-convergence no longer holds.  This issue is
not only theoretical in nature: when applied to problems of this type,
many standard estimators for $\taustar$---being motivated by classical
considerations---no longer perform well.

At a high level, there are at least two types of phenomena that can
invalidate classical $\sqrt{\numobs}$-consistency.  First, if both the
regression function and behavioral policy need to be estimated from
classes with high complexity, the difficulty of doing so---as opposed
to only the fluctuations intrinsic to the target functional---can
become dominant.  For instance, the
paper~\cite{robins2009semiparametric} studies a variety of such cases
involving H\"{o}lder classes; see also the
paper~\cite{kennedy2022minimax} for related results on CATE
estimation.  Second, the semi-parametric efficiency bound involves
certain moments of the ratio $\tfrac{d \weightfunc}{d \propscore}$.
This so-called \emph{importance ratio} measures how well the
observational data, as controlled by the behavioral policy
$\propscore$, ``covers'' the regions of space relevant for estimating
the functional.  If this coverage is especially bad, then the
importance ratio need not have finite moments, or might even fail to
exist.  This latter cause of breakdown in semi-parametric efficiency
is the primary motivation for the theory and methodology put forth in
this paper.

In the literature on causal inference with observational studies, it is common to impose the
so-called \emph{strict overlap
condition}~\cite{hirano2003efficient,chernozhukov2018double,su2022when}.
The strict overlap condition amounts to imposing a uniform bound on
the importance ratio, and so precludes the possibility of infinite
moments.  Such uniform boundedness conditions also appear frequently
in the closely related literature on bandits and reinforcement
learning.  On one hand, this condition is known to be necessary in a
worst-case sense: as shown by Khan and Tamer~\cite{khan2010irregular},
when neither strict overlap nor structural conditions on the
regression function are imposed, then it is no longer possible to
obtain $\sqrt{\numobs}$-consistency.  It should be noted, however,
that the uniform boundedness condition can be quite stringent.  For
instance, in some recent work, D'Amour et al.~\cite{d2021overlap} show
that it rules out many interesting cases of practical interest,
especially when the model involves high-dimensional covariates.
Motivated by this dilemma, there is a line of past and on-going work
(e.g.,~\cite{chen2004semiparametric,hirshberg2021augmented,ma2020robust})
that proposes estimators that exploit some kind of structure in the
regression function.  Despite this progress, we currently have a
relatively limited understanding of \emph{optimal methods} for
estimating linear functionals based on observational data without
imposing the strict overlap condition.

With this context, the main contributions of this paper are to provide
some insight into the nature of optimal methods for estimating linear
functionals without (strict) overlap.  Our first main result is a
general non-asymptotic lower bound on the mean-squared error of any
estimator.  This lower bound involves a novel variance functional, which depends both on the function class $\Fclass$ used to model the
regression function and the behavior of the importance ratio.  Turning
to upper bounds, we focus on the class of reproducing kernel Hilbert
spaces (RKHSs) as models for the regression function, and provide a
computationally simple procedure that achieves our local minimax lower
bound.  Thus, for RKHS-based models of the regression function, we are
able to identify the instance-dependent and non-asymptotic local
minimax risk up to a constant pre-factor.  As we illustrate by a range
of examples, this mean-squared error can exhibit a range of scalings,
from the classic $\sqrt{\numobs}$-consistency for well-behaved
problems to much slower non-parametric rates in cases where the
importance ratio is badly behaved.

\subsection{An illustrative simulation}
\label{SecIllustrative}

So as to provide intuition for the results to follow, let us consider
a simple family of problems for which the strict overlap assumption is
violated, and compare the performance of the estimator proposed in
this paper to other alternatives.  More specifically, we consider a
missing data problem where the action $\action \in \actionspace =
\{0,1 \}$ is an indicator of ``missingness''.  With the state space
$\statespace = [0,1]$, we construct a weight function $\weightfunc$
and behavioral policy $\propscore$ for which the importance ratio
takes the form
\begin{align}
\label{EqnIllustrative}
\frac{d \weightfunc}{d \propscore}(0 \mid \state) = 1, \quad
\mbox{and} \quad \frac{d \weightfunc}{d \propscore}(1 \mid \state)
\sim (1 - \state)^{-\overexp} \qquad \mbox{where $\overexp \geq 0$ is
  a parameter.}
\end{align}
The parameter $\overexp$ controls the heaviness of the tails exhibited
by the importance ratio: when $\overexp = 0$, the importance ratio is
simply a constant, whereas as $\overexp$ increases, its tails become
increasingly heavy.  Above $\overexp > 1$, it no longer has a finite
second moment, and this transition point turns out to be interesting.

\begin{figure}[h!]
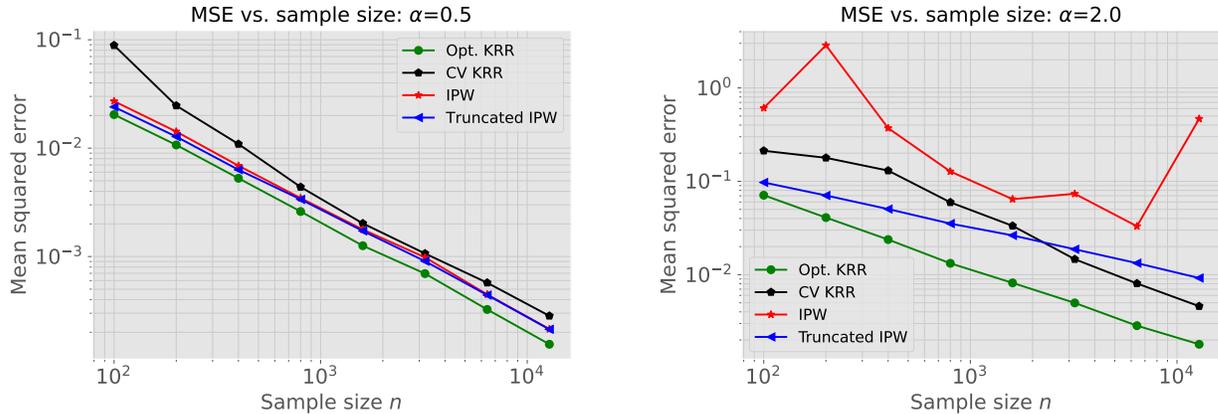

  \begin{center}
    \begin{tabular}{ccc}
      \widgraph{0.5\textwidth}{\figdir/fig_05_ntrial_2000_maxn_128000} &&
      \widgraph{0.5\textwidth}{\figdir/fig_20_ntrial_2000_maxn_128000}
      \\
      (a) Light tails: $\overexp = 0.5$ && (b) Heavy tails: $\overexp = 2.0$
    \end{tabular}
    \caption{Log-log plots of mean-squared error versus sample size
      $\numobs$ for four different estimators fo $\taustar$: procedure
      Opt-KRR is analyzed in this paper, whereas CV-KRR is a related
      method with regularization parameter chosen by cross-validation.
      We also compare to the classical IPW estimate along with a
      truncated version of IPW. (a) Setting $\overexp = 0.5$ yields a
      propensity score with light tails, and our theory predicts 
      classical $\numobs^{-1}$-decay of the MSE for Opt-KRR.  (b) Setting
      $\overexp = 2.0$ yields a heavy-tailed problem, and our theory
      guarantees consistency of Opt-KRR at the rate $\numobs^{-3/4}$.}      
    \label{FigOverexp}    
  \end{center}
\end{figure}

In Figure~\ref{FigOverexp}, we compare the performance of four
different methods: the classical \emph{inverse propensity weighting}
(IPW) estimator~\cite{rosenbaum1987model}, described in more detail
in~\Cref{subsec:simulation}); a truncated version of
IPW~\cite{khan2010irregular}; the optimal kernel-based procedure
proposed in this paper (Opt-KRR); as well as a sub-optimal
kernel-based procedure where the regularization parameter is chosen by
cross validation (CV-KRR).  In each panel and for each estimate
$\tauhat_\numobs$, we plot the mean-squared error
$\Exs[(\tauhat_\numobs - \taustar)^2]$ versus the sample size
$\numobs$ on a log-log scale.  Panel (a) corresponds to the setting
$\overexp = 0.5$: in this case, our theory predicts that the minimax
mean-squared error should decay as $\numobs^{-1}$, as expected for MSE
in the classical regime of $\sqrt{\numobs}$-consistency.  All four
methods are relatively well-behaved for this problem; for our proposed
method (Opt-KRR), performing a linear regression of log-MSE on $\log
\numobs$ gives a slope estimate of $-1.00 \pm 0.01$.  Panel (b), in
contrast, exhibits very different behavior: by setting $\overexp =
2.0$, we obtain a much harder problem.  Here the IPW performance is
very erratic due to the heavy tails of the importance ratio; the
truncated version is better behaved, but still has larger error than
Opt-KRR.  In fact, the theory given in this paper, when specialized to
this family, predicts that for any $\overexp > 1$, the optimal
mean-squared error should decay at the rate $\numobs^{-
  \frac{3}{\overexp + 2 }}$.  Thus, if we set $\overexp = 2$, then we
expect to see an error decay with exponent $-3/4$.  In order to
estimate the decay rate of Opt-KRR, we again perform a linear
regression of the log MSE on $\log(\numobs)$, and obtained an
estimated slope $-0.76 \pm 0.01$.  Once again, we see excellent
agreement with the theory.


\subsection{Our contributions}

We summarize the contributions made in
the remainder of the paper.

\bcar
\item First, working in the setting where the regression function
  belongs to a known function class, we establish a general
  non-asymptotic local minimax lower bound for estimating linear
  functionals from observational data.  The lower bound is defined by
  a variational problem the captures the interplay between the
  geometry of the function class and the importance ratio.  As the
  proof is based on Le Cam's two point approach, portions of the bound
  involve a certain modulus of continuity.

\item Second, specializing to the regression function belonging to a
  ball within a reproducing kernel Hilbert space (RKHS), we analyze a
  class of multi-stage outcome regression estimators. Under certain
  regularity conditions on the RKHS and the conditional covariance
  function, we establish a non-asymptotic upper bound that matches our
  minimax lower bound (up to a constant factor).
\item Third, we illustrate our general result by applying it to a
  range of problems, thereby obtaining a variety of novel minimax
  rates.  For treatment effect estimation when the importance ratio
  diverges at certain points, we show that minimax risk depends on the
  interaction between this singularity and the geometry of the Hilbert
  space.  In the setting of contextual bandits with continuous states
  and actions, we give results on off-policy evaluation of
  deterministic policies, thereby obtaining novel minimax rates that
  are adaptive to the complexity of the state-action space.
\ecar

\medskip

Notably, our multi-stage kernel-based estimator requires \emph{no
knowledge} of the underlying behavioral policy or propensity score
$\propscore$.  This property is very attractive from the
implementation point of view.  At the same time, its performance in
terms of MSE matches our minimax lower bound, which applies to a
broader family of estimators including those that \emph{know} the
behavioral policy.  Thus, we see an interesting implication of our
results: as long the regression function is a member of a RKHS,
knowledge of the behavioral policy plays no role in determining the
minimax risk rate. This is in sharp contrast with H\"{o}lder classes of the
non-Donsker type, where both parts of the model play an important
role~\cite{robins2009semiparametric,kennedy2022minimax}. Finally,
although the value of minimax risk itself depends on the behavioral
policy, the tuning parameter in our kernel-based estimator does not.

\paragraph{Notation:} Throughout this paper, we define
infinite-dimensional vectors (respectively matrices) in a formal sense
as mappings from the positive integers (respectively pairs of positive
integers ) to reals. Given an infinite-dimensional matrix $A$ and a
vector $z$, we define their product pointwise as
\begin{align*}
  [Az]_i \mydefn \sum_{j = 1}^{+ \infty} A_{i, j} z_j, \quad \mbox{for
    each $i = 1, 2, \ldots$,}
\end{align*}
assuming that each summation is absolutely convergent.

Let $\ell_0(\naturals)$ be the set of infinite-dimensional vectors
with finite support---that is, with finitely many non-zero entries.
We say an infinite-dimensional symmetric matrix $A$ is positive
semi-definite, denoted by $A \succeq 0$, when
\begin{align*}
x^\top A x \geq 0,\quad \mbox{for any $x \in \ell_0 (\naturals)$.}
\end{align*}
Note that for any vector space $\mathbb{V}$ in which $\ell_0
(\naturals)$ is dense, if the matrix $A$ maps from $\mathbb{V}$ to
$\mathbb{V}^*$, the definition can be easily extended to ensure that
$x^\top A x \geq 0$ for any $x \in \mathbb{V}$. Given this notation,
we can furthermore define the positive semi-definite ordering $A
\succeq B$ if $A - B \succeq 0$.

Similarly, we can define the inverse of infinite-dimensional
matrix. We call $B = A^{-1}$ if $B \cdot (A x) = A (B x) = x$ for any
$x \in \ell_0 (\naturals)$. Once again, such definition can be easily
extended to larger vector spaces by density arguments, assuming that
both $A$ and $B$ are bounded linear operators acting on suitably
defined spaces.

\subsection{Related work}

Now let us discuss various bodies of related work so as to situate our
work within a broader context.

\paragraph{Instance-optimality for non/semi-parametric estimation:}

For regular parametric models, the classical local asymptotic minimax
(LAM) framework of Le Cam and
Hajek~\cite{lecam1960locally,hajek1972local} specifies the
instance-optimal behavior of estimators as $\numobs \rightarrow
\infty$.  Levit~\cite{levit1978infinite} extended this framework to
semi-parametric settings by considering the collection of all
finite-dimensional sub-models. For the specific class of linear
functional estimation problems considered here,
Hahn~\cite{hahn1998role} laid out the asymptotic lower bounds, whereas
our recent work~\cite{mou2022off} studies the same question within a
non-asymptotic framework.

Beyond the classical $\sqrt{\numobs}$-regime, instance-dependent
optimality for semi-parametric and non-parametric estimation has been
established under various settings. In the literature, exact local
asymptotic minimax risks are obtained for Sobolev space
regression~\cite{pinsker1980optimal,brown1997superefficiency},
spectral density estimation~\cite{hasminskii1986asymptotically}, and
shape-constrained estimation~\cite{han2020limit}.  For estimation
problems involving linear functionals,
Donoho~\cite{donoho1994statistical} establishes information-theoretic
optimality (up to constant factors) of certain class of minimax linear
estimators; in the regression setting, this framework applies to fixed
design problems as opposed to the random design setting of interest
here. Also studying fixed design regression using spline methods, the
unpublished work of Speckman~\cite{speckman1979minimax} is based on a
class of under-smoothed estimators.  These spline-based estimators are
a special case of the more general RKHS set-up considered here for the
random design setting, and we also find that a form of under-smoothing
is optimal.  While all the preceding results are stated as global
minimax risks, due to the location-family structure of the underlying
model and simplifying noise assumptions, the bounds are also
instance-optimal, albeit in a less refined manner.


\paragraph{Overlap and coverage assumptions for off-policy estimation:}

The {\em overlap assumption}, first proposed by Rosenbaum and
Rubin~\cite{rosenbaum1983central}, requires that the behavioral policy
or propensity score takes value within the open interval $(0,
1)$.\footnote{In the classical binary treatment setup, the action
space is $\ActionSpace = \{0, 1\}$, and the propensity score is
defined as $\propscore (x) \mydefn \Prob (\Action = 1 | X = x)$.} In
our general set-up, the overlap assumption is equivalent to requiring
that the importance ratio $\frac{d \weightfunc}{d \propscore}$ exists
everywhere. Such a condition, along with the unconfoundedness
assumption, together imply identifiability of the average treatment
effect~\cite{rosenbaum1983central}, but could lead to arbitrarily slow
rates. In the literature, a popular choice is the much stronger {\em
  strict overlap
  assumption}~\cite{hirano2003efficient,chernozhukov2018double}, which
requires the importance ratio to be uniformly bounded. Khan and
Tamer~\cite{khan2010irregular} shows that strict overlap is a
necessary condition for uniform $\sqrt{\numobs}$-consistency in the
worst case. On the other hand, recent work~\cite{d2021overlap}
revealed that the strict overlap condition can be stringent in some
natural high-dimensional problem setups.  By making stronger
assumptions about the regression function, the strict overlap
condition can be
relaxed~\cite{chen2004semiparametric,hong2020inference}, while still
achieving the semi-parametric efficiency bound in the $\sqrt{\numobs}$
regime. The case when the semi-parametric efficiency bound becomes
infinite, known as the \emph{irregular identification} regime, has
been studied in some past work~\cite{khan2010irregular,ma2020robust},
where truncated versions of IPW estimators are proposed and analyzed
in some special cases.  Moreover, instability in the behavior of
various estimators has been documented in both simulation and
real-data experimental
studies~\cite{lunceford2004stratification,kang2007demystifying,frolich2004finite}.

Uniform boundedness of the importance ratio is also a canonical
assumption in the bandit and reinforcement learning literature.
Focusing on off-policy evaluation for bandit algorithms, Wang et
al.~\cite{wang2017optimal} proposed a ``switch estimator'' that
involves truncating the importance ratio.  Ma et
al.~\cite{ma2022minimax} showed that this procedure is
worst-case optimal for multi-arm bandits.  For off-policy
reinforcement learning problems, uniform bounds on the importance
ratio, known as coverage or concentrability coefficients, appear in
various
papers~\cite{kallus2020double,yin2020asymptotically,xie2021batch}.
Most closely related to our results are the bounds in the
paper~\cite{zanette2021provable}, which apply to MDPs with linear
function approximation and involve a finer-grained measure of the
overlap between the behavioral and target policies.

\paragraph{Kernel and nonparametric methods for off-policy estimation:}

There is also a line of past work on studying various non-parametric
procedures for estimating the average treatment effect under
different structural conditions. Under the strict overlap condition
combined with H\"{o}lder conditions imposed on both the importance
ratio and regression function, minimax rates for ATE estimation have
been established~\cite{robins2009semiparametric,robins2008higher},
albeit with pre-factors depending on the instance that need not be
optimal.  In the Donsker regime considered here, these minimax rates
coincide with the classical $\sqrt{\numobs}$-rate, due to the presence
of strict overlap. Our results reveal different phenomena that can
arise without strict overlap---more specifically, the optimal rate is
determined not only by the complexities of the importance ratio and
the regression function classes, but also by any \emph{singularity} in the
importance ratio, and how it interacts with the functional to be
estimated.  Additionally, our results also apply to estimation of
one-point linear functionals, a generalization of conditional or
heterogeneous average treatment effects. Again with the focus on
H\"{o}lder classes, some recent
work~\cite{gao2020minimax,kennedy2022minimax} has exhibited
rate-optimal non-parametric procedures. In recent years, due to their
flexibility and computational tractability, kernel-based approaches
have been the focus of research in the causal estimation
literature~\cite{singh2020kernela,singh2020kernelb,nie2021quasi},
where kernel-based estimators have been developed for various
functionals.

Recent work on off-policy estimation has explored the use of minimax
linear estimators.  In the fixed-design setup, the
papers~\cite{armstrong2021finite,kallus2018balanced,kallus2020generalized}
apply the classical framework of minimax linear
estimators~\cite{donoho1994statistical,speckman1979minimax} to the
estimation problem for the {\em sample average treatment effect}
(SATE), and establish guarantees of both the asymptotic and
non-asymptotic flavors.  Hirshberg et al.~\cite{hirshberg2019minimax}
studied a minimax linear estimator for the treatment effect when the
regression function belongs to an RKHS; as in the classical
work~\cite{speckman1979minimax}, this estimator can be reformulated in
terms of a standard kernel ridge regression estimate, as can the
two-stage procedure that we analyze in the simpler homoskedastic
setting.  Under the strict overlap condition and some additional
regularity assumptions, they prove asymptotic efficiency as well as
non-asymptotic bounds on the empirical loss function. In a more
general set-up, Hirshberg and Wager~\cite{hirshberg2021augmented}
proposed an augmented minimax linear estimator and established
non-asymptotic normal approximation results. When specialized to
off-policy estimation, their results yield non-asymptotic normal
approximation in the classical regime with finite semi-parametric
efficiency bound, but independent of the strict overlap assumption. An
important contrast with our results is that their bounds involve both
estimation error and approximation error associated with the
importance ratio (via the Riesz representer); in contrast, such terms
do not arise in our approach.


\section{Problem set-up and preview}

We begin in~\Cref{SecSetup} with a precise
formulation of the problem and discussion of some examples.
In~\Cref{SecPreview}, we describe the classical semi-parametric
efficiency bound, and detail how our analysis moves beyond it.


\subsection{Problem set-up and some examples}
\label{SecSetup}

Given some probability distribution $\probxstar$ over the state space
$\StateSpace$, suppose that we observe $\numobs$ i.i.d. triples
$(\State_i, \Action_i, \Outcome_i)$ in which $\State_i \sim
\probxstar$, and
\begin{align}
\Action_i \mid \State_i \sim \wenprop{\cdot}{\State_i}, \quad
\mbox{and} \quad \Exs \big[ \Outcome_i \mid \State_i, \Action_i \big]
= \treateff(\State_i, \Action_i), \qquad \mbox{for $i = 1, 2, \ldots,
  \numobs$.}
\end{align}
In addition to the regression function $\treateff$, our analysis also
involves the conditional variance function
\begin{align}
\label{EqnDefnCondVar}
\sigma^2(\state, \action) \mydefn \Exs \Big[ \abss{\outcome -
    \treateff(\State, \Action)}^2 \mid \State = \state, \Action =
  \action \Big],
\end{align}
which is assumed to exist for any pair $(\state, \action) \in \Xspace
\times \actionspace$.

As previously described, given a collection of signed Radon measures
$\wenweight{\cdot}{\state}$ over the action space $\actionspace$, one
for each $\state \in \Xspace$, our goal is to estimate the value
$\taustar = \LinFunc(\treateff, \probxstar)$ of the bilinear
functional
\begin{align}
  \label{EqnAveFunctional}
(\plainmu, \probx) \mapsto \LinFunc(\plainmu, \probx) & \mydefn
  \int_{\statespace} \int_\actionspace \plainmu(\state, \action) d
  \wenweight{\action}{\state} d \probx(\state)
\end{align}
evaluated at the pair $\plainmu = \treateff$ and $\probx =
\probxstar$.  We require that the signed measures defining $\LinFunc$
satisfy the condition
\begin{align}
\label{eq:signed-radon-msr-bounded}  
\int_\actionspace d \abss{\wenweight{\action}{\state}} \leq 1
\qquad \mbox{for each $\state \in \statespace$.}
\end{align}
This holds automatically when each $\wenweight{\cdot}{\state}$ is a
conditional probability distribution, as in off-policy evaluation for
contextual bandits.

\medskip


\noindent Various types of weight functions $\weightfunc$ arise in practice:
\begin{description}
\item[Average treatment effect (ATE):] This linear functional arises
  with the binary action action space $\actionspace = \{0, 1\}$, and
  weight function $\wenweight{\action}{\state} = \action -
  \frac{1}{2}$ for all $\state$.  With this choice, we have
  \begin{align*}
\taustar & = \frac{1}{2} \Exs_{\State \sim \probxstar} \Big
         [\treateff(\State, 1) - \treateff(\State, 0) \Big],
  \end{align*}
  so that $\taustar$ is proportional to the usual average treatment
  effect (i.e., equal up to the pre-factor $1/2$ that arises from our
  choice of normalization).
\item[Off-policy evaluation for multi-arm contextual bandits:] In the
  multi-arm setting of a contextual bandit, we have a finite action
  space $\actionspace$, and each weight function
  $\wenweight{\cdot}{\state}$ defines a conditional probability over
  the action space, which can be interpreted as a stochastic policy.
  We say that the weight functions $\weightfunc$ define the
  \emph{target policy} whereas the conditional distributions
  $\propscore$ define the \emph{behavioral policy}.
\item[Contextual bandits with continuous arms:] In this case, we take
  the action space $\actionspace$ to be a compact subset of
  $\real^d$. For a deterministic target policy $\pitarget: \Xspace
  \rightarrow \actionspace$, we let $\wenweight{\cdot}{\state}$ be the
  unit atomic mass at $\pitarget(\state)$.
\end{description}

\medskip

Also of interest---in addition to the
functional~\eqref{EqnAveFunctional}---is the variant obtained by
replacing the expectation over $\State \sim \probxstar$ with
evaluation at a known state $\xzero$---namely
\begin{align}
\label{EqnLinFuncXzero}  
\LinFunc(\treateff, \delta_{\xzero}) \mydefn \int_{\actionspace}
\treateff(\xzero, \action) d \wenweight{\action}{\xzero},
\end{align}
where $\delta_\xzero$ can be thought of as a point mass at $\xzero$.
While the functional is determined by $\delta_\xzero$, the samples
$\State_i$ themselves are still drawn from the distribution $\probxstar$
over the state space.

Particular examples of the functional~\eqref{EqnLinFuncXzero} include
the conditional average treatment effect (CATE) in the causal
inference literature, whereas in off-policy reinforcement learning, it
includes the problem of evaluating the policy at a fixed state
$\xzero$ based on off-policy observations.


\subsection{Moving beyond classical semi-parametric efficiency}
\label{SecPreview}

In this section, we explain how this paper moves beyond classical
semi-parametric efficiency.

\paragraph{Recap of classical results:}

We begin by explaining the usual semi-parametric efficiency bound,
which is meaningful when the importance ratio $\frac{d \weightfunc}{d
  \propscore}$ exists and has suitably controlled moments. Under these
conditions, it is possible to obtain estimates $\tauhat_\numobs$ of
$\taustar$ that converge at a $\sqrt{\numobs}$-rate.  We can thus ask
about the variance associated with the rescaled error $\sqrt{\numobs}
(\tauhat_\numobs - \taustar)$, and in particular the smallest one that
can be achieved. For estimating $\taustar = \LinFunc(\treateff,
\probxstar)$, it is known~\cite{hahn1998role} that the smallest
variance achievable, in the sense of semi-parametric efficiency, is
given by
\begin{align}
\label{eq:semi-parametric-efficiency}
\vsemi^2(\treateff, \smallimpratio) & = \myundermath{\var_{\State \sim
    \probxstar} \Big( \int_\actionspace \treateff(\State, \action) d
  \wenweight{\action}{\State} \Big)}{\vprobx^2(\treateff)} +
\myundermath{\Exs_{ \Joint{\probxstar}{\propscore}} \Big \{ \big[\tfrac{d
      \weightfunc }{d \propscore}(\Action \mid \State) \big]^2 \;
  \sigma^2(\State, \Action) \Big \}}{\vclassical^2(\smallimpratio)},
\end{align}
where $\Exs_{\Joint{\probxstar}{\propscore}}$ denotes expectation over
a pair $\State \sim \probxstar$ and $\Action \sim
\wenprop{\cdot}{\State}$.

This optimal variance consists of two term.  The first term
$\vprobx^2(\treateff)$ captures the fluctuations in an estimate of
$\taustar$ due to the randomness in sampling the states from
$\probxstar$.  This term depends on the regression function
$\treateff$, but not on the conditional variance function $\sigma^2$.
In contrast, the second term $\vclassical^2(\impratio)$ depends on
both the importance ratio $\smallimpratio$ and the conditional
variance function $\sigma^2$ but \emph{not} on the regression
function: it captures the interaction between the noise and the
importance ratio $\impratio$.  It is this latter term that can diverge
if the importance ratio is ill-behaved, and accordingly, it is the
term that takes a more refined form in our analysis. \\

\paragraph{Non-asymptotic bounds:}
With this context, our main contributions are to move beyond classical
(asymptotic) semi-parametric efficiency in the following
ways:
\bcar
\item We use Le Cam's method to prove a general non-asymptotic minimax
  lower bound on estimating functionals from observational data
  without the overlap condition, but with $\treateff$ belonging to a
  convex function class $\Fclass$.
\item When $\Fclass$ is a reproducing kernel Hilbert space (RKHS), we
  show that this lower bound can be achieved by a four-stage kernel
  regression procedure, and we compute an explicit representation of
  the minimax risk (sharp up to constant pre-factors).
  \ecar

\medskip
  
Let us describe our explicit representation of the non-asymptotic
minimax risk in the RKHS setting.  Consider an RKHS $\Hil$ that is a
subset of $\Ltwospace(\probxstar \cdot \propscore)$, and suppose that the
regression function $\treateff$ belongs to the Hilbert ball
$\Ball_\Hil(\Rbound)$ of radius $\Rbound$ in this space.  We show that
the non-asymptotic minimax risk replaces the second term
$\vclassical^2(\impratio)$ in the classical semi-parametric efficiency
bound~\eqref{eq:semi-parametric-efficiency} with a novel quantity
associated with the eigenvalues and eigenfunctions associated with the
RKHS.  More precisely, any RKHS of the Mercer type is associated with
a sequence $\{\eigval_j\}_{j =1}^\infty$ of positive eigenvalues, and
associated eigenfunctions $\{\eigfun_j\}_{j =1}^\infty$.  We let
$\EigenMat = \diag \{\eigval_j \}_{j=1}^\infty$ be a diagonal matrix
defined by the eigenvalues, and using the eigenfunctions, we define an
infinite-dimensional vector $\featweighbar =
\featweighbar(\probxstar)$ with elements
\begin{subequations}
  \begin{align}
\label{EqnFirstFeature}    
(\featweighbar)_j & \mydefn \Exs_{\State \sim \probxstar} \Big[
  \int_\actionspace \eigfun_j(\State, \action) d
  \wenweight{\action}{\State} \Big] \qquad \mbox{for $j = 1, 2,
  \ldots$,}
\end{align}
along with the infinite-dimensional matrix $\bigcovsig$ with
elements
\begin{align}
\label{EqnFirstBigCovSig}  
[\bigcovsig]_{jk} \mydefn \Exs_{\Joint{\probxstar}{\propscore}} \Big[
  \tfrac{1}{\sigma^2(\State, \Action)} \eigfun_j(\State, \Action)
  \eigfun_k(\State, \Action) \Big] \qquad \mbox{for $j, k = 1, 2,
  \ldots$.}
\end{align}
\end{subequations}

We prove that when the regression function $\treateff$ lies within a
ball of radius $R$ within this RKHS, then the minimax mean-squared
error for estimating $\taustar$ is proportional to
\mbox{$\frac{1}{\numobs} \big \{ \vprobx^2(\treateff) + \hackhilsq
  \big \}$,} where
\begin{align}
\label{EqnFirstKernel}
\hackhilsq & \defn \featweighbar^\top \Big( \bigcovsig +
\tfrac{1}{\Rbound^2 \numobs} \EigenMat^{-1} \Big)^{-1} \featweighbar.
\end{align}
Note that $\hackhilsq$ depends (among other quantities) on the
sample size $\numobs$, and it can actually diverge as $\numobs
\rightarrow \infty$.  This type of divergence leads to non-parametric
rates for estimating the functional $\taustar$.  Indeed panel (b) in
Figure~\ref{FigOverexp} provides an illustration of this phenomenon in
one particular setting.

\paragraph{Connection to classical semi-parametric efficiency:}

To understand the connection between our result and the the classical
semi-parametric efficiency
bound~\eqref{eq:semi-parametric-efficiency}, let us
consider\footnote{Note that our theory does not require these assumptions,
but imposing them makes clear the connection to classical
semi-parametric efficiency.}  the following special case:
\bcar
\item The importance ratio $\impratio$ exists, and the classical
  semi-parametric efficiency bound is finite.
\item The problem is homoskedastic, with constant conditional variance
  function $\sigma^2(\state, \action) = \sigbar^2$ for all pairs
  $(\state, \action)$.
\ecar Under homoskedasticy, the matrix $\bigcovsig$ is diagonal with
$\frac{1}{\sigbar^2}$ along its diagonal, using the fact that the
eigenfunctions are orthonormal in $\Ltwospace(\probxstar \cdot \propscore)$.
Since the matrix $\EigenMat^{-1}$ is also diagonal, we find that
\begin{subequations}
\begin{align}
\label{EqnSpecialCase}
\hackhilsq & = \sum_{j=1}^\infty
\frac{\featweighbar_j^2}{\frac{1}{\sigbar^2} + \frac{1}{\Rbound^2
    \numobs \eigval_j}} \; \leq \; \sigbar^2 \sum_{j=1}^\infty
\featweighbar_j^2.
\end{align}
When the importance ratio $\impratio$ exists, we can write
\begin{align}
[\featweighbar]_j & \mydefn \Exs_{\State \sim \probxstar} \Big[
  \int_\actionspace \eigfun_j(\State, \action) d
  \wenweight{\action}{\State} \Big] \; = \;
\Exs_{\Joint{\probxstar}{\propscore}} \big[ \eigfun_j(\State, \Action)
  \smallimpratio(\Action \mid \State) \big],
\end{align}
so that $\featweighbar_j$ is the basis coefficient of $\impratio$ when
expanded in the eigenbasis $\{\eigfun_j\}_{j \geq 1}$.  Thus, by
Parseval's theorem, we see that equation~\eqref{EqnSpecialCase}
implies that
\begin{align}
\hackhilsq & \leq \Exs_{\Joint{\probxstar}{\propscore}} \Big[
  \big(\impratio (\Action \mid \State)\big)^2 \sigbar^2 \Big] \; = \;
\vclassical^2(\smallimpratio),
\end{align}
\end{subequations}
so that the Hilbert-restricted functional is always upper bounded by
the classical semi-parametric quantity
$\vclassical^2(\smallimpratio)$.  In fact, when the semi-parametric
efficiency bound is finite and the RKHS is suitably rich---that is,
``universal''---then $\hackhilsq$ converges to
$\vclassical^2(\smallimpratio)$ as $\numobs$ tends to infinity.  All
of these facts hold more generally for heteroskedastic noise, as we
detail in~\Cref{PropExplicit,PropUniversal} to follow in~\Cref{SecExplicit}.


\section{Main results and their consequences}
\label{SecMain}

We now turn to precise statements of our main results, along with
discussion of their consequences for various examples.
In~\Cref{SecGeneral}, we state and prove non-asymptotic lower bounds
that hold for any convex and symmetric function class $\Fclass$
used to model the regression function.  We specialize these lower
bounds to reproducing kernel Hilbert spaces in
\Cref{SecExplicit}, where we derive the
functional~\eqref{EqnFirstKernel} discussed in the previous section.

In~\Cref{SecUpper}, we turn to the complementary question of deriving
upper bounds for reproducing kernel Hilbert spaces.  We begin with the
simpler homoskedastic case in~\Cref{SecUpperHomo} before turning to
the more challenging heteroskedastic case in~\Cref{SecUpperHetero}.
Finally, \Cref{SecExamples} is devoted to the consequences of these
results for various specific examples.

\subsection{Non-asymptotic lower bounds}
\label{SecLower}

Suppose that we model the regression function $\treateff$ using a
class $\funcClass$ of real-valued functions defined on the
state-action space $\statespace \times \actionspace$.  In this
section, we prove some non-asymptotic minimax lower bounds for both
the averaged quantity $\taustar = \LinFunc(\treateff, \probxstar)$ and
the one-point quantities $\taustarzero = \LinFunc(\treateff,
\delta_\xzero)$.  In order to cover both cases in a unified way, for
any distribution $\newxdist$ over the state space $\statespace$, let
us define\footnote{To explain our notational choices, in the special
case that $\newxdist = \probxstar$ and $\Fclass =
\Ball_\Hil(\Rbound)$, this functional is proportional to the quantity
$\hackhilsq$ that we defined previously, as shown in the sequel
(cf.~\Cref{PropExplicit} in~\Cref{SecExplicit}).}
\begin{align}
\label{EqnVsuper}  
\vsuper{\newxdist} & \defn \sqrt{\numobs} \sup_{\optfun \in
  \funcClass} \Big\{ \LinFunc(\optfun, \newxdist) \mid \;
\Exs_{\Joint{\probxstar}{\propscore}} \big[ \tfrac{\optfun^2(\State,
    \Action)}{\sigma^2(\State, \Action)} \big] \leq \tfrac{1}{4
  \numobs} \Big \}.
\end{align}

\subsubsection{General lower bounds}
\label{SecGeneral}
Local minimax bounds describe the behavior of optimal estimators in a local
neighborhood of a given instance.  For the problem at hand, we define
a given problem instance via the pair $\probInstanceStar = (\treateff,
\probxstar)$.  The behavioral policy $\propscore$, conditional
variance function $\sigma^2$, and the weight function $\weightfunc$
are shared across all instances.  Our local neighborhood of a given
instance $\probInstanceStar$ is given by
\begin{align}
  \Nplain{\numobs}(\treateff, \probxstar) & \mydefn \Big\{
  \mbox{$\probx$ s.t.  $\chisqdiv{\probx}{\probxstar} \leq
    \tfrac{1}{\numobs}$, and $\plainmu \in \funcClass$ s.t.
    $\vecnorm{\plainmu -
      \treateff}{\Ltwospace(\Joint{\probxstar}{\propscore})}^2 \leq
    \tfrac{\varbound^2}{\numobs}$} \Big\},
\end{align}
and it defines the local minimax risk
\begin{subequations}
\begin{align}
\label{EqnLocalMinimax}  
\minimaxRisk(\probInstanceStar; \funcClass) & \mydefn
\inf_{\tauhat_\numobs} \sup_{(\plainmu, \probx) \in
  \Nplain{\numobs}(\treateff, \probxstar)} \Exs \big[
  \abss{\tauhat_\numobs - \LinFunc( \plainmu, \probx)}^2 \big], \quad
\mbox{and} \\
\label{EqnLocalMinimaxZero}  
\minimaxRisk(\probInstanceStar, \xzero; \funcClass) & \mydefn
\inf_{\tauhat_\numobs} \sup_{\plainmu \in \Nplain{\numobs}(\treateff,
  \probxstar)} \Exs \big[ \abss{\tauhat_\numobs - \LinFunc(\plainmu,
    \delta_\xzero)}^2 \big].
\end{align}
\end{subequations}
With a slight abuse of notation, in the
definition~\eqref{EqnLocalMinimaxZero}, we have written $\plainmu \in
\Nplain{\numobs}(\treateff, \probxstar)$ to mean $(\plainmu,
\probxstar) \in \Nplain{\numobs}(\treateff, \probxstar)$. For the rest of this paper, we will drop $\probInstanceStar$ in the notation when it is clear from the context.

In stating lower bounds for estimating $\LinFunc(\treateff,
\probxstar)$, we require that the effective noise in the
observations---namely \mbox{$Z(\State) \defn \big(
  \int_{\actionspace}\treateff(\State, \action) d
  \wenweight{\action}{\State} - \taustar \big)$} for $X \sim \probxstar$---has a bounded
kurtosis:
\begin{align}
\label{EqnLtwoLfour}  
  \ltwolfour{Z} \defn \frac{\sqrt{\Exs[Z^4]}}{\Exs[Z^2]} \leq
  \ctwofour(\probxstar) < \infty.
\end{align}
With this condition in place, we are ready to state a lower bound.
\begin{theorem}
\label{ThmLower}
There exists a universal constant $c$ such that for any problem
instance \mbox{$\probInstanceStar = (\treateff, \probxstar)$} with
$\treateff \in \tfrac{1}{2} \funcClass$:
\begin{enumerate}
\item[(a)] Under the moment condition~\eqref{EqnLtwoLfour} and given a
  sample size $\numobs \geq 16 \ctwofour^2(\probxstar)$, the local
  minimax risk~\eqref{EqnLocalMinimax} is lower bounded as
\begin{subequations}  
\begin{align}
\label{eq:local-minimax-statement-general}    
\minimaxRisk(\probInstanceStar, \funcClass) & \geq \frac{c}{\numobs}
\Big\{ \vprobx^2(\treateff) + \vsupersq{\probxstar} \Big\}.
\end{align}
\item[(b)] Given a sample size $\numobs \geq 16$, the local minimax
  risk~\eqref{EqnLocalMinimaxZero} is lower bounded as
  \begin{align}
    \label{EqnLocalMinimaxZeroBound}    
    \minimaxRisk(\probInstanceStar, \xzero; \funcClass) & \geq
    \frac{c}{\numobs} \: \vsupersq{\delta_\xzero}.
\end{align}
\end{subequations}
\end{enumerate}
\end{theorem}
\noindent See~\Cref{SecProofThmLower} for the proof. \\

\medskip

The lower bound~\eqref{eq:local-minimax-statement-general} consists of
two terms.  The first term $\vprobx^2(\treateff)$ captures uncertainty
induced by not knowing the distribution $\probxstar$; in our lower
bound, we obtain it by applying the Le Cam argument to allowable
perturbations of $\probxstar$.  The second term captures the effective
noise induced by a combination of the additive noise, and potential
lack of coverage of the behavioral policy $\propscore$.

The reader should observe
the contrast between the
bound~\eqref{eq:local-minimax-statement-general}, applicable to a
\mbox{$\probxstar$-averaged} functional, and the
bound~\eqref{EqnLocalMinimaxZeroBound} that applies to a one-point
functional. The latter bound takes a similar form, except that the
term $\vprobx^2(\treateff)$ no longer appears.  Here knowledge of
$\probxstar$ is irrelevant, because the functional to be estimated is
known, and does not depend on it.

It is worth emphasizing that~\Cref{ThmLower} and its corollaries are
all stated for a \emph{fixed} behavior policy $\propscore$.
Accordingly, the stated lower bounds apply even to ``oracle''
estimators that know the behaviorial policy.  In practice, this
function often not known, especially for observational studies in
causal inference.  However, as we show in the following section, when
we specialize to reproducing kernel Hilbert spaces
(cf.~\Cref{ThmHomo,ThmHetero} to follow), this lower bound can
achieved (up to universal constants) via a simple procedure that operates without any knowledge of
the policy $\propscore$. Thus, a surprising consequence of our theory
is that, at least for the RKHS case, knowledge of the behavioral
policy $\propscore$ has no effect on the minimax risk.  This statement
is \emph{not} true in general, as demonstrated by past work on
H\"{o}lder classes~\cite{robins2009semiparametric,kennedy2022minimax}.

We also note that the lower bounds in~\Cref{ThmLower} are related to
past work for estimating linear functionals in fixed design regression
(e.g.,~\cite{speckman1979minimax,stander1995minimax,
  donoho1994statistical}).  As in this work, the
quantity~\eqref{EqnVsuper} can be seen as a modulus of continuity for
the functional $\optfun \mapsto \LinFunc(\optfun, \newxdist)$.  Our
work deals instead with a random design setting, so that the proof
techniques are different.  Moreover, it is not always possible to
achieve the lower bounds in~\Cref{ThmLower}; in particular, as we
noted above, for certain types of H\"{o}lder classes and unknown
behavioral policies, sharp lower bounds require an argument that
involves mixtures (as opposed to the two-point Le Cam argument that
underlies~\Cref{ThmLower}).


\subsubsection{Explicit representation for reproducing kernels}
\label{SecExplicit}
  
As noted in~\Cref{SecPreview}, when $\Fclass$ is a reproducing kernel
Hilbert space (RKHS), our minimax lower bounds are sharp (up to a
constant pre-factor), and the optimal risk has an explicit expression. To set up the problem, we consider functions belonging to
a subset of $\Ltwospace(\pstar)$ where $d \pstar(\state, \action) = d
\probxstar(\state) d \wenprop{\action}{\state}$ is a distribution over
$\Uspace \defn \statespace \times \actionspace$. In particular, let
$\kernelfunc$ be a real-valued kernel function defined on the
Cartesian product space $\Uspace \times \Uspace$.  We assume that the
kernel function is continuous and positive semi-definite, and we let
$\rkhs$ be the associated reproducing kernel Hilbert space (RKHS).
Associated with the kernel function is the kernel integral operator
\begin{align*}
\myoptfun \mapsto \kernelfunc(\myoptfun)(z) \defn \int_{\Uspace}
\kernelfunc(\uvar, \uvar') \myoptfun(z') d \pstar(\uvar')
\end{align*}
By Mercer's theorem~\cite{mercer1909functions}, under mild regularity conditions, this operator has
real eigenvalues $\{\eigen_j\}_{j=1}^\infty$, all of which are
non-negative due to the assumption of positive semidefiniteness, along
with eigenfunctions $\{\phi_j\}_{j=1}^\infty$ that are orthonormal in
$\Ltwospace(\pstar)$. Under such notations, the minimax risk can be represented
in terms of these sequences.

Recall the definition~\eqref{EqnVsuper} of the quantity
$\vsuper{\newxdist}$, where $\newxdist = \probxstar$ or $\newxdist =
\delta_{\xzero}$ are the two cases of primary interest in this paper.
For any distribution $\newxdist$ over the state space, we define the
infinite-dimensional vector $\featweighbar(\newxdist)$ with components
\begin{subequations}
\begin{align}
\featweighbar_j(\newxdist) & \mydefn \Exs_{\State \sim \newxdist}
\Big[ \int_\actionspace \phi_j(\State, \action) d
  \wenweight{\action}{\State} \Big].\label{eq:defn-featureweightbar-general}
\end{align}
This vector is a generalization of our previous
definition~\eqref{EqnFirstFeature}, which was specialized to
$\newxdist = \probxstar$.  We also recall from
equation~\eqref{EqnFirstBigCovSig} the infinite-dimensional matrix
$\bigcovsig$ with elements
\begin{align}
[\bigcovsig]_{jk} & \mydefn \Exs_{\Joint{\probxstar}{\propscore}}
\Big[ \tfrac{1}{\sigma^2(\State, \Action)} \phi_j(\State, \Action)
  \phi_k(\State, \Action) \Big],
\end{align}
\end{subequations}
and the diagonal matrix $\EigenMat = \diag
\{\eigval_j\}_{j=1}^\infty$.  With these definitions, we have
\begin{proposition}
\label{PropExplicit}
For the RKHS ball $\Ball_\Hil(\Rbound) \defn \big\{ f \; \mid \;
\rkhsnorm{f} \leq \Rbound \big\}$ and any distribution $\newxdist$
over the state space, we have
\begin{align}
\label{EqnExplicit}
  \frac{1}{2} \, \vdoublesq{\newxdist}{\Ball_\Hil(\Rbound)}
  \stackrel{(a)}{\leq} \featweighbar^T(\newxdist) \Big( \bigcovsig +
  \tfrac{1}{\Rbound^2 \numobs} \EigenMat^{-1} \Big)^{-1}
  \featweighbar(\newxdist) \stackrel{(b)}{\leq} 4 \;
  \vdoublesq{\newxdist}{\Ball_\Hil(\Rbound)}.
\end{align}
\end{proposition}
\noindent See~\Cref{SecProofPropExplicit} for the proof.

\medskip

As discussed in~\Cref{SecPreview}, the functional is closely related
to the classical semi-parametric efficiency bound.  The following result
makes this connection precise:
\begin{proposition}
\label{PropUniversal}
Under the setup of Proposition~\ref{PropExplicit}, if the RKHS $\rkhs$ is dense in
$\Ltwospace(\Joint{\probxstar}{\propscore})$ and $\vsemi(\treateff, \smallimpratio) <
+ \infty$, then we have
\begin{align*}
\lim_{\numobs \rightarrow \infty} \vwenlonghil = \frac{1}{2}
\sqrt{\Exs_{\Joint{\probxstar}{\propscore}} \Big[ \big( \frac{d
      \weightfunc }{d \propscore} (\Action \mid \State) \big)^2 \cdot
    \sigma^2 (\State, \Action) \Big]}.
    \end{align*}
\end{proposition}
\noindent See~\Cref{SecProofPropUniversal} for the proof.


\subsection{Achieving the lower bounds for kernel classes}
\label{SecUpper}

We now show how the lower bounds in~\Cref{ThmLower} can be achieved
when the regression function $\treateff$ is assumed to lie within some
reproducing kernel Hilbert space (RKHS). The setup for RKHS can be found in~\Cref{SecExplicit}.
Our theory involves these eigenvalues and eigenfunctions via the
following notion of \emph{effective dimension}:
\begin{align}
\label{eq:defn-effdim} 
\effdim(\regu) & \mydefn \sup_{(\state, \action) \in \statespace
  \times \actionspace} \sum_{j = 1}^{\infty} \frac{\eigen_j \phi_j^2
  (\state, \action)}{\eigen_j + \regu} \qquad \mbox{for any scalar
  $\regu > 0$}.
\end{align}
Similar notions of effective dimension have been used in past
work~\cite{zhang2002effective,caponnetto2007optimal}. Roughly
speaking, the quantity $\effdim (\regu)$ provides a characterization
of the global complexity of the RKHS at the scale $\regu > 0$.


\subsubsection{Homoskedastic case}
\label{SecUpperHomo}

Let us warm up by describing a simpler (but possibly sub-optimal)
bound that ignores any possible heteroskedasticity.  More
specifically, we suppose that the conditional variance function is
uniformly bounded as $\sigma^2(\state, \action) \leq \varbound^2$ for
all pairs $(\state, \action)$, and prove results in terms of
$\varbound$.

\medskip 
\noindent In this case, the procedure is very simple to describe, and
consists of two steps.
\paragraph{Two-stage procedure:}
Given a data set of size $2 \numobs$, we split it evenly into two sets
$(\dsone{\State}_i, \dsone{\Action}_i, \dsone{\outcome}_i)_{i \in
  [\numobs]}$ and $ (\dstwo{\State}_i, \dstwo{\Action}_i,
\dstwo{\outcome}_i)_{i \in [\numobs]}$, each of size $\numobs$.  Each
step in our procedure uses one of the data splits.
\ALGSTEPS{
\begin{subequations}
\label{eq:two-stage-homoskedastic}  
\paragraph{Stage I:} Given a regularization parameter $\regu_\numobs > 0$,
compute the kernel ridge regression (KRR) estimate on data split I:
\begin{align}
\label{eq:two-stage-1}  
\muhat_{\numobs} \mydefn \arg\min_{ \myoptfun \in \rkhs} \Big\{
\frac{1}{\numobs} \sum_{i = 1}^\numobs \big( \dsone{\outcome}_i -
\myoptfun(\dsone{\State}_i, \dsone{\Action}_i) \big)^2 +
\reguhomo \rkhsnorm{\myoptfun}^2 \Big\}.
\end{align}
\paragraph{Stage II:}  Use the estimate $\muhat_\numobs$ and
split II to compute the empirical average
\begin{align}
\label{eq:two-stage-2}  
\tauhat_\numobs = \frac{1}{\numobs} \sum_{i = 1}^\numobs \Big \{
\int_{\actionspace} \muhat_{\numobs} (\dstwo{\State}_i, \action) d
\wenweight{\action}{\dstwo{\State}_i} \Big \}
\end{align}
\end{subequations}
}

In practice, so as to make most efficient use of the data, one could
also perform a form of cross-fitting
(e.g.,~\cite{chernozhukov2018double}).  However, given that our main
goal is to show that the lower bounds from~\Cref{ThmLower} are
achieved up to constant factors, it suffices to focus attention on the
simpler procedure given here.

\medskip

\paragraph{Assumptions:} In
our analysis of this method, we assume that the kernel function
$\kernelfunc$ is $\kappa$-uniformly bounded:
\begin{align}
  \label{EqnKerBou}
\tag{\mbox{{{Kbou$(\kappa)$}}}}    
\sup_{(\state, \action) \in \Xspace \times \actionspace} \kernelfunc
\big( (\state, \action), (\state, \action) \big) \leq \kappa^2.
\end{align}
This condition is frequently used in the literature on kernel methods.
It is satisfied, for instance, for any continuous kernel function
$\kernelfunc$ on a compact domain $\Xspace \times \actionspace$.

In addition, we assume that the zero-mean noise variables
$\noisevar(\state, \action) \defn \outcome - \treateff(\state,
\action)$ are uniformly $\subgaussian$-sub-Gaussian, meaning that
for all pairs $(\state, \action)$, we have
\begin{align}
  \label{EqnCleanSubGauss}
  \tag{\mbox{{{subG$(\subgaussian)$}}}} \Exs \big[e^{t
      \noisevar(\state, \action)} \big] & \leq e^{\frac{t^2
      \subgaussian^2}{2}} \qquad \mbox{for all $t \in \real$.}
\end{align}

\medskip

We are now equipped to state our first main upper bound.  It requires
that the ratio of sample size and effective dimension at scale
$\regu_\numobs$ is lower bounded as
\begin{subequations}
\begin{align}
\label{eq:sample-size-req-based-on-eigendecay}  
\frac{\numobs}{\effdim(\regu_\numobs)} & \geq c \frac{
  \subgaussian^2}{\varbound^2} \log \big(\frac{\numobs \Rbound
  \kappa}{\varbound \delta} \big) \cdot \log^2(\numobs) \qquad
\mbox{where $\regu_\numobs = \frac{\varbound^2}{\Rbound^2 \numobs}$,}
\end{align}
for a universal constant $c > 0$. We discuss this condition at more length following the statement.

Our result involves the \emph{higher-order term}
\begin{align}
\label{eq:high-order-term}  
\highorder_\numobs(\delta) \mydefn c \frac{\log(1/\delta)}{\numobs}
\Big \{ \kappa \Rbound + \subgaussian \sqrt{\effdim(\regu_\numobs)}
\log(\numobs) \Big \},
\end{align}
where $\delta \in (0,1)$ is a user-specified failure probability. It
can be verified that under the sample size
condition~\eqref{eq:sample-size-req-based-on-eigendecay}, we have
$\highorder_\numobs(\delta) = o(\numobs^{-1/2})$, so this term is of
higher order in the analysis.  Finally, the dominant term in our upper
bound is the quantity
\begin{align}
\vwenlonghilbar^2 \defn \sigbar^2 \sum_{j=1}^\infty \frac{\eigval_j
  \featweighbar_j^2}{\eigval_j + \frac{\sigbar^2}{\Rbound^2 \numobs}}.
\end{align}
\end{subequations}
\begin{theorem}
  \label{ThmHomo}
  Under the~\eqref{EqnKerBou} and~\eqref{EqnCleanSubGauss} conditions,
  suppose that $\treateff \in \Ball_\Hil(\Rbound)$, there exists a universal constant $c > 0$, such that for any
  $\delta \in (0, 1)$ and sample size $2\numobs$ satisfying the
  bound~\eqref{eq:sample-size-req-based-on-eigendecay}.  Then the
  two-stage estimate $\tauhat_\numobs$ computed with regularization
  \mbox{$\regu_\numobs = \frac{\varbound^2}{\Rbound^2 \numobs}$}
  satisfies
\begin{align}
\label{eq:main-homosked}  
\abss{\tauhat_\numobs - \taustar} \leq c \Big \{ \vprobx (\treateff)
+ \vwenlonghilbar \Big \} \sqrt{\tfrac{\log (1 / \delta)}{\numobs} } +
\highorder_\numobs(\delta),
\end{align}
with probability at least $1 - \delta$.
\end{theorem}
\noindent See \Cref{SecProofThmHomo} for the proof.

\medskip

\noindent Let us make a few comments about this result, and its
connection to our lower bounds.

\paragraph{Comparison with~\Cref{ThmLower}:}  As noted, 
the dominant term in the bound~\eqref{eq:main-homosked} is the first
one.  When the noise is homoskedastic (i.e., constant conditional
variance), then this first term matches the lower bound given
in~\Cref{ThmLower} up to constants and the logarithmic
factor\footnote{While we have stated a high probability guarantee, a
simple modification yields an estimator with mean-squared error
guarantees. In particular, since $|\taustar| \leq \sup_{\state,
  \action} \abss{\treateff(\state, \action)} \leq \Rbound \sqrt{\kappa} $ by
the Cauchy--Schwarz inequality, we can construct a truncated estimator
\begin{align*}
\widetilde{\tau}_\numobs \mydefn \mathrm{sgn}
\big(\tauhat_\numobs\big) \cdot \min \big\{ |\tauhat_\numobs|, \Rbound
\sqrt{\kappa} \big\}.
\end{align*}
By construction, we have $|\widetilde{\tau}_\numobs - \taustar| \leq
|\tauhat_\numobs - \taustar|$ almost surely, and since
$\widetilde{\tau}_\numobs$ is a bounded random variable, the
high-probability bounds established in~\Cref{ThmHomo} can be converted
to a MSE bound whose leading term matches \Cref{ThmLower} up to a
constant factor.\label{ftnt:high-prob-to-mse} } in the failure
probability $\delta$.  When the conditional variance function is not
constant, then our bound~\eqref{eq:main-homosked} no longer
matches~\Cref{ThmLower}.  We rectify this shortcoming
in~\Cref{SecUpperHetero}, where we analyze a more refined four-stage
procedure that adapts to heteroskedasticity.

It should be emphasized that the two-stage estimator analyzed
in~\Cref{ThmHomo} does not require any knowledge of the behavioral
policy $\propscore$.  At the same time, as we just described, for
homoskedastic noise, it matches the lower bound
from~\Cref{ThmLower}, which applies even to oracle estimators that
know the policy.  Thus, we conclude that at least in the special case
of an RKHS, knowledge of the behavior policy does not alter minimax
risks (apart from possibly in constant factors).

\paragraph{Tuning parameter:}  The only tuning parameter in the estimator
is the regularization weight $\regu_\numobs =
\frac{\varbound^2}{\Rbound^2 \numobs}$.  This choice depends on the
signal-to-noise-ratio, as measured by the ratio
$\Rbound^2/\varbound^2$, but does not depend on kernel eigenvalues or
other aspects of the problem.  We note that this choice also appears
in the classical work on linear functional estimation in fixed design
settings~\cite{speckman1979minimax,stander1995minimax}, but the
analysis leading to it in our random design setting is quite
different.

The decay rate $\regu_\numobs \asymp \numobs^{-1}$ of the
regularization parameter is much faster than the standard one required
to achieve optimal mean-squared error when estimating the full
regression function (c.f.~\cite{wainwright2019high}, Chapter 13).
Consequently, the first stage of our procedure outputs an
\emph{under-smoothed} estimate of the regression function $\treateff$,
and using this estimate in the second stage produces an optimal
estimate of the functional.  This difference arises because the
bias-variance trade-off that underlies estimating the functional of
$\treateff$ is very different from that associated with estimating the
full regression function $\treateff$.  In particular, when estimating
a functional, we pay for variance only at the direction of the target functional, whereas
the bias induced by regularization wholly appears in the estimation
error.


\paragraph{Lower bound on sample size:}  Finally, let us comment
on the required lower
bound~\eqref{eq:sample-size-req-based-on-eigendecay} on the sample
size.  There are various conditions that ensure~\eqref{eq:sample-size-req-based-on-eigendecay}.  For example, in various examples, it is possible to show
that the effective dimension satisfies the bound
\begin{align}
\label{EqnEffdimBound}  
\effdim(\regu) \leq \frac{D_0}{\regu^{1-\offpar}} \qquad \mbox{for some
    scalar $\offpar \in (0, 1]$.}
\end{align}
In~\Cref{SecExamples}, we discuss various concrete applications in
which this growth condition holds.  Under the
bound~\eqref{EqnEffdimBound}, the sample size
condition~\eqref{eq:sample-size-req-based-on-eigendecay} is satisfied
as long as
\begin{align*}
\frac{\numobs}{ \log^{3 /\offpar} \big(\frac{\numobs \Rbound
    \kappa}{\varbound \delta} \big)} \geq c D_0^{1/ \offpar}
\Rbound^{\frac{2}{\offpar} - 2} \subgaussian^{\frac{2}{\offpar}}
\varbound^{2 - \frac{4}{\offpar}}.
\end{align*}

In~\Cref{AppSupNorm}, we present various conditions under which the
effective dimension satisfies a growth condition that ensures the
sample size condition~\eqref{eq:sample-size-req-based-on-eigendecay}
can be satisfied.  Moreover, in~\Cref{subsubsec:near-optimal-rates},
we present alternative guarantees that do not rely on any additional
growth conditions.


\subsubsection{Extension to heteroskedasticity}
\label{SecUpperHetero}

We now turn to the more challenging problem of achieving the minimax
optimal risk in the heteroskedastic case. In this case, we propose
and analyze a four-stage procedure. Since the conditional variance
function is non-constant and unknown, we need to estimate it, and the
first two steps of our four-stage procedure are devoted to this task.

Let us provide a high-level perspective.  The first stage generates a
rough estimate $\mutil_\numobs$ of the regression function
$\treateff$.  In the second stage, we first use $\mutil_\numobs$ to
compute estimates \mbox{$Z_i \defn \big \{ \Outcome_i -
  \mutil_\numobs(\State_i, \Action_i) \big\}^2$} of the squared noise
associated with a new set $\Set$ of triples $\{(\Outcome_i, \State_i,
\Action_i)\}_{i \in \Set}$.  We then compute an estimate of the
conditional variance function of the form
\begin{align}
\sighat_{\numobs}^2 \mydefn \abstrEst \Big(\{(\State_i, \Action_i,
Z_i)\}_{i \in \Set} \Big),
\end{align}
for a suitably chosen estimator $\abstrEst$.  We allow the conditional
variance estimator $\abstrEst$ to take different forms depending on
the application, so our set-up provides a family of possible
procedures, indexed by this choice. In our theory, we require only a
relatively mild form of accuracy from the estimator, which we refer to
as {\em robust pointwise accuracy}.

Now let us specify all four stages in more detail.  Given sample size
$4 \numobs$, we split the data evenly into four pieces, and perform
the following four steps:
\ALGSTEPS{
\begin{subequations}
\label{eq:four-stage-framework}
\paragraph{Stage I:}
Using the first dataset $ (\dsone{\State}_i, \dsone{\Action}_i,
\dsone{\outcome}_i)_{i = 1}^\numobs$ and regularization parameter
$\regu_\numobs > 0$, compute the pilot estimate
\begin{align}
\label{EqnPilot}  
\mutil_{\numobs} \mydefn \arg\min_{\plaintreateff \in \rkhs} \Big\{
\frac{1}{\numobs} \sum_{i = 1}^\numobs \big( \dsone{\outcome}_i -
\plaintreateff(\dsone{\State}_i, \dsone{\Action}_i) \big)^2 +
\regustageone \rkhsnorm{\plaintreateff}^2 \Big\}.
\end{align}

\paragraph{Stage II:}
Using the second dataset $ (\dstwo{\State}_i, \dstwo{\Action}_i,
\dstwo{\outcome}_i)_{i = 1}^\numobs$ and the procedure $\abstrEst$,
compute the squared noise estimates $\dstwo{Z}_i \defn
(\dstwo{\outcome}_i - \mutil_\numobs (\dstwo{\State}_i,
\dstwo{\Action}_i))^2$ based on the pilot estimate~\eqref{EqnPilot},
and then compute the estimate
\begin{align}
\sighat_{\numobs}^2 \mydefn \abstrEst \Big( \Big\{ \dstwo{\State}_i,
\dstwo{\Action}_i, \dstwo{Z}_i \Big\}_{i \in [\numobs]} \Big) \qquad
\mbox{of the conditional variance function.}
\end{align}

\paragraph{Stage III:}
Using the third dataset $(\dsthree{\State}_i, \dsthree{\Action}_i,
\dsthree{\outcome}_i)_{i = 1}^\numobs$, regularization parameter
$\regustagethree > 0$, and the estimated function $\sighat_\numobs^2$,
compute the weighted regression estimate
\begin{align}
\label{EqnWeightedRegression}  
\muhat_\numobs \mydefn \arg \min_{\plaintreateff \in \rkhs} \Big\{
\frac{1}{\numobs} \sum_{i = 1}^\numobs \frac{1}{\sighat^2_\numobs
  (\dsthree{\State}_i, \dsthree{\Action}_i)} \big(
\dsthree{\outcome}_i - \plaintreateff(\dsthree{\State}_i,
\dsthree{\Action}_i) \big)^2 + \regustagethree
\rkhsnorm{\plaintreateff}^2 \Big\}.
\end{align}
\paragraph{Stage IV:} 
Using the fourth dataset $(\dsfour{\State}_i, \dsfour{\Action}_i,
\dsfour{\outcome}_i)_{i = 1}^\numobs$ and the weighted regression
estimate~\eqref{EqnWeightedRegression}, compute the empirical average
\begin{align}
\label{eq:four-stage-framework-4}  
\tauhat_\numobs = \frac{1}{\numobs} \sum_{i = 1}^\numobs \int
\muhat_{\numobs} (\dsfour{\State}_i, \action) d
\wenweight{\action}{\dsfour{\State}_i}
\end{align}
\end{subequations}
} 

We remark that the idea of re-weighting with estimated conditional variance has been utilized in literature, in the context of parameter estimation for linear models. See the paper~\cite{romano2017resurrecting} and references therein for detailed discussion. 
We now turn to the analysis of the $4$-stage procedure.  Rather than
analyze a particular estimator $\abstrEst$ of the conditional variance
function, let us lay out an abstract condition that handles a variety
of different estimators

\paragraph{Robust pointwise variance estimators:}

This property is a way of certifying that the estimator $\abstrEst$
provides an $\varepsilon$-accurate estimate in a pointwise sense: for
any fixed pair $(\state_0, \action_0)$, with probability at least $1-\delta$, we have
\begin{align}
\label{eq:conclusion-for-robust-ptwise-risk}  
  \abss{\abstrEst \Big( \big\{ \State_i, \Action_i, Z_i \big\}_{i \in
      [\numobs]} \Big) (\state_0, \actionzero) - \sigma^2 (\state_0,
    \actionzero)} \leq \varepsilon.
\end{align}
The key is to quantify how
errors in the inputs $Z_i$ as approximations of the squared noise
$(\Outcome_i - \treateff(\State_i, \Action_i))^2$ affect this
guarantee.  We do so via a pair of functions on the inputs
$(\varepsilon, \delta)$, known as the tolerance function $\wenbias$
and sample threshold $\Wensamp$ respectively.

\begin{definition}
The procedure $\abstrEst$ is $(\wenbias, \Wensamp)$-pointwise-robust
if for any pair $(\varepsilon, \delta) \in [0,1]^2$, any dataset
$\{Z_i\}_{i=1}^\numobs$ of size $\numobs \geq \Wensamp(\varepsilon,
\delta)$, consisting of variables such that for each $i \in [\numobs]$\footnote{In the
definition~\eqref{eq:observation-assumption-for-robust-ptwise-risk},
the quantity $\subgaussian$ is the sub-Gaussian parameter
(cf. condition~\eqref{EqnCleanSubGauss}), whereas
$\vecnorm{\cdot}{\psi_1}$ is the Orlicz$(1)$-norm, or sub-exponential
parameter.}
\begin{align}
\label{eq:observation-assumption-for-robust-ptwise-risk}
  \abss{\Exs [Z_i \mid \State_i, \Action_i] - \sigma^2(\State_i,
    \Action_i)} \leq \wenbias (\varepsilon, \delta) \quad
  \mbox{and}\quad \vecnorm{Z_i \mid \State_i, \Action_i }{\psi_1} \leq
  4 \big(\subgaussian^2 + \wenbias (\varepsilon, \delta) \big),
\end{align}
then for any fixed pair $(\xzero, \action_0)$, the bound~\eqref{eq:conclusion-for-robust-ptwise-risk} holds with
probability $1 - \delta$.
\end{definition}
\noindent There are various estimators that satisfy the robust
pointwise risk property; see~\Cref{subsubsec:cond-var-est-example} for
further discussion.

In order to analyze the $4$-stage procedure, we require one additional
condition on the conditional variance function: there are scalars $0 <
\lowersigma \leq \varbound < \infty$ such that
\begin{align}
\label{EqnSigmaInt}  
\tag{$\sigma$-INT} \sigma(\state, \action) \in [\lowersigma,
  \varbound] \quad \mbox{for all $(\state, \action)$.}
\end{align}
We also require that the sample size satisfies the lower bounds
\begin{subequations}
  \begin{align}
\label{eq:hetero-sample-size-req-first-stage}     
\frac{\numobs}{\effdim(\regustageone)} & \geq \Big\{ c \frac{
  \subgaussian^2}{\lowersigma^2} \log \big(\frac{\numobs \Rbound
  \kappa}{\lowersigma \delta} \big) \cdot \log^2 \numobs \Big\}, \\
\label{eq:hetero-sample-size-req-second-stage}
\numobs \geq \Wensamp \big( \lowersigma / 2, \delta / (2\numobs) \big)
& \quad \mbox{and} \quad \frac{\numobs}{\subgaussian \effdim
  (\regustageone) \log (\numobs/ \delta)} \geq \frac{1}{\wenbias \big(
  \tfrac{\lowersigma^2}{2}, \tfrac{\delta}{2\numobs} \big)}.
\end{align}
\end{subequations}

With this set-up, we are now ready to state a guarantee on our 4-state
procedure.  Note that it has two tuning parameters: the regularization
parameter $\regustageone$ from the first stage regression, and the
regularization parameter $\regustagethree$ from the weighted
regression in the third stage.  Our guarantee applies to the procedure
using the parameters
\begin{align}
\label{eq:optimal-lambda-four-stage}  
\regustageone = \frac{\varbound^2}{\Rbound^2 \numobs} \quad \mbox{and}
\quad \regustagethree = \frac{1}{\Rbound^2 \numobs}.
\end{align}

\begin{theorem}
\label{ThmHetero}
In addition to the assumptions of~\Cref{ThmHomo}, suppose
that the conditional variance function satisfies the interval
condition~\eqref{EqnSigmaInt}, the estimator $\abstrEst$ is
$(\wenbias, \Wensamp)$-robust, and the sample size $4\numobs$ satisfies
the lower bounds~\eqref{eq:hetero-sample-size-req-first-stage}
and~\eqref{eq:hetero-sample-size-req-second-stage}.  Then using
regularization parameters from
equation~\eqref{eq:optimal-lambda-four-stage}, the $4$-stage procedure
yields an estimate $\tauhat_\numobs$ such that
\begin{align}
\label{eq:main-heterosked}  
\abss{\tauhat_\numobs - \taustar} \leq c \Big\{ \vprobx (\treateff) +
\vwenlonghil \Big\} \sqrt{\frac{\log (1 / \delta)}{\numobs} } +
\highorder_\numobs(\delta),
\end{align}
with probability $1 - \delta$, where the higher-order term
$\highorder_\numobs(\delta)$ was previously
defined~\eqref{eq:high-order-term}.\footnote{We take $\regu_\numobs = \regustageone$ in its expression.}
\end{theorem}
\noindent See~\Cref{SecProofThmHetero} for the proof.

\medskip

A few remarks are in order. First,~\Cref{ThmHetero}
is adaptive to the heteroskedastic nature of the observation noise ---
the term $\vwenlongsq$ involves the actual conditional variance
$\sigma$, instead of its uniform upper bound $\varbound$. With such a
fine-grained variance, the bound~\eqref{eq:main-heterosked} achieves
the instance-dependent optimality result in~\Cref{ThmLower}, up
to universal constants and high-order terms. (See
footnote~\ref{ftnt:high-prob-to-mse} for the connection between
high-probability bounds and mean-squared error bounds.)

We note that the sample size
condition~\eqref{eq:hetero-sample-size-req-first-stage} is slightly
stronger than the
condition~\eqref{eq:sample-size-req-based-on-eigendecay} used in
\Cref{ThmHomo}, with the variance upper bound
$\varbound^2$ in the denominator replaced by the lower bound
$\lowersigma^2$. \Cref{ThmHetero} further requires an
additional sample size
condition~\eqref{eq:hetero-sample-size-req-second-stage}, which comes
from the sample complexity of the robust pointwise estimator
$\sighat_\numobs$.

\paragraph{Extension to estimating one-point functionals:}

Now we extend our results to estimating the one-point functional
$\taustar(\xzero) \defn \LinFunc(\treateff, \delta_{\xzero})$. In this
case, the target functional is known, so that the fourth stage of the
four-stage procedure is not necessary.  It suffices to split the data
into three folds in total, and we plug in the regression function
$\muhat_\numobs$ directly to obtain the estimate
\begin{align}
\label{eq:cate-estimator}  
\tauhat_{\numobs} (\xzero) \mydefn
\LinFunc(\muhat_\numobs, \delta_\xzero) & = \int_{\actionspace}
\muhat_\numobs (\xzero, \action) d \wenweight{\action}{\xzero}.
\end{align}
This estimate satisfies optimal guarantees matching~\Cref{ThmLower}(b)
up to a constant factor. In particular, under the setup
of~\Cref{ThmHetero}, for any $\xzero \in \Xspace$, we have
\begin{align}
\label{EqnCATEclaim}  
\abss{\tauhat_{\numobs} (\xzero) - \taustar(\xzero)} \leq
c \cdot \vsuperhil{\delta_{\xzero}}
\sqrt{\frac{\log (1 / \delta)}{\numobs}},
\end{align}
with probability $1 - \delta$.  See~\Cref{subsubsec:proof-cor-cate}
for the proof.

\medskip

A few remarks are in order. First, the upper bound in
equation~\eqref{EqnCATEclaim} matches the local minimax lower bound
in~\Cref{ThmLower}(b) up to universal constant factors, exhibiting its
optimality in an instance-dependent sense.\footnote{Following the
discussion in footnote~\ref{ftnt:high-prob-to-mse}, the
high-probability bound can be readily converted into a mean-squared
error bound using a simple truncation method.} As opposed
to~\Cref{ThmHomo,ThmHetero}, the optimal instance-dependent risk is
achieved (up to universal constants) without additional high-order
terms. The choice of parameters and sample size requirement
in~\eqref{EqnCATEclaim} is exactly the same as the one
in~\Cref{ThmHetero}, and does not depend on the query point $\xzero$.
Such an adaptive property makes the estimator useful in practice,
allowing for a plug-and-play approach: one only needs to run stages
I--III of the four-stage framework~\eqref{eq:four-stage-framework},
and generate an estimator $\muhat_\numobs$. By substituting such an
estimator in equation~\eqref{eq:four-stage-framework-4} using another
fold of data, or in equation~\eqref{eq:cate-estimator} for any query
point $\xzero$, optimal and adaptive guarantees can always be
achieved.


\subsection{Consequences for some concrete examples}
\label{SecExamples}

In this section, we develop some consequences of our general theory
for some specific classes of problems, including the missing data
problem without overlap (\Cref{SecMissingOverlap}), for which we
presented an illustrative simulation in~\Cref{SecIntro}.

\subsubsection{A missing data example without overlap assumption}
\label{SecMissingOverlap}

Consider the classical missing data setting, where the action space is
$\actionspace = \{0, 1\}$ and the weight is given by
$\wenweight{\action}{\state} = \action$ for any $\state
\in \Xspace$. Assume without loss of generality that $\outcome = 0$
whenever $\Action = 0$. We slightly abuse the notation to use
$\treateff : \Xspace \rightarrow \real$ to denote the outcome function
$\treateff(\cdot, 1)$ and use $\propscore: \Xspace \rightarrow [0, 1]$
to denote the propensity score $\propscore (1 \mid \cdot)$. Similarly, we
use $\kernelfunc(\state, \state')$ to denote $\kernelfunc (\state, 1),
(\state', 1))$, and let the kernel function be $0$ if one of the
arguments has action equal to $0$. Under this simplified notation, the
inner product of $\Ltwospace (\Joint{\probxstar}{\propscore})$ takes the form
\begin{align*}
 \inprod{f_1}{f_2} \mydefn \int_{\Xspace} f_1(\state) f_2(\state)
 \propscore (\state) d \probxstar (\state),
\end{align*}
and we are interested in estimating the average treatment effect and
its conditional analogue
\begin{align*}
    \taustar \mydefn \Exs_\probx \big[ \treateff(\State) \big], \quad
    \mbox{and} \quad \LinFuncZero \mydefn \treateff(\xzero).
\end{align*}
For concreteness, we let the state space be a unit interval $\Xspace =
[0, 1]$ and take the input distribution $\probx$ be the uniform
distribution on $\Xspace$. In order to illustrate the effect of the lack
of the overlap condition on the risk, given a scalar $\overexp > 0$, we
construct the following propensity score function
\begin{align}
\label{EqnParticular}
\propscore(\state) = (1 - \state)^\overexp \quad \mbox{for any $\state
  \in [0, 1]$.}
\end{align}
Our goal is to understand the effect of a singularity in the
importance ratio with local $\overexp$-th order polynomial
growth. The specific location of such singularity, and any properties
apart from the existence of this $\overexp$-th order singularity are
not germane to our comparison, so that we have chosen the particular
form~\eqref{EqnParticular} for technical convenience.

We consider an RKHS $\rkhs$ corresponding to the first-order Sobolev
space on $[0, 1]$ (see~\cite{wainwright2019high}, Chapter 12). Its kernel function is given by
$\kernelfunc(\state, \state) = \min \{ \state, \state' \}$, and the
corresponding RKHS $\rkhs$ consists of functions $f$ satisfying $f(0)
= 0$, and
\begin{align*}
\rkhsnorm{f}^2 \mydefn \int_0^1 \big( f'(\state) \big)^2 d \state < +
\infty.
\end{align*}
We assume that the regression function $\treateff: \Xspace
\rightarrow \real$ belongs to this RKHS, with $\rkhsnorm{\treateff}
\leq 1/2$. Finally, we take the conditional variance as
$\sigma^2(\state, \action) \equiv 1$ for any state-action pair
$(\state, \action)$, and assume that the sub-Gaussian parameter
$\subgaussian$ is of order one.

With this set-up, we are ready to compute minimax rates for various
linear functionals. Throughout this section, we use the notation $a_n \asymp b_n$ to denote that the ratio $a_n / b_n$ satisfies finite positive upper and lower bounds depending on the constants $(\alpha, \xzero)$ but independent of $n$. As mentioned before, we omit the problem instance $\probInstanceStar$ in the notations $\minimaxRisk \big( \probInstanceStar, \funcClass \big)$ and $\minimaxRisk \big( \probInstanceStar, \xzero; \funcClass \big)$ for minimax risk rates.
\begin{corollary}
\label{cor:missing-data-rates}
Under the above set-up, for any function $\treateff \in
\ball_\rkhs(1/2)$, the minimax risk for estimating the linear
functional $\taustar$ is given by
\begin{subequations}    
    \begin{align}
        \minimaxRisk \big( \ball_\rkhs (1) \big) \asymp
        \begin{cases}
        \numobs^{-1} & \alpha < 1,\\ \numobs^{-1} \log \numobs &
        \alpha = 1,\\ \numobs^{- \frac{3}{\alpha + 2}} & \alpha >
        1.\label{eq:missing-data-example-ate-rate}
        \end{cases}
    \end{align}
    For the one-point functional $\LinFuncZero$, we have
    \begin{align}
        \minimaxRisk \big( \xzero; \ball_\rkhs (1) \big)
        \asymp \begin{cases} 0 & \xzero =
          0,\\ \numobs^{-1/2} & \xzero \in (0, 1),\\ \numobs^{\frac{-
              1}{2 + \alpha}} & \xzero = 1,
        \end{cases} \label{eq:missing-data-example-cate-rate}
    \end{align}
\end{subequations}    
\end{corollary}
\noindent See~\Cref{app:subsec-proof-missing-data-rates} for the proof.

\medskip

A few remarks are in order. For the average treatment effect
$\taustar$, the optimal rate of estimation exhibits a phase transition
depending on the local growth exponent $\alpha$. In the regime $\alpha
\in [0,1)$, the importance ratio is sufficiently well-behaved that the
  classical quantity $\vsemi$ is finite, so that we obtain convergence
  at the classical $\sqrt{\numobs}$-rate. Slower rates arise once
  $\alpha \geq 1$, where the variance $\vsemi^2$ is infinite. A large
  value of $\alpha$ yields fewer observations in the neighborhood of
  $x = 1$, which in turn leads to slower rate of convergence. Note
  that even if the target $\taustar$ is defined as a global average
  over the interval $[0, 1]$, the optimal rate of convergence is still
  affected by the singularity within the interval.\footnote{The proofs
  in \Cref{app:subsec-proof-missing-data-rates} can be easily extended
  to propensity score functions with zeros at any finite subset of
  $[0, 1]$, with arbitrary behavior except for the local growth
  conditions around the zeros.} Finally, we note that
  although~\Cref{cor:missing-data-rates} exhibits a wide spectrum of
  rates, they all can be achieved adaptively---that is, using an
  estimator that requires no knowledge of the behavioral policy
  $\propscore$ nor the exponent $\alpha$.

Let us make a few comments on the conditional average treatment effect
$\LinFuncZero$.  In this case, the problem becomes trivial at $\xzero
= 0$, as the functions in the Sobolev space $\rkhs$ satisfy
$\treateff (0) = 0$. The optimal rate is $\numobs^{-1/2}$ for any
$\xzero \in (0, 1)$, which corresponds to the minimax one-point rate
for Sobolev regression in
literature~\cite{tsybakov2008introduction}. A much slower minimax rate
is observed at $\xzero = 1$, where the scarcity of outcome
observations is controlled by the exponent $\alpha$. These rates, just
as in the ATE case, can be achieved using an estimator without any
knowledge of the function $\propscore$.


\subsubsection{Off-policy evaluation with continuous actions}
\label{subsubsec:continuum-bandit}

Now we consider a continuum-arm bandit setup.  For simplicity, we work
with the state space $\Xspace = [0, 1]^{\dimx}$ and the action space
$\actionspace = [0, 1]^{\dima}$, and let the distributions $\probx$,
$\propscore (\cdot \mid x)$ be the uniform distribution on the spaces
$\Xspace$ and $\actionspace$, respectively, for any $x
\in \Xspace$. Given a scalar $\smoothorder > (\dimx + \dima) / 2$, we
let the RKHS $\rkhs = \rkhs^\smoothorder$ be the Sobolev space of order $\smoothorder$,
with periodic boundary conditions (so that the state-actions spaces are seen as tori).
    
Under this setup, the eigenfunctions are given by the standard
(complex) Fourier bases on the torus $\torus^{\dimx + \dima}$, which
can be written in a product form
\begin{align*}
  \Big\{ (\state, \action) \mapsto \phi_j(\state) \psi_k
  (\action) \Big\}_{j, k \geq 0},
    \end{align*}
where $\{\phi_j\}_{j \geq 0}$ and $\{\psi_k\}_{k \geq 0}$ are the
Fourier bases on the tori $\torus^{\dimx}$ and $\torus^{\dima}$,
respectively. Note that these eigenfunctions are uniformly
bounded in sup norm.

Throughout this section, we view the problem parameters $(d_x, d_a,
\smoothorder)$ as universal constants, and suppress any constant
factor depending only on them. For the Sobolev space
$\rkhs^\smoothorder$, let $\eigen_{j, k}$ be the eigenvalue associated
to the eigenfunction indexed by $j, k$, which satisfies the decay
condition (see~\cite{berlinet2011reproducing})
\begin{align}
  \eigen_{j, k} \asymp \min \Big\{ j^{- 2 \smoothorder / \dimx},
  k^{- 2 \smoothorder / \dima}
  \Big\}.\label{eq:sobo-eigendecay-general}
\end{align}
Combining the eigendecay assumption and the boundedness condition on
the eigenfunctions, we can verify that condition~\eqref{EqnKerBou} holds;
in particular, we have
\begin{align*}
\kappa^2 \defn \sup_{\state, \action} \sum_{j, k \geq 0} \eigen_{j, k}
\phi_j^2 (\state) \psi_k^2 (\action) \leq \sup_{\state, \action}
\sum_{j, k \geq 0} \eigen_{j, k} < \infty,
\end{align*}
where the last inequality follows from the fact $\smoothorder > (\dima + \dimx) / 2$.
      
Given a deterministic target policy $\pitarget: \Xspace
\rightarrow \actionspace$, we let $\wenweight{\cdot}{\state}$ be
the atomic measure on $\pitarget(\state)$ for $\state
\in \statespace$, so that the linear functionals of interest take the
following form:
\begin{align*}
  \taustar \mydefn \int_\Xspace \treateff(\state, \pitarget(\state)) d
  \state, \quad \mbox{and} \quad \LinFuncZero \mydefn
  \treateff(\xzero, \pitarget(\xzero)).
\end{align*}
Finally, we let the conditional variance function be unity
$\sigma^2 \equiv 1$, and assume that the sub-Gaussian parameter
$\subgaussian$ is of order one.  Let $\treateff$ be any function
lying in the Hilbert ball $\ball_\rkhs(1/2)$.

Note that in this example, the importance ratio
$\tfrac{d\weightfunc}{d \propscore}$ is not well-defined, as the
measure $\wenweight{\cdot}{\state}$ is atomic, for any $\state
\in \statespace$. Nevertheless, estimation is still possible, and our
general frameworks provide precise characterization of the minimax
risks, stated as follows.
\begin{corollary}
\label{cor:continuum-bandit-rates}
Under the above setup, we have
\begin{subequations}
\begin{align}
  \minimaxRisk \big( \ball_\rkhs (1) \big) &\asymp \frac{\var_\probx
    \big( \treateff(\State, \pitarget(\State)) \big)}{\numobs} +
  \sum_{j, k \geq 1} \frac{\abss{\inprod{\phi_j}{\psi_k \circ
        \pitarget}}^2}{\numobs + j^{2\smoothorder / \dimx} +
    k^{2\smoothorder/ \dima}},\label{eq:continuum-arm-ate-minimax} \\
\label{eq:continuum-arm-cate-minimax}  
\minimaxRisk \big(  \xzero; \ball_\rkhs (1) \big) &\asymp
\numobs^{\frac{\dima + \dimx}{2 \smoothorder} - 1}, \quad \mbox{for
  any $\xzero \in \torus^{d_x}$}.
\end{align}
Furthermore, under the worst-case target policy, we have
\begin{align}
  \sup_{\pitarget} \minimaxRisk \big( \ball_\rkhs (1) \big)
  \asymp \numobs^{\frac{\dima}{2 \smoothorder} - 1}.
\end{align}
\end{subequations}
\end{corollary}
\noindent See~\Cref{app:subsec-proof-continuum-rates} for the proof.

\medskip

A few remarks are in order. For the one-point functional
$\LinFuncZero$, the optimal rate given by
equation~\eqref{eq:continuum-arm-cate-minimax} is exactly the optimal
rate for estimating a $(d_a + d_x)$-dimensional Sobolev function at
the point $(\xzero, \pitarget (\xzero))$. For the averaged functional
$\taustar$, in the worst case, we only need to pay for the dimension
$d_a$ of the action space, due to the averaging effect in the state
space. Moreover, the precise complexity for estimation is
characterized by equation~\eqref{eq:continuum-arm-ate-minimax}, which
depends on the behavior of the target policy $\pitarget$. Such an
instance-optimal risk is achieved by the estimator
$\tauhat_\numobs$. Finally, we remark that though the statement
of~\Cref{cor:continuum-bandit-rates} is for a deterministic target
policy $\pitarget$, the result naturally extends to general randomized
target policies.


\section{Simulation studies}
\label{subsec:simulation}

In this section, we present the results of some simulation studies in
which we compare our procedures with other methods. In particular, we
perform experiments on two classes of missing data problems, one
defined by the family of singular importance ratios discussed
in~\Cref{SecMissingOverlap} and the heavy-tailed example proposed by
Khan and Tamer~\cite{khan2010irregular}, with some
generalizations. These two examples allow us to explore two different
ways in which unbounded importance ratios can arise.

Concretely, we perform experiments in which the goal is to estimate
the treatment effect based on missing data. Let the state space be the real line $\Xspace = \real$, and let the action space
be binary, $\actionspace = \{0, 1\}$. We use the action $\action
\in \actionspace$ to model missingness, so that we only observe the
outcome $\outcome$ if and only if $\Action = 1$.  For simplicity, we
slightly abuse notation, and let $\treateff$ denote the function
$\treateff (\cdot, 1)$. Similarly, we use $\propscore$ to denote the
function $\propscore (1 \mid \cdot)$.  The goal is to estimate the linear
functional~\eqref{EqnDefnFunctional} with $\wenweight{\action}{\state}
= \action$, i.e.,
\begin{align*}
\taustar = \Exs_{\probx} \big[ \treateff(\State) \big].
\end{align*}
Throughout this section, we consider the homoskedastic case with
$\sigma^2 (\state, \action) \equiv 1$.  By
equation~\eqref{eq:semi-parametric-efficiency}, the semi-parametric
efficiency bound for this problem takes the form
\begin{align}
\label{eq:semiparametric-efficiency-bound-in-simulation-instances}  
\vsemi^2 = \var_\probxstar \big( \treateff(\State) \big) +
\int_{- \infty}^{\infty} \frac{\probx(\state)}{\propscore(\state) }
d\state,
\end{align}
which may or may not be finite.

For the rest of this section, we describe and discuss the construction
of simulation problem instances, as well as various choices of
estimators under our consideration. We then present the simulation
results.

\paragraph{Four possible estimators:}
We compare the performance of four possible estimators for the average
treatment effect --- two of which are based on inverse propensity
weights, while the other two (including our estimator) are based on
outcome regression.

First, we consider the na\"{i}ve inverse propensity weighting (IPW)
estimator, defined as
\begin{align}
  \tauhat_{\numobs, \mathrm{ipw}} = \frac{1}{\numobs} \sum_{i =
    1}^\numobs \frac{\outcome_i \Action_i}{\propscore (\State_i)}.
\end{align}
Note that the estimator $\tauhat_{\numobs, \mathrm{ipw}}$ always has
finite expectation, with $\Exs [\tauhat_{\numobs, \mathrm{ipw}}] =
\taustar$. Assuming that the outcome functions are bounded, the
variance of $\tauhat_{\numobs, \mathrm{ipw}}$, if exists, is given by
\begin{align*}
  \Exs \big[ \abss{\tauhat_{\numobs, \mathrm{ipw}} - \taustar}^2 \big]
  \asymp \numobs^{-1} \Big( \vsemi^2 + \Exs_{\probx} \big[
    \frac{1 + [\treateff(\State)]^2}{\propscore(\State)} \big] \Big)
  \asymp \frac{1}{\numobs} \int_{- \infty}^{\infty}
  \frac{\probx(\state)}{\propscore(\state)} d \state.
\end{align*}
In general, if the second moment does not exist, the na\"{i}ve IPW
estimator may converge to a heavy-tailed stable law, at a rate slower
than $\sqrt{\numobs}$. (c.f.~\cite{kallenberg1997foundations}, Chapter
14)

Khan and Tamer~\cite{khan2010irregular} suggested improving the
na\"{i}ve IPW by removing data with extremely small propensity
scores. Given a truncation level $\gamma_\numobs$, we define the
estimator
\begin{align}
  \tauhat_{\numobs, \mathrm{trunc}} = \frac{1}{\numobs} \sum_{i =
    1}^\numobs \frac{\outcome_i \Action_i}{\propscore (\State_i)}
  \bm{1}[\propscore(\State_i) \geq \gamma_\numobs],
\end{align}
where $\bm{1}[\propscore(X_i) \geq \gamma_\numobs]$ is equal to $1$
when $\propscore(\State_i) \geq \gamma_\numobs$, and zero otherwise.

Now we turn to the outcome-regression estimators based on kernel ridge
regression, as defined in the two-stage
framework~\eqref{eq:two-stage-homoskedastic}.  In order to improve the
universal constant factors (which are not covered by our theory), we
use a cross-fit procedure, i.e., we generate an estimator
$\dsone{\tauhat}_\numobs$ from the
framework~\eqref{eq:two-stage-homoskedastic}. By switching the role of
$(\dsone{\State}_i, \dsone{\Action}_i, \dsone{\outcome}_i)$ and
$(\dstwo{\State}_i, \dstwo{\Action}_i, \dstwo{\outcome}_i)$ and
applying the same two-stage framework, we obtain another estimator
$\dstwo{\tauhat}_\numobs$, and the final estimator is given by
\begin{align}
  \tauhat_\numobs = \frac{1}{2} \big( \dsone{\tauhat}_\numobs +
  \dstwo{\tauhat}_\numobs \big).
\end{align}

In terms of the regularization parameter $\reguhomo$, we consider two
possible choices:
\begin{itemize}
\item Optimal choice: based on the optimal theoretical prediction in
  equation~\eqref{eq:sample-size-req-based-on-eigendecay}, we set
  $\reguhomo = \tfrac{0.5}{\numobs}$ for any $\numobs > 0$. We
  call this estimator $\tauhat_{\numobs, \mathrm{opt}}$.
\item Cross validation: for each sample size $\numobs$, we use cross
  validation to find the regularization parameter that minimizes the
  mean-squared error in predicting $\treateff$.  We call this
  estimator $\tauhat_{\numobs, \mathrm{cv}}$.
\end{itemize}
It is worth noticing that the optimal choice of the regularization
parameter $\reguhomo$ for estimating the scalar $\taustar$ does
\emph{not} correspond to the optimal choice in estimating the
regression function $\treateff$. Indeed, as we will see in the
simulation results, the common cross-validation approach in
non-parametric estimation leads to sub-optimal semi-parametric
performance under our framework, and under-smoothing is crucial to the
optimal guarantees.

\subsection{Simulation results with heavy-tailed covariates}
\label{subsec:simu-heavytail}

We first present the simulation setup and results on the heavy-tailed
covariate examples proposed by Khan and
Tamer~\cite{khan2010irregular}.

\begin{figure}[ht!]
\begin{tabular}{ccc}
  \widgraph{0.45\textwidth}{\figdir/fig_logistic_logistic_ntrial_2000_maxn_128000}
  &&
  \widgraph{0.45\textwidth}{\figdir/fig_logistic_norm_ntrial_2000_maxn_128000}
  \\ (a) && (b) \\
  \widgraph{0.45\textwidth}{\figdir/fig_norm_logistic_ntrial_2000_maxn_128000}
  &&
  \widgraph{0.45\textwidth}{\figdir/fig_norm_norm_ntrial_2000_maxn_128000}
  \\
(c) && (d) \\

  \widgraph{0.45\textwidth}{\figdir/fig_cauchy_logistic_ntrial_2000_maxn_128000}
  &&
  \widgraph{0.45\textwidth}{\figdir/fig_cauchy_norm_ntrial_2000_maxn_128000}
  \\
  (e) && (f)
\end{tabular}
\caption{Plots of the mean-squared error $\Exs \big[
    \abss{\tauhat_{\numobs, \diamond} - \taustar}^2 \big]$ versus
  sample size $\numobs$.  Each curve corresponds to a different
  algorithm $\diamond \in \big\{ \mathrm{ipw}, \mathrm{trunc},
  \mathrm{opt}, \mathrm{cv} \big\}$. Each marker corresponds to a
  Monte Carlo estimate based on the empirical average of $2000$
  independent runs. For the cross-validated estimator
  $\tauhat_{\numobs, \mathrm{cv}}$, the choice of regularization
  parameter is based on averaging the cross validation results of the
  first $50$ runs. As indicated by the sub-figure titles, each panel
  corresponds to a problem setup $(\probx, \propscore) \in
  \big\{\probxlogistic, \probxnormal, \probxcauchy \big\} \times
  \big\{ \proplogistic, \propnormal \big\}$. Both axes in the plots
  are given by logarithmic scales. Some of the curves may overlap with
  each other.}
\label{fig:simulation}
\end{figure}

\paragraph{Model set-up:}

We consider the following choices for the distribution $\probxstar$
over data:
\begin{align*}
\mbox{Standard normal:}&\quad \probxnormal (x) = \frac{1}{\sqrt{2
    \pi}} \exp \big( - x^2 / 2\big), \\
\mbox{Standard logistic:} &\quad \probxlogistic (x) = \big( e^{x / 2}
+ e^{- x / 2} \big)^{-2},\\
\mbox{Standard Cauchy:}& \quad \probxcauchy (x) = \frac{1}{\pi (1 +
  x^2)}.
\end{align*}
Among these choices, the normal distribution possesses the lightest
tail, while the tail of the Cauchy distribution is the heaviest.

We carry out our simulation studies using the regression function
\begin{align*}
\treateff(\state) = 1 + \cos(\state) \quad \mbox{for $\state \in
  \real$.}
\end{align*}
This specific choice is not essential to our study; we have simply
chosen a bounded and smooth regression function.  Note that many
estimators under our consideration involve shrinkage, regularization,
or truncation steps, which make the output contract towards $0$. In
order to ensure a fair comparison, we include offset $1$ so that the
regression function is non-negative, and the target functional is
bounded away from zero.

In order to implement the kernel-based procedures, we use a Laplacian
kernel
\begin{align*}
\kernelfunc(\uvar, \uvar') & \mydefn \exp \big( - 2 |\uvar - \uvar'|
\big), \quad \mbox{for $\uvar, \uvar' \in \real$.}
\end{align*}
For the behavioral policy $\propscore$, we use the cumulative
distribution functions of logistic and normal distributions,
respectively.
\begin{align*}
    \propnormal (x) = \int_x^{\infty} \frac{1}{\sqrt{2 \pi}} \exp
    \big( - y^2 / 2\big) dy, \quad \mbox{and} \quad \proplogistic (x)
    = \frac{1}{1 + e^x}.
\end{align*}
Note that for both choices, the value $\propscore(\state)$ approaches
$0$ as $\state$ increases; the rate of decay is faster under the
normal model than the logistic model.

The paper~\cite{khan2010irregular} considers the density function
$\probx = \probxlogistic$, along with propensity score functions
$\propscore \in \big\{ \propnormal, \proplogistic \big\}$. Under both
setups, the semi-parametric efficiency bound derived in
equation~\eqref{eq:semiparametric-efficiency-bound-in-simulation-instances}
are infinite, while certain rates of convergence are still achieved
via truncation-based estimators (see Section 4.1
of~\cite{khan2010irregular} for details). Khan and Tamer proposed truncating at the threshold
$X_i \leq \sqrt{\log \numobs}$ for $\propscore = \propnormal$,
and $X_i \leq \log \numobs$ for $\propscore = \proplogistic$. Indeed,
they are the thresholds that ensures that the truncated inverse
propensity weight is uniformly bounded by a polynomial of $\numobs$,
under propensity scores $\proplogistic$ and $\propnormal$,
respectively. In our simulation studies,
we consider all possible combinations of $\probx \in \big\{
\probxnormal, \probxlogistic, \probxcauchy \big\}$ and $\propscore \in
\big\{ \propnormal , \proplogistic \big\} $. We choose $\gamma_\numobs =
\propscore(\log \numobs)$ under the logistic propensity score
$\proplogistic$, and $\gamma_\numobs = \propscore( \sqrt{\log
  \numobs})$ under the normal propensity score $\propnormal$, which yield near-optimal truncation levels,
regardless of the choice of data distribution $\probx$. Intuitively, heavier tail
of the distribution $\probx$ and lighter tail of the propensity score
$\propscore$ together lead to less regular behavior for estimators
based on important weighting.

Among our simulation setups, the only case that yields finite $\vsemi$
is that of $(\propscore = \proplogistic, \probx = \probxnormal)$;
other recent work~\cite{su2022when,jiang2022new} has also studied this
particular configuration. For the other five setups, the classical
theories for $\sqrt{\numobs}$-consistency and semi-parametric
efficiency are not available, due to the singular behavior of
propensity scores.

\paragraph{Simulation results:}
In Figure~\ref{fig:simulation}, we demonstrate the simulation results
for different estimators under aforementioned setups. The sample size
varies within the range $\numobs \in \{50, 100, 200, 400, 800, 1600,
3200, 6400, 12800\}$, and the mean-squared error is estimated through
empirical average over $2000$ independent runs.

From our simulation results, it can be observed that our estimator
$\tauhat_{\numobs, \mathrm{opt}}$ consistently outperforms other three
baselines. When tuning the regularization parameter using cross
validation, however, the estimator $\tauhat_{\numobs, \mathrm{cv}}$
performs significantly worse, over all the simulation instances. This
shows that under-smoothing is crucial to the performance of
outcome-regression estimators, and that the optimal bias-variance
trade-off in function and scalar estimation problems are drastically
different. The truncated IPW estimator also yields a reasonable and
robust performance, but in most settings, its rate of convergence
(represented as the slope of the curve in log-log plot) is worse than
$\tauhat_{\numobs, \mathrm{opt}}$. The na\"{i}ve IPW estimator, on the
other hand, can be highly unstable, especially for heavy-tailed data
distributions $\probxlogistic$ and $\probxcauchy$. Finally, we remark
that the two classes of estimators are not comparable in general, as
they use different information --- the IPW-based estimators
$\tauhat_{\numobs, \mathrm{ipw}}$ and $\tauhat_{\numobs,
  \mathrm{trunc}}$ use the information of the true propensity score
$\propscore$, which is not needed for $\tauhat_{\numobs,
  \mathrm{opt}}$ and $\tauhat_{\numobs, \mathrm{cv}}$; on the other
hand, the outcome regression estimators $\tauhat_{\numobs,
  \mathrm{opt}}$ and $\tauhat_{\numobs, \mathrm{cv}}$ require the
treatment effect function to lie in an RKHS, while the truncated IPW
estimator $\tauhat_{\numobs, \mathrm{trunc}}$ only requires it to be
bounded.

It is also useful to discuss the difference in the performance of
estimators under various setups. In the classical
$\sqrt{\numobs}$-regime with $\propscore = \proplogistic$ and $\probx
= \probxnormal$, the truncation does not happen with high probability,
and the na\"{i}ve IPW estimator yields the same MSE as the truncated
one, as shown in Figure~\ref{fig:simulation}(c). In other five cases,
the estimation error of $\tauhat_{\numobs, \mathrm{ipw}}$ is unstable,
and worse than the truncated analogue. It can be observed from that
the slopes of the green curves are around $1$ in the log-log plots in
panels (a)--(d) of Figure~\ref{fig:simulation}, but are much flatter
in panels (e) and (f).  This observation suggests that the optimal
rate of convergence may be near-parametric under the logistic and
normal model, while a slower minimax rate could be unavoidable in the
Cauchy setting.


\subsection{Simulation results with singular importance ratio}
\label{SecAlpha}

In this section, we report complete simulation results for the missing
data problem, but with singular importance ratios, as previously
described in~\Cref{SecIllustrative}---in particular, see
equation~\eqref{EqnIllustrative}.  We run the four estimators
discussed above, and compare their performance. The simulation setup
is essentially the same as Section~\ref{subsec:simu-heavytail}, with
the only difference being that the sample size varies within the range
$\numobs \in \{100, 200, 400, 800, 1600, 3200, 6400,
12800\}$.\footnote{We made this slight modification so as to avoid the
rare event that no outcome is observed.}  In defining the
truncation-based estimator $\tauhat_{\numobs, \mathrm{trunc}}$, we use
the truncation level $\gamma_\numobs = 1 / \sqrt{\numobs}$; this
choice yields the optimal rate of convergence among truncated IPW
estimators.

In~\Cref{fig:simulation-singular}, we present the results of our
simulations.  The problem instances are generated from the family of
singular models~\eqref{EqnIllustrative} with exponents $\overexp \in
\{0.5, 1,2, 3\}$. It can be seen that the simulation results match
well with our theoretical prediction: in the classical regime with
$\alpha = 0.5$, all the four estimators yield the same rate, while
$\tauhat_{\numobs, \mathrm{opt}}$ achieves slightly better
instance-dependent behavior; in the critical regime $\alpha = 1$, the
four estimators start to exhibit diverging behavior; in the harder
regimes of $\alpha \in \{2, 3\}$, the optimal estimator
$\tauhat_{\numobs, \mathrm{opt}}$ achieves the sharpest slope,
significantly outperforming the other three alternatives.

\begin{figure}[ht!]
\begin{tabular}{ccc}
  \widgraph{0.45\textwidth}{\figdir/fig_05_ntrial_2000_maxn_128000} &&
  \widgraph{0.45\textwidth}{\figdir/fig_10_ntrial_2000_maxn_128000}
  \\ (a) && (b) \\
  \widgraph{0.45\textwidth}{\figdir/fig_20_ntrial_2000_maxn_128000}
  &&
  \widgraph{0.45\textwidth}{\figdir/fig_30_ntrial_2000_maxn_128000}
  \\
(c) && (d)
\end{tabular}
\caption{Plots of the mean-squared error $\Exs \big[
    \abss{\tauhat_{\numobs, \diamond} - \taustar}^2 \big]$ versus
  sample size $\numobs$ for $\diamond \in \big\{ \mathrm{ipw}, \mathrm{trunc},
  \mathrm{opt}, \mathrm{cv} \big\}$. The simulation parameters are exactly the same as~\Cref{fig:simulation}, except for the underlying problem instances. As indicated by the sub-figure titles, each panel
  corresponds to an exponent $\alpha \in \{0.5, 1,2,3\}$. We have already presented part of the results (the cases of $\alpha = 0.5$ and $\alpha = 2$) in \Cref{SecIntro}.}
\label{fig:simulation-singular}
\end{figure}

\section{Proofs}
\label{sec:proofs}

In this section, we collect the proofs of our main results, with some
auxiliary results deferred to the appendices.


\subsection{Proof of~\Cref{ThmLower}}
\label{SecProofThmLower}

Throughout this section, we adopt the shorthand
$\minimaxRisk(\probInstanceStar, \funcClass) \equiv
\minimaxRisk(\funcClass)$, since $\probInstanceStar$ remains fixed
throughout.

\subsubsection{Proof of~\Cref{ThmLower}(a)}

This proof exploits some techniques introduced in our
previous paper~\cite{mou2022off}. Recalling that $c$ is a universal
constant, it suffices to prove the following two claims:
\begin{align}
\label{eq:minimax}
\minimaxRisk(\probInstanceStar, \funcClass) \stackrel{(a)}{\geq}
\frac{c}{\numobs} \vprobx^2(\treateff), \quad \mbox{and} \quad
\minimaxRisk(\probInstanceStar, \funcClass) \stackrel{(b)}{\geq}
\frac{c}{\numobs} \vwenlongsq.
\end{align}

Beginning with the bound~\eqref{eq:minimax}(a), we first observe that
the minimax risk over the class $\Nplain{\numobs}(\treateff, \probxstar)$
is lower bounded by the risk with fixed outcome function $\treateff$
and underlying distribution in the neighborhood of $\probxstar$, i.e.,
\begin{align*}
\minimaxRisk(\funcClass) \geq \minimaxRisk(\{\treateff\}) =
\inf_{\tauhat_\numobs} \sup_{(\probx, \treateff) \in \Nplain{\numobs} (\treateff, \probxstar)}
\Exs \big[ \abss{\tauhat_\numobs - \avgtreat( \treateff, \probx)}^2
  \big].
\end{align*}
But by Theorem 3 in our previous paper~\cite{mou2022off}, this minimax
risk is lower bounded by $\frac{c}{\numobs} \vprobx^2(\treateff)$,
which establishes the claim. \\

We now turn to proving the bound~\eqref{eq:minimax}(b), and we do so
via a version of Le Cam's two point lower bound.  More precisely, for
a fixed underlying distribution $\probxstar$, we construct a pair
$(\plainmu_+, \plainmu_-)$ of outcome functions within the
neighborhood $\big\{ \plainmu: \vecnorm{\plainmu - \treateff}{\Ltwospace (\Joint{\probxstar}{\propscore})} \leq \tfrac{\varbound^2}{\numobs} \big\} \cap \funcClass$ such that if we let $\Prob_{\mu, \xi}$ be the distribution of observations under the ground truth $(\mu, \xi)$, there is
\begin{align}
\label{eq:key-bounds}
  \totalvariation \Big(\Prob_{\plaintreateff_{+}, \probxstar}^{\otimes
    \numobs}, \Prob_{\plaintreateff_{-}, \probxstar}^{\otimes \numobs}
  \Big) \stackrel{(a)}{\leq} \frac{1}{2}, \quad \mbox{and} \quad
  \LinFunc (\plaintreateff_{+}, \probxstar) - \LinFunc(
  \plaintreateff_{-}, \probxstar) \stackrel{(b)}{\geq}
  \frac{c}{\sqrt{\numobs}} \vwenlong
\end{align}
Le Cam's two-point lemma (see e.g.~\cite{wainwright2019high}, Chapter 15) then implies
\begin{align*}
 \minimaxRisk(\funcClass) \geq \frac{1}{4} \Big\{1 - \totalvariation
 \Big(\Prob_{\plaintreateff_{+}, \probxstar}^{\otimes
   \numobs},\Prob_{\plaintreateff_{-}, \probxstar}^{\otimes \numobs}
 \Big) \Big\} \cdot \big \{  \LinFunc (\plaintreateff_{+}, \probxstar) - \LinFunc(
  \plaintreateff_{-}, \probxstar) \big \}^2 \geq \frac{c^2}{8
   \numobs} \vwenlongsq,
\end{align*}
completing the proof of equation~\eqref{eq:minimax}(b)

In order to prove the two bounds in line~\eqref{eq:key-bounds}, we
first need to specify the problem instances.

\paragraph{Construction of problem instances:}
We consider the noisy observation model
\begin{align}
\label{eq:gaussian-noise-in-lower-bound-construct}  
\outcome_i ~\mid~ \State_i, \Action_i \sim \mathcal{N} \Big(
\treateff(\State_i, \Action_i), \sigma^2(\State_i, \Action_i) \Big)
\qquad \mbox{for $i = 1, 2, \ldots, \numobs$.}
\end{align}
We may assume that $\vwenlong > 0$ without loss of generality (otherwise the lower bound is trivial).  By the
defining equation~\eqref{EqnDefnFunctional} and~\eqref{EqnVsuper}, and the symmetry of the
function class $\funcClass$, there exists a function $q_0: \Xspace
\times \actionspace \rightarrow \real$ such that
\begin{align*}
\Exs_{\probxstar} \Big[ \int_\actionspace q_0 (\State, \action) d
  \weightfunc(\action \mid \State) \Big] & \geq \tfrac{\vwenlong}{2},
\quad \mbox{and} \\ \tfrac{q_0}{\sqrt{\numobs}} \in \funcClass, \quad
& \Exs_{\Joint{\probxstar}{\propscore}} \Big[ \tfrac{q_0^2(\State,
    \Action)}{\sigma^2(\State, \Action)} \Big] \leq \tfrac{1}{4}.
\end{align*}
Using this function, we construct the outcome functions
\begin{align*}
    \plainmu_+ \mydefn \treateff + \tfrac{1}{2 \sqrt{\numobs}} q_0,
    \quad \mbox{and} \quad \plainmu_- \mydefn \treateff - \tfrac{1}{2
      \sqrt{\numobs}} q_0.
\end{align*}
Since $\treateff \in \tfrac{1}{2} \funcClass$ and $\frac{1}{2
  \sqrt{\numobs}} q_0 \in \tfrac{1}{2} \funcClass$, we have
$\plainmu_+, \plainmu_- \in \funcClass$ by convexity and symmetry. On
the other hand, we have the distance bound
\begin{align*}
\vecnorm{\treateff - \plainmu_+}{\Ltwospace
  (\Joint{\probxstar}{\propscore})}^2 = \frac{1}{4 \numobs}
\Exs_{\Joint{\probxstar}{\propscore}} \big[q_0^2(\State, \Action)
  \big] \leq \frac{\varbound^2}{4 \numobs}
\Exs_{\Joint{\probxstar}{\propscore}} \Big[ \frac{q_0^2(\State,
    \Action)}{\sigma^2(\State, \Action)} \Big] \leq
\frac{\varbound^2}{16 \numobs}.
\end{align*}
Consequently, we have $\plainmu_+ \in \Nval{\numobs} (\treateff) \cap
\funcClass$. Similarly, we also have $\plainmu_- \in \Nval{\numobs}
(\treateff) \cap \funcClass$.


\paragraph{Proof of equation~\eqref{eq:key-bounds}(a):}
We bound the KL divergence between the product measures. Let $ \law (\outcome
    | \State, \Action)$ denote the conditional law of $\outcome$ given the pair $(\State, \Action)$, we note that
\begin{align*}
\kull{\Prob_{\plaintreateff_{+}, \probxstar}^{\otimes
    \numobs}}{\Prob_{\plaintreateff_{-}, \probxstar}^{\otimes
    \numobs}} & \overset{(i)}{=} \numobs \cdot
\kull{\Prob_{\plaintreateff_{+},
    \probxstar}}{\Prob_{\plaintreateff_{-}, \probxstar}} \\
& \overset{(ii)}{\leq} \numobs \cdot \Exs \Big[ \kull{ \law (\outcome
    | \State, \Action) \big|_{\plainmu_+} }{ \law (\outcome | \State,
    \Action) \big|_{\plainmu_-} } \Big] \\
& = \numobs \cdot \frac{1}{4 \numobs} \cdot \Exs \Big[ \frac{q_0^2(\State,
    \Action)}{\sigma^2(\State, \Action)} \Big] \leq \frac{1}{4},
\end{align*}
where we use tensorization of KL divergence in step (i), and use
convexity of KL divergence in step (ii).

Applying Pinsker's inequality  yields
\begin{align*}
\totalvariation \Big(\Prob_{\plaintreateff_{+}, \probxstar}^{\otimes
  \numobs}, \Prob_{\plaintreateff_{-}, \probxstar}^{\otimes \numobs}
\Big) \leq \sqrt{\frac{1}{2}\kull{\Prob_{\plaintreateff_{+},
      \probxstar}^{\otimes \numobs}}{\Prob_{\plaintreateff_{-},
      \probxstar}^{\otimes \numobs}}} \leq \frac{1}{2 \sqrt{2}},
\end{align*}
which proves equation~\eqref{eq:key-bounds}(a).

\paragraph{Proof of equation~\eqref{eq:key-bounds}(b):}
Straightforward calculation yields
\begin{align*}
 \LinFunc (\plaintreateff_{+}, \probxstar) - \LinFunc(
  \plaintreateff_{-}, \probxstar)
= \frac{1}{\sqrt{\numobs}} \Exs_{\probxstar} \Big[ \int_\actionspace
  q_0 (\State, \action) d \weightfunc(\action \mid \State) \Big] \geq
\frac{1}{2 \sqrt{\numobs}} \vwenlong.
\end{align*}


\subsubsection{Proof of~\Cref{ThmLower}(b)}

Similar to the proof of~\Cref{ThmLower}, we use Le Cam's
two-point lemma. By the definition~\eqref{EqnVsuper} of the variance functional $\vsuper{\delta_\xzero}$, there exists a function $q_0: \Xspace
\times \actionspace \rightarrow \real$, such that
\begin{align*}
    \int_\actionspace q_0 (\xzero, \action) d \weightfunc(\action \mid \xzero)
     \geq \frac{\vsuper{\delta_\xzero}}{2},
    \quad \frac{q_0}{\sqrt{\numobs}} \in \funcClass, \quad\mbox{and}
    \quad \Exs_{\probxstar \cdot \propscore} \Big[
      \frac{q_0^2(\State, \Action)}{\sigma^2(\State, \Action)} \Big]
    \leq \frac{1}{4}.
\end{align*}
Using this function, we construct the outcome functions
\begin{align*}
    \plainmu_+ \mydefn \treateff + \frac{1}{2 \sqrt{\numobs}} q_0,
    \quad \mbox{and}\quad \plainmu_- \mydefn \treateff - \frac{1}{2
      \sqrt{\numobs}} q_0.
\end{align*}
Under the
construction~\eqref{eq:gaussian-noise-in-lower-bound-construct},
following the derivation of equation~\eqref{eq:key-bounds}(a), we have
\begin{align*}
    \totalvariation \Big(\Prob_{\plaintreateff_{+},
      \probxstar}^{\otimes \numobs}, \Prob_{\plaintreateff_{-},
      \probxstar}^{\otimes \numobs} \Big) \leq \frac{1}{2}.
\end{align*}
On the other hand, the gap satisfies
\begin{align*}
      \LinFunc (\plainmu_+, \delta_{\xzero}) - \LinFunc (\plainmu_-, \delta_{\xzero}) = \frac{1}{\sqrt{\numobs}}
      \int_\actionspace q_0 (\xzero, \action) d \weightfunc(\action \mid \xzero) \geq \frac{1}{2 \sqrt{\numobs}} \vsuper{\delta_{\xzero}}.
\end{align*}
Applying Le Cam's lemma yields the claim.


\subsection{Proof of~\Cref{PropExplicit,PropUniversal}}
\label{SecProofKernelProps}

In this section, we prove our two propositions that characterize the
variance functional in the case of an RKHS.


\subsubsection{Proof of~\Cref{PropExplicit}}
\label{SecProofPropExplicit}

The claim consists of two inequalities, and we split our proof
accordingly.


\paragraph{Proof of inequality~\eqref{EqnExplicit}(b):}
By definition, we have
\begin{align}
\vdouble{\newxdist}{\Hball} & = \sqrt{\numobs} \sup_{f \in \Hil}
\Big\{ \abss{\Exs_{\newxdist} \Big[\int_\actionspace \optfun(\State,
    \action) d \weightfunc(\action \mid \State) \Big]} \, \mid \,
\rkhsnorm{\optfun} \leq \Rbound \mbox{ and }
\Exs_{\Joint{\probxstar}{\propscore}} \Big[ \tfrac{\optfun^2(\State,
    \Action)}{\sigma^2(\State, \Action)} \Big] \leq \tfrac{1}{4
  \numobs} \Big\} \nonumber \\
\label{eq:simplify-variance-functional-in-rkhs}
& \geq \sup_{q \in \Hil} \Big\{ \Exs_{\newxdist}
\Big[\int_\actionspace q(\State, \action) d \weightfunc(\action \mid \State) \Big] \, \mid \, \tfrac{\rkhsnorm{q}^2}{\Rbound^2 \numobs}
+ 4 \Exs_{\Joint{\probxstar}{\propscore}} \Big[ \tfrac{q^2(\State,
    \Action)}{\sigma^2(\State, \Action)} \Big] \leq 1 \Big\},
\end{align}
where we have made the change of variable $\optfun =
q/\sqrt{\numobs}$.

Any function $q \in \rkhs$ has a basis expansion of the form $q =
\sum_{j=1}^{\infty} \theta_j \phi_j$, whence
\begin{align*}
    \Exs_{\newxdist} \Big[\int_\actionspace q(\State, \action) d
      \weightfunc(\action \mid \State) \Big] & =
    \inprod{\theta}{\featweighbar}_{\ell^2}, \quad \mbox{where
      $\featweighbar \equiv \featweighbar(\newxdist)$, and} \\
\tfrac{1}{\Rbound^2 \numobs} \rkhsnorm{q}^2 + 4
\Exs_{\Joint{\probxstar}{\propscore}} \Big[ \tfrac{q^2(\State,
    \Action)}{\sigma^2(\State, \Action)} \Big] & = \theta^\top \Big\{
(\Rbound^2 \numobs)^{-1} \EigenMat^{-1} + 4 \: \bigcov_{\sigma} \Big\}
\theta,
\end{align*}
where we use the eigen-value representation $\rkhsnorm{q}^2 = \theta^\top \EigenMat^{-1} \theta$.

We make the choice
\begin{align*}
\theta = \big\{(\Rbound^2 \numobs)^{-1} \EigenMat^{-1} + 4
\bigcov_{\sigma} \big\}^{- 1} \featweighbar /
\vecnorm{\big\{(\Rbound^2 \numobs)^{-1} \EigenMat^{-1} + 4
  \bigcov_{\sigma} \big\}^{- 1 / 2} \featweighbar}{\ell^2}.
\end{align*}
Substituting this choice into
equation~\eqref{eq:simplify-variance-functional-in-rkhs} yields
\begin{align*}
\vdoublesq{\newxdist}{\Hball} & \geq \featweighbar^\top
\big\{(\Rbound^2 \numobs)^{-1} \EigenMat^{-1} + 4 \bigcov_{\sigma}
\big\}^{-1} \featweighbar \geq \frac{1}{4} \featweighbar^\top
\big\{(\Rbound^2 \numobs)^{-1} \EigenMat^{-1} + \bigcov_{\sigma}
\big\}^{-1} \featweighbar,
\end{align*}
which establishes inequality (b).


\paragraph{Proof of inequality~\eqref{EqnExplicit}(a):}

Turning to the other inequality in the claim, the same change of
variable and followed by basis expansion yields
\begin{align*}
 \vdouble{\newxdist}{\Hball} & = \sup_{q \in \Hil} \Big\{
 \abss{\Exs_{\newxdist} \Big[\int_\actionspace q(\State, \action) d
     \weightfunc(\action \mid \State) \Big]} \; \mid \;
 \tfrac{1}{\Rbound\sqrt{\numobs}} \rkhsnorm{q} \leq 1,
 ~\Exs_{\Joint{\probxstar}{\propscore}} \Big[ \tfrac{q^2(\State,
     \Action)}{\sigma^2(\State, \Action)} \Big] \leq \tfrac{1}{4}
 \Big\} \nonumber\\
    & \leq \sup_{q \in \Hil} \Big\{ \Exs_{\newxdist}
 \Big[\int_\actionspace q(\State, \action) d \weightfunc( \action \mid \State) \Big] \; \mid \; \tfrac{1}{\Rbound^2 \numobs}
 \rkhsnorm{q}^2 + \Exs_{\Joint{\probxstar}{\propscore}} \Big[
   \tfrac{q^2(\State, \Action)}{\sigma^2(\State, \Action)} \Big] \leq
 \tfrac{5}{4} \Big\} \nonumber \\
& = \sup_{\theta \in \ell^2} \Big\{
 \inprod{\theta}{\featweighbar}_{\ell^2} \; \mid \; \theta^\top \Big(
 (\Rbound^2 \numobs)^{-1} \EigenMat^{-1} + \bigcovsig \Big) \theta \leq \tfrac{5}{4}
 \Big\} \nonumber \\
& \leq \tfrac{\sqrt{5}}{2} \cdot \sqrt{\featweighbar^\top
  \big\{(\Rbound^2 \numobs)^{-1} \EigenMat^{-1} + \bigcov_{\sigma}
  \big\}^{-1} \featweighbar},
\end{align*}
which completes the proof of inequality (a).


\subsubsection{Proof of~\Cref{PropUniversal}}
\label{SecProofPropUniversal}

Throughout this proof, we use $\|\cdot\|_2$ as a shorthand for the
$\Ltwospace(\Joint{\probxstar}{\propscore})$-norm. The variational formulation~\eqref{EqnVsuper} can be re-written as
\begin{align}
	\vsuperhil{\probxstar} = \sup \Big\{ \inprod{\optfun}{\smallimpratio} \mid \;
\Exs_{\Joint{\probxstar}{\propscore}} \big[ \tfrac{\optfun^2(\State,
    \Action)}{\sigma^2(\State, \Action)} \big] \leq \tfrac{1}{4}, \rkhsnorm{\optfun} \leq R \sqrt{\numobs} \Big\}.\label{eq:vsuper-rewritten-for-universal}
\end{align}

Clearly, the function
$\numobs \mapsto \vwenlonghil$ is non-decreasing, and since $\vsemi (\treateff, \smallimpratio) < + \infty$, it is uniformly bounded
from above. Therefore, by taking $\numobs \rightarrow + \infty$, the limit exists.  Moreover, we have
\begin{align*}
\lim_{\numobs \rightarrow \infty} \vwenlonghil \leq \tfrac{1}{2}
\sqrt{\Exs_{\Joint{\probxstar}{\propscore}} \Big[ \big(
    \smallimpratio(\Action \mid \State) \big)^2 \cdot \sigma^2(\State,
    \Action) \Big]}.
\end{align*}
Now the function $(\state, \action) \mapsto \sigma(\state, \action) \;
\tfrac{d \weightfunc}{d \propscore}(\action \mid \state)$ belongs to
$\Ltwospace (\Joint{\probxstar}{\propscore})$. Combined with the
uniform upper bound $\sup_{\state, \action} \sigma^2 (\state, \action)
\leq \varbound^2$, it follows that the function $\sigma^2 \; \tfrac{d
  \weightfunc}{d \propscore}$ also belongs to
$\Ltwospace(\Joint{\probxstar}{\propscore})$.  Since the Hilbert space
$\rkhs$ is universal (and hence dense in $\Ltwospace
(\Joint{\probxstar}{\propscore})$), it follows that for any
$\varepsilon > 0$, we can find a function $h_\varepsilon \in \rkhs$
such that $\|h_\varepsilon - \sigma^2 \smallimpratio \|_2 \leq
\varepsilon$.

Now define the rescaled function $q_\varepsilon \mydefn h_\varepsilon
/ (2 \|\sigma \tfrac{d \weightfunc}{d \propscore}\|_2)$.  With this
definition, we have
\begin{align*}
\Exs \Big[ \tfrac{q_\varepsilon^2(\State, \Action)}{\sigma^2(\State,
    \Action)} \Big] & \leq  \tfrac{1}{4} \|\sigma \tfrac{d \weightfunc}{d
  \propscore}\|_2^{-2} \Exs \Big[ \Big\{ \abss{ \sigma (\State, \Action) \cdot \tfrac{d
    \weightfunc}{d \propscore} (\Action \mid \State)} + \frac{\abss{h_\varepsilon - \sigma^2 \smallimpratio}}{\sigma} (\State, \Action) \Big\}^2 \Big]  \\
    & \leq\|\sigma \tfrac{d \weightfunc}{d
  \propscore}\|_2^{-2} \Big \{ \tfrac{1 + \varepsilon}{4} \Exs \big[
 \sigma^2 (\State, \Action) \cdot \big( \tfrac{d
    \weightfunc}{d \propscore} (\Action \mid \State) \big)^2 \big] +
\tfrac{1}{\varepsilon} \Exs \big[ \big( \tfrac{h_\varepsilon - \sigma^2
    \tfrac{d \weightfunc}{d \propscore}}{\sigma} \big)^2 (\State,
  \Action) \big] \Big \} \\
& \leq \frac{1 + \varepsilon}{4} + \frac{\varepsilon}{\lowersigma^2},
\end{align*}
along with the bound $\rkhsnorm{q_\varepsilon} \leq \|\sigma \tfrac{d
  \weightfunc}{d \propscore} \|_2^{-1} \; \rkhsnorm{h_\varepsilon} <
\infty$.

We now define the rescaled function $q \mydefn \frac{1}{1 +
  \varepsilon + \tfrac{4 \varepsilon}{\lowersigma^2}} q_\varepsilon$.
Given a sample size lower bounded as $\numobs \geq
\rkhsnorm{h_\varepsilon}^2 / (R^2 \vecnorm{\sigma \tfrac{d \weightfunc}{d
    \propscore}}{\Ltwospace(\Joint{\probxstar}{ \propscore})}^2)$, the
above inequalities imply that the rescaled function $q$ satisfies the
constraints in the optimization
problem~\eqref{eq:vsuper-rewritten-for-universal}. Substituting this choice into the objective
function, we find that
\begin{align*}
 \Exs_{\probxstar} \Big[\int_\actionspace q(\State, \action) d
   \weightfunc(\action \mid \State) \Big] & \geq \frac{1 - \varepsilon -
   4\varepsilon / \lowersigma^2}{2 \|\sigma \tfrac{d \weightfunc}{d
     \propscore}\|_2} \Exs_{\probxstar} \Big[\int_\actionspace
   h_\varepsilon (\State, \action) d \weightfunc(\action \mid \State)
   \Big] \\
 & \geq \frac{1 - \varepsilon - 4\varepsilon / \lowersigma^2}{2
   \|\sigma \tfrac{d \weightfunc}{d \propscore}\|_2} \left\{
 \|\smallimpratio \sigma\|_2^2 - \Exs \Big[ \big(h_\varepsilon -
   \sigma^2 \tfrac{d \weightfunc}{d \propscore}\big) \cdot \tfrac{d
     \weightfunc}{d \propscore}(\Action \mid \State) \Big] \right\}.
\end{align*}
The Cauchy--Schwarz inequality implies that
\begin{align*}
    \Exs \Big[ \big(h_\varepsilon - \sigma^2 \tfrac{d \weightfunc}{d
        \propscore}\big) \cdot \frac{d \weightfunc}{d \propscore}
      (\Action \mid \State) \Big] \leq \| h_\varepsilon - \sigma^2 \frac{d
      \weightfunc}{d \propscore}\|_2 \cdot \frac{1}{\lowersigma} \:
    \|\sigma \tfrac{d \weightfunc}{d \propscore}\|_2,
\end{align*}
where we have used the fact that $\sigma(\state, \action) \geq \lowersigma$
for all pairs $(\state, \action)$.

Combining the two bounds together and taking the limit, we conclude
that
\begin{align*}
   \lim_{\numobs \rightarrow \infty} \vwenlonghil & \geq
   \Exs_{\probxstar} \Big[\int_\actionspace q(\State, \action) d
     \weightfunc(\action | \State) \Big] \\
& \geq \tfrac{1 - \varepsilon - \varepsilon / \lowersigma^2 }{2} \;
   \|\sigma \tfrac{d \weightfunc}{d \propscore}\|_2 - \tfrac{1 -
     \varepsilon - \varepsilon / \lowersigma^2 }{2 \lowersigma}
   \varepsilon.
\end{align*}
Since the choice of $\varepsilon$ is arbitrary, this concludes the
proof of this proposition.

\subsection{Proof of~\Cref{ThmHomo} and variants}
\label{SecProofThmHomo}

Let us first introduce some notation used in the proof.  Our proof
involves the diagonal operator $\EigenMat^{-1} \mydefn \mathrm{diag}
\big( \{\eigen_k^{-1}\}_{k = 1}^{\infty} \big)$, and the weighted
$\ell^2$-norms
\begin{subequations}
\begin{align}
 \vecnorm{\plainvec}{\eigen}^2 \mydefn \sum_{j = 1}^{\infty} \eigen_j
 \plainvec_j^2, \quad \mbox{and} \quad
 \vecnorm{\plainvec}{\eigen^{-1}}^2 \mydefn \sum_{j = 1}^{\infty}
 \eigen_j^{-1} \plainvec_j^2.\label{eq:weighted-norm-eigen-defn}
\end{align}
\end{subequations}


\subsubsection{Set-up for auxiliary results}\label{subsubsec:setup-for-aux-results}

We define the empirical feature vector
\mbox{$\widehat{\featweigh}_{\numobs} \mydefn \frac{1}{\numobs}
  \sum_{i = 1}^\numobs \int_\actionspace \phi(\dstwo{\State}_i,
  \action) d \weightfunc( \action \mid \dstwo{\State}_i)$.}  By the kernel
boundedness assumption~\ref{EqnKerBou}, we have
$\vecnorm{\featweighbar}{\eigen} < \infty$ and
$\vecnorm{\featweighhat_{\numobs}}{\eigen} < \infty$ almost surely.
We also define a linear operator $\Psi$ from the Hilbert space $\rkhs$
to the sequence space $\ell^2$ with components $[\Psi(f)]_j =
\inprod{f}{\phi_j}_{\Ltwospace(\Joint{\probxstar}{\propscore})}$.  Since
$\treateff$ and $\muhat_{\numobs}$ belong to the Hilbert space
$\rkhs$, it is meaningful to define
\begin{align*}
    \projparam_* & \mydefn \Psi(\treateff), \quad \mbox{and} \quad
    \widehat{\projparam}_{\numobs} \mydefn \Psi(\muhat_{\numobs}).
\end{align*}
Note that for any function $f \in \rkhs$, we have
\begin{align*}
\vecnorm{\Psi
    f}{\eigen^{-1}}^2 = \sum_{j = 1}^{+ \infty} \eigen_j^{-1}\inprod{f}{\phi_j}^2_{\Ltwospace (\Joint{\probxstar}{\propscore})} = \sum_{j = 1}^{\infty}
  \rkhsinprod{f}{\sqrt{\eigen_j} \phi_j}^2 = \rkhsnorm{f}^2 <
  \infty.
\end{align*}
Consequently, the inner products
$\inprod{\featweigh}{\projparam} \leq \vecnorm{\featweigh}{\eigen} \cdot \vecnorm{\projparam}{\eigen^{-1}}$ are well-defined for $\featweigh \in
\{\featweighbar, \featweighhat_{\numobs}\}$ and $\projparam \in \{
\projparam_*, \projparamhat_\numobs \}$. Fubini's theorem guarantees
that target functional $\taustar$ and the estimator $\tauhat_\numobs$
can be written as
\begin{align*}
    \taustar & = \sum_{k = 0}^{\infty}
    \inprod{\treateff}{\basisfunc_k}_{\Ltwospace
      (\Joint{\probxstar}{\propscore})} \cdot \Exs_{\probxstar} \left[
      \int_{\actionspace} \basisfunc_k(\State, \action) d
      \weightfunc( \action \mid \State) \right] =
    \inprod{\featweighbar}{\projparam_*}, \quad \mbox{and}\\
     \tauhat_\numobs &
    = \numobs^{-1} \sum_{k = 0}^{\infty} \inprod{\muhat_{\numobs}^{(1)}}{\basisfunc_k}_{\Ltwospace
      (\Joint{\probxstar}{\propscore})} \sum_{i = 1}^\numobs
    \int_\actionspace \basisfunc_k (\dstwo{\State}_i, \action)
    d\weightfunc(\action \mid \dstwo{\State}_i) \; = \;
    \inprod{\featweighhat_{\numobs} }{\projparamhat_{\numobs}}.
\end{align*}
We therefore have the following error decomposition:
\begin{align}
\label{eq:linear-regression-error-decomp}
  \tauhat_{ \numobs} - \taustar = \inprod{\featweighhat_{\numobs} -
    \widebar{\featweigh}}{\projparam_*} +
  \inprod{\widebar{\featweigh}}{\projparamhat_\numobs - \projparam_*}
  + \inprod{\featweighhat_\numobs - \widebar{\featweigh}
  }{\projparamhat_\numobs - \projparam_*}.
\end{align}
The rest of this section is devoted to bounds on each terms appearing
in this decomposition.  In particular, we require two auxiliary
results.  Recall the shorthand notation \mbox{$\vprobx^2(f) \defn \var_{X \sim \probxstar}
  \big( \int_\actionspace f(\State, \action) d \wenweight{\action}{\State} \big)$.}
\begin{lemma}
\label{lemma:avg-in-rkhs-directional}
Under Assumptions~\eqref{EqnKerBou} and~\eqref{EqnCleanSubGauss}, for
any function $f \in \rkhs$, we have
\begin{subequations}
\begin{align}
\label{eq:avg-directional-bound}    
\abss{ \inprod{\featweighhat_\numobs - \widebar{\featweigh}}{\Psi f}}
\leq 2 \vprobx (f) \sqrt{\frac{\log (1 / \delta)}{\numobs}} + 6 \kappa
\rkhsnorm{f} \frac{\log (1 / \delta)}{\numobs}.
\end{align}
with probability at least $1 - \delta$.  Furthermore, given a sample size
$\numobs \geq \log (1 / \delta)$ and a scalar $\regu > 0$, we have
\begin{align}
\label{eq:avg-preconditioned-bound}  
  \vecnorm{(\IdMat + \regu \EigenMat^{-1})^{- 1/2}
    \big(\featweighhat_\numobs - \widebar{\featweigh} \big)}{\ell^2}
  \leq \sqrt{\frac{\effdim (\regu)}{\numobs} \log (1 / \delta)}
\end{align}
 with probability $1 - \delta$.
\end{subequations}    
\end{lemma}
\noindent See~\Cref{subsubsec:proof-lemma-avg-directional} for the
proof.

\medskip

\begin{lemma}
\label{lemma:krr-directional-bound}
 Suppose that the kernel bound~\eqref{EqnKerBou} and tail
 condition~\eqref{EqnCleanSubGauss} are in force.  Then for any
 infinite-dimensional vector $\plainvec$ and scalar $\delta \in (0,
 1)$, with the regularization parameter \mbox{$\regu_\numobs =
   \frac{\varbound^2}{\Rbound^2 \numobs}$,} and under the sample-size
 condition~\eqref{eq:sample-size-req-based-on-eigendecay}, we have
 \begin{align}
   \abss{\inprod{\plainvec}{ \projparamhat_\numobs - \projparam_*
   }} \leq c \vecnorm{\big( \IdMat + \regu_\numobs \EigenMat^{-1}
     \big)^{-1 / 2} \plainvec}{\ell^2} \cdot \Big\{ \varbound
   \sqrt{\frac{\log (1 / \delta)}{\numobs}} + \subgaussian
   \sqrt{\effdim(\regu_\numobs)} \cdot \frac{\log (1 / \delta)
     \log \numobs}{\numobs} \Big\},
 \end{align}
with probability at least $1 - \delta$.
\end{lemma}
\noindent See~\Cref{subsubsec:proof-lemma-krr-directional} for the
proof.


\subsubsection{Main argument}

Taking these two lemmas as given, we now prove~\Cref{ThmHomo} by
bounding each term in the decomposition
result~\eqref{eq:linear-regression-error-decomp}.  First, recalling
that $\projparam_* = \Psi \treateff$, we can apply the
bound~\eqref{eq:avg-directional-bound} to find that
\begin{align}
\label{eq:ubar-beta-noise-bound}  
\abss{ \inprod{\featweighhat_\numobs -
    \widebar{\featweigh}}{\projparam_*}} \leq 2 \vprobx (\treateff)
\sqrt{\frac{\log (1 / \delta)}{\numobs}} + 6 \kappa \Rbound \frac{\log
  (1 / \delta)}{\numobs}
\end{align}
with probability at least $1 - \delta / 3$.  Second, by the boundedness of
basis functions, we have $\featweighbar \in \ell^\infty$, so that
\Cref{lemma:krr-directional-bound} can be applied to obtain
\begin{align}
\abss{\inprod{\featweighbar}{ \projparamhat_\numobs - \projparam_* }}
& \leq c \vecnorm{\big( \IdMat + \regu_\numobs \EigenMat^{-1}
  \big)^{-1 / 2} \featweighbar}{\ell^2} \cdot \Big\{ \varbound
\sqrt{\frac{\log (1 / \delta)}{\numobs}} + \subgaussian
\sqrt{\effdim(\regu_\numobs)} \cdot \frac{\log (1 / \delta) \log
  \numobs}{\numobs} \Big\} \nonumber\\ &\leq 2c \vecnorm{\big( \IdMat
  + \regu_\numobs \EigenMat^{-1} \big)^{-1 / 2} \featweighbar}{\ell^2}
\varbound \sqrt{\frac{\log (1 /
    \delta)}{\numobs}}, \label{eq:u-noise-betastar-bound}
\end{align}
with probability $1 - \delta / 3$, where the last step follows from the sample-size
condition~\eqref{eq:sample-size-req-based-on-eigendecay}, as the condition ensures that $\subgaussian \sqrt{\effdim(\regu_\numobs)} \cdot \tfrac{\log (1 / \delta) \log
  \numobs}{\numobs} \leq  \varbound
\sqrt{\tfrac{\log (1 / \delta)}{\numobs}}$.

Next, we apply equation~\eqref{eq:avg-preconditioned-bound} in
combination with~\Cref{lemma:krr-directional-bound}, and obtain
the following inequality with probability $1 - \delta / 3$
\begin{align}
    \abss{\inprod{\featweighhat_\numobs - \widebar{\featweigh}
      }{\projparamhat_\numobs - \projparam_*}} &\leq c \vecnorm{\big(
      \IdMat + \regu_\numobs \EigenMat^{-1} \big)^{-1 / 2} \big(
      \featweighhat_\numobs - \widebar{\featweigh}\big)}{\ell^2} \cdot
    \Big\{ \varbound \sqrt{\frac{\log (1 / \delta)}{\numobs}} +
    \subgaussian \sqrt{\effdim(\regu_\numobs)} \cdot \frac{\log (1 /
      \delta) \log \numobs}{\numobs} \Big\} \nonumber\\ &\leq c
    \sqrt{\frac{\effdim(\regu_\numobs)}{\numobs} \log (1 / \delta)}
    \cdot \Big\{ \varbound \sqrt{\frac{\log (1 / \delta)}{\numobs}} +
    \subgaussian \sqrt{\effdim(\regu_\numobs)} \cdot \frac{\log (1 /
      \delta) \log \numobs}{\numobs} \Big\}
    \nonumber\\ &\overset{(i)}{\leq} c \subgaussian
    \sqrt{\effdim(\regu_\numobs)} \frac{ \log (1 / \delta) \log
      \numobs}{\numobs} \Big\{ 1 + \sqrt{\frac{\effdim(\regu_\numobs)
        \log (1 / \delta)}{\numobs}}
    \Big\}\nonumber\\ &\overset{(ii)}{\leq} 2 c \subgaussian
    \sqrt{\effdim(\regu_\numobs)} \frac{ \log (1 / \delta) \log
      \numobs}{\numobs}. \label{eq:u-noise-beta-noise-bound}
\end{align}
where step (i) follows from the relation $\varbound \leq
\subgaussian$, whereas the final step follows from the sample size
condition~\eqref{eq:sample-size-req-based-on-eigendecay}, as it ensures $\tfrac{\effdim(\regu_\numobs)
        \log (1 / \delta)}{\numobs} \leq 1$.

Combining the
inequalities~\eqref{eq:ubar-beta-noise-bound},~\eqref{eq:u-noise-betastar-bound},
and~\eqref{eq:u-noise-beta-noise-bound} completes the proof
of~\Cref{ThmHomo}.


\subsubsection{Proof of~\Cref{lemma:avg-in-rkhs-directional}}
\label{subsubsec:proof-lemma-avg-directional}

We simplify notation by omitting the superscript $\dstwo{}$, and using
$(\State_i, \Action_i, \outcome_i)$ to denote the data.

\paragraph{Proof of the directional bound~\eqref{eq:avg-directional-bound}}:
By definition, we have
\begin{align*}
\inprod{\featweighhat_\numobs }{\Psi f} & = \frac{1}{\numobs} \sum_{i
  = 1}^\numobs \int_{\actionspace} \inprod{\phi(\State_i,
  \action)}{\Psi f} d \wenweight{\action}{\State_i}
\\
& = \frac{1}{\numobs} \sum_{i=1}^\numobs Z_i \qquad \mbox{where $Z_i
  \defn \int_{\actionspace} f(\State_i, \action) d
  \wenweight{\action}{\State_i}$.}
\end{align*}
Similarly, the population-level vector $\featweighbar$ satisfies
$\inprod{ \featweighbar}{\Psi f} = \Exs_{\probxstar} [Z]$, so that our
problem amounts to bounding the fluctuations of the sample average
$\frac{1}{\numobs} \sum_{i=1}^\numobs Z_i$ around its mean.  Our
approach is via Bernstein's inequality, and applying it requires
control on both the variance and absolute value of $Z_i$.  By
inspection, we have \mbox{$\var(Z_i) = \vprobx^2(f) = \var \big(
  \int_\actionspace f(\State, \action) d \wenweight{\action}{\State}
  \big)$} and moreover, we claim that
\begin{align}
\label{EqnZbound}  
|Z_i| & \leq \kappa \cdot \rkhsnorm{f}.
\end{align}
With these two bounds in hand, invoking Bernstein's inequality (see e.g.~\cite{massart2007concentration}) yields
\begin{align*}
\Prob \left( \abss{\inprod{\featweighhat_\numobs - \featweighbar}{\Psi
    f}} \geq t \right) \leq 2 \exp \left(\frac{- \numobs t^2}{2
  \vprobx^2(f) + 3 \kappa \rkhsnorm{f} t} \right),
\end{align*}
setting the right hand side as $\delta$ and solving for $t$, we complete the proof of the directional
bound~\eqref{eq:avg-directional-bound}.

It remains to prove the claim~\eqref{EqnZbound}.  Using our
boundedness condition~\eqref{eq:signed-radon-msr-bounded} on $\weightfunc$, we have
\begin{align*}
|Z_i| = \abss{\int_{\actionspace} f(\State_i, \action) d
  \wenweight{\action}{\State_i}} & \leq \sup_{(\state, \action)
  \in \statespace \times \actionspace} \abss{f(\state, \action)} \\
& \stackrel{(i)}{\leq} \sup_{(\state, \action) \in \statespace
  \times \actionspace} \sum_{k \geq 1} \eigen_k \abss{\basisfunc_k
  (\state, \action)} \cdot \abss{ \rkhsinprod{f}{\basisfunc_k}} \\
& \overset{(ii)}{\leq} \sup_{(\state, \action) \in \statespace
      \times \actionspace} \Big(\sum_{k \geq 1} \eigen_k
    \phi_k^2(\state, \action) \Big)^{1/2} \cdot \Big( \sum_{k \geq 1}
    \eigen_k \rkhsinprod{f}{\basisfunc_k}^2 \Big)^{1/2}
\end{align*}
where step (i) follows by expanding $f$ into a basis representation $f = \sum_{k = 1}^{+ \infty} \eigen_k \rkhsinprod{f}{\basisfunc_k} \basisfunc_k$;
and step (ii) follows from the Cauchy--Schwarz inequality.  Now by Mercer's theorem, we have the relation
\begin{align*}
\big( \sum_{k \geq 1} \eigen_k \rkhsinprod{f}{\basisfunc_k}^2
\big)^{1/2} = \rkhsnorm{f},
\end{align*}
whereas the boundedness condition~\eqref{EqnKerBou} implies that $
\sum_{k \geq 1} \eigen_k \phi_k^2(\state, \action) \leq \kernelfunc \big( (\state, \action), (\state, \action) \big) \leq
\kappa^2$, for any $(\state, \action) \in \Xspace \times \ActionSpace$.  Putting together the pieces yields the claimed
bound~\eqref{EqnZbound}.

\paragraph{Proof of the preconditioned bound~\eqref{eq:avg-preconditioned-bound}:}

Defining $U_i \mydefn \int_\actionspace \phi(\State_i, \action) d
\wenweight{\action}{\State_i}$ so that $\featweighhat_\numobs = \tfrac{1}{\numobs} \sum_{i = 1}^\numobs U_i$, the norm on the left-hand-side of Eq~\eqref{eq:avg-preconditioned-bound} can be equivalently written
as an empirical process supremum
\begin{align*}
\vecnorm{(\IdMat + \regu \EigenMat^{-1})^{- 1/2} \big(\featweighhat_\numobs
  - \widebar{\featweigh} \big)}{\ell^2} = \sup_{z^\top (\IdMat + \regu
  \EigenMat^{-1}) z \leq 1} \frac{1}{\numobs} \sum_{i= 1}^\numobs z^\top
(U_i - \featweighbar) =: H_\numobs
\end{align*}
In order to bound the expected supremum, we simply use the Cauchy--Schwarz
inequality to arrive at the bound
\begin{align*}
    \Exs \big[H_\numobs \big] \leq \Big\{ \Exs \Big[ \vecnorm{(\IdMat
        + \regu \EigenMat^{-1})^{- 1/2} \cdot \frac{1}{\numobs} \sum_{i =
          1}^\numobs (U_i - \featweighbar)}{\ell^2}^2 \Big]
    \Big\}^{1/2} \leq \sqrt{\numobs^{-1} \Exs \big[ \vecnorm{(\IdMat +
          \regu \EigenMat^{-1})^{-1/2} U_i}{\ell^2}^2 \big]},
\end{align*}
where the second inequality comes from the fact that $U_i$'s are $\mathrm{i.i.d.}$

In order to bound this quantity, we use Talagrand's concentration
inequality (c.f.~\cite{wainwright2019high}, Theorem 3.8 and
remarks). With probability $1 - \delta$, we have
\begin{multline}
\label{eq:conc-in-avg-preconditioned-bound-proof}  
H_\numobs \leq 2 \Exs \big[ H_\numobs \big] + c \Big(\sup_{z^\top
  (\IdMat + \regu \EigenMat^{-1}) z \leq 1} \Exs[ (z^\top (U_i -
  \featweighbar))^2] \tfrac{\log (1 / \delta)}{\numobs} \Big)^{1/2} \\
+ c \sup_{(\state, \action) \in \Xspace \times \actionspace}
\vecnorm{(\IdMat + \regu \EigenMat^{-1})^{- 1/2} \phi(\state,
  \action)}{\ell^2} \cdot \tfrac{\log (1/\delta)}{\numobs}.
\end{multline}

Since the Radon measure $\weightfunc$ satisfies the
bound~\eqref{eq:signed-radon-msr-bounded}, the summand $U_i$ satisfies
the almost sure upper bound
\begin{align*}
    \vecnorm{(\IdMat + \regu \EigenMat^{-1})^{- 1/2} U_i}{\ell^2} \leq
    \sup_{(\state, \action) \in \statespace \times \actionspace}
    \vecnorm{(\IdMat + \regu \EigenMat^{-1})^{- 1/2} \phi(\state,
      \action)}{\ell^2} = \sqrt{\effdim (\regu)}
\end{align*}
The bound~\eqref{eq:conc-in-avg-preconditioned-bound-proof} then
becomes
\begin{align*}
    H_\numobs \leq c \sqrt{\effdim (\regu)} \cdot \Big\{
    \sqrt{\frac{\log (1 / \delta)}{\numobs}} + \frac{\log (1 /
      \delta)}{\numobs}\Big\},
\end{align*}
which completes the proof of
equation~\eqref{eq:avg-preconditioned-bound}.

\subsubsection{Proof of~\Cref{lemma:krr-directional-bound}}\label{subsubsec:proof-lemma-krr-directional}
For notational simplicity, we omit the supscript $\dsone{}$ in
$(\State_i, \Action_i, \outcome_i)$.  Using the basis expansion
$\plaintreateff = \sum_{k \geq 0} \projparam (k) \basisfunc_k$, we
have the equivalence
\begin{align*}
\projparamhat_\numobs &= \Psi \cdot \arg\min_{\plaintreateff \in
  \rkhs} \Big\{ \frac{1}{\numobs} \sum_{i = 1}^{\numobs} \big(
\outcome_i - \plaintreateff(\State_i, \Action_i) \big)^2 +
\regu_\numobs \rkhsnorm{\plaintreateff}^2 \Big\} \\
& = \arg\min_{\projparam \in \ell^2(\naturals)} \Big\{
\frac{1}{\numobs} \sum_{i = 1}^{\numobs} \big( \outcome_i -
\inprod{\projparam}{\phi(\State_i, \Action_i)} \big)^2 + \regu_\numobs
\vecnorm{\projparam}{\eigen^{-1}}^2 \Big\}.
\end{align*}
Define the noise variable $\noise_i \mydefn \outcome_i -
\treateff(\State_i, \Action_i)$ along with the empirical covariance
operator $\widehat{\bigcov}_{\numobs} \mydefn \frac{1}{\numobs}
\sum_{i = 1}^\numobs \phi(\State_i, \Action_i) \phi(\State_i,
\Action_i)^\top$.  Using this notation, we can write
\begin{align*}
\projparamhat_\numobs - \projparam_* = \big( \widehat{\bigcov}_\numobs
+ \regu_\numobs \EigenMat^{-1} \big)^{-1} \cdot \frac{1}{\numobs}
\sum_{i = 1}^\numobs \Big\{ \noise_i \phi(\State_i, \Action_i) -
\regu_\numobs \EigenMat^{-1} \projparam_* \Big\}.
\end{align*}
Our approach to controlling the projection of this quantity in any
fixed direction $\plainvec$ consists of two steps:
\bcar
\item First, conditionally on the state-action pairs $(\State_i,
  \Action_i)_{i = 1}^\numobs$, we exhibit a high-probability upper
  bound on the error $\plainvec^\top \big(\projparamhat_\numobs -
  \projparam_* \big)$ with respect to the randomness in the outcomes
  $\outcome_i$. The bound depends on the behavior of the empirical
  covariance operator $\widehat{\bigcov}_{\numobs}$ of feature
  vectors; see
 ~\Cref{lemma:conditional-high-prob-bound-for-krr-directional}
  for details.
\item Second, we relate the empirical covariance operator
  $\widehat{\bigcov}_{\numobs}$ with its population analogue (which is
  the identity operator $\IdMat$, since $(\phi_j)_{j = 1}^{\infty}$
  forms an orthonormal basis). The form of infinite-dimensional
  concentration results is exactly the form required in the first
  step. See~\Cref{lemma:emp-cov-concentration} for details. \\
\ecar

\noindent Let us give precise statements of the two auxiliary results
needed in the proof:
\begin{lemma}
\label{lemma:conditional-high-prob-bound-for-krr-directional}
Conditionally on the state-action sequence $(\State_i, \Action_i)_{i =
  1}^\numobs$, for any $\plainvec \in \ell^\infty (\naturals)$, we
have
\begin{multline}
  \abss{\plainvec^\top (\projparamhat_\numobs - \projparam_*)} \leq c
  \vecnorm{\big( \widehat{\bigcov}_{\numobs} + \regu_\numobs \EigenMat^{-1}
    \big)^{-1 / 2} \plainvec}{\ell^2} \\
  \times \Big\{ \sqrt{\regu_\numobs} \rkhsnorm{\treateff} + \varbound
  \sqrt{\tfrac{\log (1 / \delta)}{\numobs}} + \subgaussian
  \sup_{(\state, \action)} \vecnorm{\big( \widehat{\bigcov}_\numobs +
    \regu_\numobs \EigenMat^{-1} \big)^{-1/2} \phi(\state,
    \action)}{\ell^2} \cdot \tfrac{\log (1 / \delta) \log
    \numobs}{\numobs} \Big\},
\end{multline}
with probability at least $1 - \delta$.
\end{lemma}
\noindent
See~\Cref{app:subsec-proof-conditional-high-prob-bound-for-krr-directional}
for the proof.

\medskip

Our next auxiliary result relates the sample covariance operator
$\widehat{\bigcov}_\numobs$ with the population one.   Here we
state a somewhat general result, since we use it both here and
in our later proof of~\Cref{ThmHetero}.

Consider a weight function $(\state, \action) \mapsto q(\state,
\action) \in [\underline{q}, \overline{q}]$, where $(\underline{q},
\overline{q})$ are a pair of positive scalars. Define the empirical
operator
\begin{align*}  
\widehat{\bigcov}_{\numobs, q} \mydefn \numobs^{-1} \sum_{i =
  1}^\numobs q(\State_i, \Action_i) \phi(\State_i, \Action_i)
\phi(\State_i, \Action_i)^\top,
\end{align*}
along with its its population version $\bigcov_{*, q} \mydefn \Exs
\big[ \widehat{\bigcov}_{\numobs, q} \big]$.  For the current proof,
it suffices to take $q(\state, \action) = 1$.
\begin{lemma}
\label{lemma:emp-cov-concentration}
For scalars $\delta, \offpar \in (0, 1)$, consider a regularization
parameter $\regu_\numobs$ satisfying the relation
\begin{align}
\label{eq:regu-requirement-in-emp-cov-concentration-general}  
(\overline{q} / \underline{q}) \log \big( \tfrac{\kappa^2
}{\regu_\numobs \delta} \big) \cdot \tfrac{\effdim \big(\regu_\numobs
  / \underline{q} \big)}{\numobs} & \leq \tfrac{\offpar}{16}.
\end{align}
Then we have
\begin{align}
\label{eq:matrix-dominance-in-emp-cov-concentration}  
  (1 - \offpar) \big( \bigcov_{*, q} + \regu_\numobs \EigenMat^{-1} \big)
\preceq \widehat{\bigcov}_{\numobs, q}+ \regu_\numobs \EigenMat^{-1} \preceq
(1 + \offpar) \big( \bigcov_{*, q} + \regu_\numobs \EigenMat^{-1} \big)
\end{align}
with probability at least $1 - \delta$.
\end{lemma}
\noindent See~\Cref{subsubsec:proof-lemma-emp-cov-concen} for the
proof.

\medskip

Taking these two lemmas as given, we now proceed with the proof
of~\Cref{lemma:krr-directional-bound}.  We define the event
\begin{align*}
\Event \mydefn \Big\{ \widehat{\bigcov}_\numobs + \regu_\numobs
\EigenMat^{-1} \succeq \tfrac{1}{2} (\IdMat + \regu_\numobs
\EigenMat^{-1}) \Big\}.
\end{align*}
With the given choice $\regu_{\numobs} = \frac{\varbound^2}{\Rbound^2
  \numobs}$, for a sample size $\numobs$ satisfying the
requirement~\eqref{eq:sample-size-req-based-on-eigendecay}, we have
\begin{align*}
\log \big(\frac{\kappa^2}{\regu_\numobs \delta}\big)
\frac{\effdim(\regu_\numobs)}{\numobs} \leq \frac{1}{32}.
\end{align*}
By applying~\Cref{lemma:emp-cov-concentration} with $q(\state,
\action) \equiv 1$, we are guaranteed that $\Prob (\Event) \geq 1 -
\delta$.

Conditioned on the event $\Event$, the definition of effective
dimension guarantees that
\begin{align*}
  \sup_{(\state, \action) \in \statespace \times \actionspace}
  \vecnorm{\big( \widehat{\bigcov}_\numobs + \regu_\numobs \EigenMat^{-1}
    \big)^{-1/2} \phi(\state, \action)}{\ell^2} \leq \sqrt{2}
  \sup_{(\state, \action) \in \statespace \times \actionspace}
  \vecnorm{\big( \IdMat + \regu_\numobs \EigenMat^{-1} \big)^{-1/2}
    \phi(\state, \action)}{\ell^2} \leq \sqrt{2 \effdim
    (\regu_\numobs)}.
\end{align*}
Consequently, conditioned on the event $\Event$,
\Cref{lemma:conditional-high-prob-bound-for-krr-directional}
guarantees that
\begin{align*}
\abss{\plainvec^\top (\projparamhat_\numobs - \projparam_*)} & \leq c
\vecnorm{\big( \widehat{\bigcov}_{\numobs} + \regu_\numobs
  \EigenMat^{-1} \big)^{-1 / 2} \plainvec}{\ell^2} \cdot \Big\{
\sqrt{\regu_\numobs} \rkhsnorm{\treateff} + \varbound
\sqrt{\tfrac{\log (1 / \delta)}{\numobs}} + \subgaussian \sqrt{2
  \effdim(\regu_\numobs)} \cdot \tfrac{\log (1 / \delta) \log
  \numobs}{\numobs} \Big\}\\
& \leq 2 c \vecnorm{\big( \IdMat + \regu_\numobs \EigenMat^{-1}
  \big)^{-1 / 2} \plainvec}{\ell^2} \cdot \Big\{ \sqrt{\regu_\numobs}
\rkhsnorm{\treateff} + \varbound \sqrt{\tfrac{\log (1 /
    \delta)}{\numobs}} + \subgaussian \sqrt{\effdim(\regu_\numobs)}
\cdot \tfrac{\log (1 / \delta) \log \numobs}{\numobs} \Big\},
\end{align*}
with probability at least $1 - \delta$.

Substituting the choice $\regu_\numobs = \varbound^2 / (\Rbound^2
\numobs)$, we note that $\sqrt{\regu_\numobs} \rkhsnorm{\treateff}
\leq \sqrt{\regu_\numobs} \Rbound \leq \varbound / \sqrt{\numobs}$,
leading to the bound
\begin{align*}
  \abss{\plainvec^\top (\projparamhat_\numobs - \projparam_*)} \leq 2
  c \vecnorm{\big( \IdMat + \regu_\numobs \EigenMat^{-1} \big)^{-1 / 2}
    \plainvec}{\ell^2} \cdot \Big\{ \varbound \sqrt{\tfrac{\log (1 /
      \delta)}{\numobs}} + \subgaussian \sqrt{\effdim(\regu_\numobs)}
  \cdot \tfrac{\log (1 / \delta) \log \numobs}{\numobs} \Big\}, \quad
  \mbox{w.p. $1- \delta$},
\end{align*}
which completes the proof of~\Cref{lemma:krr-directional-bound}.

\subsection{Proof of~\Cref{ThmHetero} and corollaries}
\label{SecProofThmHetero}

The proof consists of three parts: we first establish guarantees on
the auxiliary estimators $\mutil_\numobs$ and $\sighat^2_\numobs$, and
then use these guarantees to bound the error of the two-stage
estimator $\tauhat_\numobs$. Concretely, we prove the following claims
in turn.
\bcar
\item For any fixed state-action pair $(\state_0, \actionzero)$ and
  any $\delta \in (0, 1)$, the first-stage estimator $\mutil_\numobs$
  satisfies the bound
\begin{subequations}  
  \begin{align}
\label{eq:four-stage-first-part-bound}    
\abss{\mutil_\numobs (\state_0, \actionzero) - \treateff(\state_0,
  \actionzero)} \leq c \subgaussian
\sqrt{\frac{\effdim(\regustageone)}{\numobs} \log (1 / \delta)},
  \end{align}
with probability $1 -
\delta$. See~\Cref{subsubsec:proof-four-stage-first-part-bound} for
the proof.
\item For any fixed state-action pair $(\state_0, \actionzero)$, the
  second-stage estimator $\sighat_\numobs$ satisfies the bound
\begin{align}
\frac{1}{2} \sigma^2(\state_0, \actionzero) \leq
\sighat_\numobs^2(\state_0, \actionzero) \leq 2 \sigma^2(\state_0,
\actionzero),\label{eq:four-stage-second-part-bound}
\end{align}
with probability $1 - \delta /
\numobs$. See~\Cref{subsubsec:proof-four-stage-second-part-bound} for
the proof.
\item Using an approach analogous to that in the proof
  of~\Cref{ThmHomo}, we represent the target functionals
  using basis functions, and recall the error decomposition
  \begin{align}
\label{eq:error-decomp-copied-in-heterosked}    
\tauhat_{\numobs} - \taustar & =
\inprod{\featweighhat_{\numobs}-\widebar{\featweigh}}{\projparam_*} +
\inprod{\widebar{\featweigh}}{\projparamhat_\numobs - \projparam_*} +
\inprod{\featweighhat_\numobs-\widebar{\featweigh}}{\projparamhat_\numobs
  - \projparam_*},
  \end{align}
where we denote $\projparam_* \mydefn \Psi \treateff$ and
$\projparamhat_\numobs \mydefn \Psi \muhat_\numobs$. The errors in the
sample average feature vector $\featweighhat_\numobs$ can be
controlled using~\Cref{lemma:avg-in-rkhs-directional} just as in the
proof of~\Cref{ThmHomo}, while bounding the error for the
weighted least-square estimator $\projparamhat_\numobs$ requires new
ingredients; see~\Cref{lemma:weighted-krr-directional-bound} to
follow.
\end{subequations}
\ecar

\begin{lemma}
\label{lemma:weighted-krr-directional-bound}
Under the conditions of~\Cref{ThmHetero}, with probability
$1 - \delta$, for any infinite-dimensional vector $\plainvec$, we have
\begin{align}
  \abss{\plainvec^\top (\projparamhat_\numobs - \projparam_*)} \leq c
  \vecnorm{\big(\bigcovsig + \regustagethree \EigenMat^{-1} \big)^{-1/2}
    \plainvec}{\ell^2} \sqrt{\frac{\log (1 / \delta)}{\numobs}},
    \end{align}
where $c > 0$ is a universal constant.
\end{lemma}
\noindent See~\Cref{subsubsec:proof-weighted-krr-directional-bound}
for the proof.

\medskip

Having set up the basic ingredients, we are now ready to prove the
main claims of~\Cref{ThmHetero}. We bound each terms in the
decomposition result~\eqref{eq:error-decomp-copied-in-heterosked} as
follows.

Applying the bound~\eqref{eq:avg-directional-bound}
from~\Cref{lemma:avg-in-rkhs-directional} with $f = \treateff$, with
probability $1 - \delta$, we have
\begin{align}
\label{eq:u-noise-bound-in-hetero-proof}  
\abss{\inprod{\featweighhat_{\numobs} -
    \widebar{\featweigh}}{\projparam_*}} \leq 2 \vprobx (\treateff)
\sqrt{\frac{\log (1 / \delta)}{\numobs}} + 6 \kappa \Rbound \frac{\log
  (1 / \delta)}{\numobs}.
\end{align}

Applying~\Cref{lemma:weighted-krr-directional-bound} with $\plainvec =
\featweigh$ yields the bound
\begin{align}
\label{eq:beta-noise-bound-in-hetero-proof}  
  |\inprod{\featweighbar}{\projparamhat_\numobs - \projparam_*}| \leq
  c \vecnorm{\big(\bigcovsig + \regustagethree \EigenMat^{-1}
    \big)^{-1/2} \featweighbar}{\ell^2} \sqrt{\frac{\log
      (1/\delta)}{\numobs}} & \leq 2 c \vwenlonghil \:
  \sqrt{\frac{\log(1/\delta)}{\numobs}}.
\end{align}
valid with probability $1 - \delta$.  Here the second step follows
from \Cref{PropExplicit}.

Finally, applying~\Cref{lemma:weighted-krr-directional-bound}
with $\plainvec = \featweighhat_\numobs -
\featweighbar$,\footnote{Note that the vector $ \featweighhat_\numobs$
is independent of $(\dsthree{\State}_i, \dsthree{\Action}_i,
\dsthree{\outcome}_i)_{i = 1}^\numobs$, so that
\Cref{lemma:weighted-krr-directional-bound} is applicable.} as
well as equation~\eqref{eq:avg-preconditioned-bound}
in~\Cref{lemma:avg-in-rkhs-directional}, with probability $1 -
\delta$, we have the upper bound
\begin{align}
\abss{\inprod{\featweighhat_\numobs -
    \widebar{\featweigh}}{\projparamhat_\numobs - \projparam_*}} &
\leq c \vecnorm{\big(\bigcovsig + \regustagethree \EigenMat^{-1}
  \big)^{-1/2} (\featweighhat_\numobs - \widebar{\featweigh})
}{\ell^2} \sqrt{\frac{\log (1 / \delta)}{\numobs}} \nonumber \\
& \leq c \varbound \vecnorm{\big(\IdMat + \varbound^2 \regustagethree
  \EigenMat^{-1} \big)^{-1/2} (\featweighhat_\numobs -
  \widebar{\featweigh}) }{\ell^2} \sqrt{\frac{\log (1 /
    \delta)}{\numobs}} \nonumber \\
& \leq c \varbound \sqrt{\effdim(\varbound^2 \regustagethree)}
\frac{\log (1/\delta)}{\numobs} \nonumber\\
\label{eq:cross-term-bound-in-hetero-proof}
& = c \varbound \sqrt{\effdim(\regustageone)}
\frac{\log (1/\delta)}{\numobs}.
\end{align}
Combining
equations~\eqref{eq:u-noise-bound-in-hetero-proof},~\eqref{eq:beta-noise-bound-in-hetero-proof},
and~\eqref{eq:cross-term-bound-in-hetero-proof} completes the proof
of~\Cref{ThmHetero}.


\subsubsection{Proof of equation~\eqref{eq:four-stage-first-part-bound}}
\label{subsubsec:proof-four-stage-first-part-bound}

Define the infinite-dimensional vectors
\begin{align*}
\widetilde{\projparam}_\numobs \mydefn \Psi \mutil, \quad \mbox{and}
\quad \projparam_* \mydefn \Psi \treateff.
\end{align*}
The error can be written in the form of the basis function representation
\begin{align*}
    \mutil_\numobs (\state_0, \actionzero) - \treateff(\state_0,
    \actionzero) = \inprod{\widetilde{\projparam}_\numobs -
      \projparam_*}{\phi(\state_0, \actionzero)}.
\end{align*}
Invoking~\Cref{lemma:krr-directional-bound} with $\plainvec =
\phi(\state_0, \actionzero)$, we have
\begin{align*}
     \abss{\inprod{\phi(\state_0, \actionzero)}{ \projparamhat_\numobs
         - \projparam_* }} \leq c \vecnorm{\big( \IdMat +
       \regustageone \EigenMat^{-1} \big)^{-1 / 2} \phi(\state_0,
       \actionzero)}{\ell^2} \cdot \Big\{ \varbound \sqrt{\frac{\log
         (1 / \delta)}{\numobs}} + \subgaussian \phimax
     \sqrt{\effdim(\regustageone)} \cdot \frac{\log (1 / \delta) \log
       \numobs}{\numobs} \Big\},
\end{align*}
holding true with probability $1 - \delta$.

Recall the definition~\eqref{eq:defn-effdim} of effective dimension,
we have the uniform upper bound
\begin{align*}
    \vecnorm{\big( \IdMat + \regustageone \EigenMat^{-1} \big)^{-1 /
        2} \phi(\xzero, \actionzero)}{\ell^2} \leq \sup_{(\state,
      \action)} \vecnorm{\big( \IdMat + \regustageone \EigenMat^{-1}
      \big)^{-1 / 2} \phi(\state, \action)}{\ell^2} \leq \sqrt{\effdim
      (\regustageone)}
\end{align*}
Substituting back and using the sample-size
condition~\eqref{eq:hetero-sample-size-req-first-stage} completes the
proof of equation~\eqref{eq:four-stage-first-part-bound}.


\subsubsection{Proof of equation~\eqref{eq:four-stage-second-part-bound}}
\label{subsubsec:proof-four-stage-second-part-bound}

Define the $\sigma$-field $\datablock_1 \mydefn \sigma \big( \{
\dsone{\State}_i, \dsone{\Action}_i, \dsone{\outcome}_i\}_{i=
  1}^\numobs \big)$. Clearly, the first-stage regression function
$\mutil_\numobs$ is measurable in $\datablock_1$. Since each
data-point $(\dstwo{\State}_i, \dstwo{\Action}_i, \dstwo{\outcome}_i)$
at the second stage is independent of $\datablock_1$,
equation~\eqref{eq:four-stage-first-part-bound} guarantees that
\begin{align*}
    \Prob \left\{ \abss{\mutil_\numobs (\dstwo{\State}_i,
      \dstwo{\Action}_i) - \treateff(\dstwo{\State}_i,
      \dstwo{\Action}_i)} \geq 2 c\subgaussian
    \sqrt{\frac{\effdim(\regustageone)}{\numobs} \log (\numobs /
      \delta)} \right\} \leq \frac{\delta}{2 \numobs^2}.
\end{align*}
Applying union bound and the tower property yields
\begin{align}
 &\Exs \left[ \Prob \Big\{ \max_{i \in [\numobs]} \abss{\mutil_\numobs
      (\dstwo{\State}_i, \dstwo{\Action}_i) -
      \treateff(\dstwo{\State}_i, \dstwo{\Action}_i)} \geq 2c
    \subgaussian \sqrt{\frac{\effdim(\regustageone)}{\numobs} \log
      (\numobs / \delta)} \Big\} \mid \datablock_1 \right]
  \nonumber\\ &\leq \sum_{i = 1}^\numobs \Prob \left\{
  \abss{\mutil_\numobs (\dstwo{\State}_i, \dstwo{\Action}_i) -
    \treateff(\dstwo{\State}_i, \dstwo{\Action}_i)} \geq 2c
  \subgaussian \sqrt{\frac{\effdim(\regustageone)}{\numobs} \log
    (\numobs / \delta)} \right\} \leq \frac{\delta}{2
    \numobs}.\label{eq:union-bound-in-four-stage-second-part-proof}
\end{align}
For the noisy observation we construct in this step, the conditional
expectation takes the form
\begin{align*}
    \Exs \big[ (\dstwo{\outcome}_i - \mutil_\numobs (\dstwo{\State}_i,
      \dstwo{\Action}_i) )^2 \mid \dstwo{\State}_i, \dstwo{\Action}_i,
      \datablock_1 \big] = \sigma^2(\dstwo{\State}_i,
    \dstwo{\Action}_i) + \underbrace{(\mutil_\numobs -
      \treateff)^2(\dstwo{\State}_i, \dstwo{\Action}_i)}_{=:
      \smallbias (\dstwo{\State}_i, \dstwo{\Action}_i)}.
\end{align*}
We further note that for any $p > 0$, we have
\begin{align*}
    &\Exs \Big[ \abss{ \big( \dstwo{\outcome}_i - \mutil_\numobs
      (\dstwo{\State}_i, \dstwo{\Action}_i) \big)^2}^p \mid
    \dstwo{\State}_i, \dstwo{\Action}_i, \datablock_1 \Big] \\ &\leq
  2^{2p} \Exs \Big[ \abss{ \big( \dstwo{\outcome}_i -
      \treateff(\dstwo{\State}_i, \dstwo{\Action}_i) \big)^2}^p \mid
    \dstwo{\State}_i, \dstwo{\Action}_i, \datablock_1 \Big] + 2^{2p}
  \abss{\mutil_\numobs (\dstwo{\State}_i, \dstwo{\Action}_i) -
    \treateff(\dstwo{\State}_i, \dstwo{\Action}_i)}^{2p}\\ &\leq
  4^{2p} p^{p} \subgaussian^{2p} + 2^{2p} \smallbias^p
  (\dstwo{\State}_i, \dstwo{\Action}_i),
\end{align*}
which verifies the tail assumption $\vecnorm{\dstwo{\outcome}_i -
  \mutil_\numobs (\dstwo{\State}_i, \dstwo{\Action}_i)
  \big)^2}{\psi_1} \leq 4 (\subgaussian^2 + \smallbias
(\dstwo{\State}_i, \dstwo{\Action}_i))$ conditionally on
$\dstwo{\State}_i, \dstwo{\Action}_i$ and $\datablock_1$.

Having verified the observation
assumption~\eqref{eq:observation-assumption-for-robust-ptwise-risk},
we are ready to apply the robust pointwise risk property satisfied by
the estimating procedure $\abstrEst$. By definition, given a sample
size $\numobs \geq \maux \big( \lowersigma^2 / 2, \delta / (2 \numobs)
\big)$, for any $\varepsilon > 0$ and $\delta \in (0, 1)$, we have
\begin{align}
     \Prob \Big\{ \abss{\sighat_{\numobs}^2(\state_0, \actionzero) -
       \sigma^2(\state_0, \actionzero)} \leq \lowersigma^2 / 2 \mid
     \datablock_1 \Big\} \geq \Prob \Big\{ \max_{i \in [\numobs]}
     |\smallbias (\dstwo{\State}_i, \dstwo{\Action}_i)| \leq
     \widebar{\smallbias} (\lowersigma^2 / 2, \delta / (2 \numobs))
     \mid \datablock_1 \Big\} - \frac{\delta}{2
       \numobs},\label{eq:robust-ptwise-in-four-stage-second-part-proof}
\end{align}
almost surely.

Given a sample size satisfying the requirement in
equation~\eqref{eq:hetero-sample-size-req-second-stage}, we have
\begin{multline}
   \Prob \Big\{ \max_{i \in [\numobs]} |\smallbias (\dstwo{\State}_i,
   \dstwo{\Action}_i)| \leq \widebar{\smallbias} (\lowersigma^2 / 2,
   \delta / (2\numobs)) \mid \datablock_1 \Big\}\\ \geq \Prob \Big\{
   \max_{i \in [\numobs]} \abss{\mutil_\numobs (\dstwo{\State}_i,
     \dstwo{\Action}_i) - \treateff(\dstwo{\State}_i,
     \dstwo{\Action}_i)} \leq c \subgaussian
   \sqrt{\frac{\effdim(\regustageone)}{\numobs} \log (\numobs /
     \delta)} \mid \datablock_1
   \Big\},\label{eq:max-bias-in-four-stage-second-part-proof}
\end{multline}
almost surely.

Combining
equations~\eqref{eq:union-bound-in-four-stage-second-part-proof},~\eqref{eq:robust-ptwise-in-four-stage-second-part-proof},
and~\eqref{eq:max-bias-in-four-stage-second-part-proof} and taking
expectations with respect to $(\dsone{\State}_i, \dsone{\Action}_i,
\dsone{\outcome}_i)$, we conclude that
\begin{multline*}
    \Prob \Big\{ \abss{\sighat_{\numobs}^2(\state_0, \actionzero) -
      \sigma^2(\state_0, \actionzero)} \leq \lowersigma^2 / 2 \Big\} =
    \Exs \Big[ \Prob \Big\{ \abss{\sighat_{\numobs}^2(\state_0,
        \actionzero) - \sigma^2(\state_0, \actionzero)} \leq
      \lowersigma^2 / 2 \mid \datablock_1 \Big\} \Big]\\ \geq \Exs
    \left[ \Prob \Big\{ \max_{i \in [\numobs]} \abss{\mutil_\numobs
        (\dstwo{\State}_i, \dstwo{\Action}_i) -
        \treateff(\dstwo{\State}_i, \dstwo{\Action}_i)} \geq 2c
      \subgaussian \sqrt{\frac{\effdim(\regustageone)}{\numobs} \log
        (\numobs / \delta)} \Big\} \mid \datablock_1 \right] -
    \frac{\delta}{2 \numobs} \geq 1 - \frac{\delta}{\numobs}.
\end{multline*}
On the event $\abss{\sighat_{\numobs}^2(\state_0, \actionzero) -
  \sigma^2(\state_0, \actionzero)} \leq \lowersigma^2 / 2$, we have
\begin{align*}
    \frac{1}{2} \sigma^2(\state_0, \actionzero) \leq
    \sigma^2(\state_0, \actionzero) - \lowersigma^2 / 2\leq
    \sighat_\numobs^2(\state_0, \actionzero) \leq \sigma^2(\state_0,
    \actionzero) + \lowersigma^2 / 2 \leq 2 \sigma^2(\state_0,
    \actionzero),
\end{align*}
completing the proof of
equation~\eqref{eq:four-stage-second-part-bound}.


\subsubsection{Proof of~\Cref{lemma:weighted-krr-directional-bound}}
\label{subsubsec:proof-weighted-krr-directional-bound}

First, by the guarantee~\eqref{eq:four-stage-second-part-bound} from
the second stage and a union bound, we have
\begin{align}
  \label{eq:second-stage-conclusion-in-weighted-krr-proof}  
  \Prob \Big\{ \exists i \in [\numobs],
  ~\frac{\sighat_\numobs^2}{\sigma^2} (\dsthree{\State}_i,
  \dsthree{\Action}_i) \notin \big(1, 2 \big) \Big\} \leq \sum_{i =
    1}^\numobs \Prob \Big\{ \frac{\sighat_\numobs^2 }{\sigma^2}
  (\dsthree{\State}_i, \dsthree{\Action}_i) \notin \big(1, 2 \big)
  \Big\} \leq \delta.
\end{align}
Defining the event
\begin{align*}
    \dsthree{\Event} \mydefn \Big\{\frac{1}{2}
    \sigma^2(\dsthree{\State}_i, \dsthree{\Action}_i) \leq
    \sighat_\numobs^2 (\dsthree{\State}_i, \dsthree{\Action}_i) \leq
    2\sigma^2(\dsthree{\State}_i, \dsthree{\Action}_i), \quad
    \mbox{for any } i \in [\numobs] \Big\},
\end{align*}
we have $\Prob \big( \dsthree{\Event} \big) \geq 1 - \delta$, with
respect to the randomness of both the state-action pairs
$(\dsthree{\State}_i, \dsthree{\Action}_i)_{i = 1}^\numobs$ and the
function $\sighat_\numobs^2$.

The remainder of the proof is analogous to that
of~\Cref{lemma:krr-directional-bound}.  For notational simplicity, we
omit the supscript $\dsthree{}$ in $(\State_i, \Action_i,
\outcome_i)$.  Under the basis function representation, we have
\begin{align*}
  \projparamhat_\numobs & = \Psi \cdot \arg\min_{\plaintreateff \in
    \rkhs} \Big\{ \frac{1}{\numobs} \sum_{i = 1}^{\numobs}
  \sighat_\numobs^{-2} (\State_i, \Action_i) \big( \outcome_i -
  \plaintreateff(\State_i, \Action_i) \big)^2 + \regustagethree
  \rkhsnorm{\plaintreateff}^2 \Big\}\\
  & = \arg\min_{\projparam \in \ell^2(\naturals)} \Big\{
  \frac{1}{\numobs} \sum_{i = 1}^{\numobs} \sighat_\numobs^{-2}
  (\State_i, \Action_i) \big( \outcome_i -
  \inprod{\projparam}{\phi(\State_i, \Action_i)} \big)^2 +
  \regustagethree \vecnorm{\projparam}{\eigen^{-1}}^2 \Big\}.
\end{align*}
Defining the noise $\noise_i \mydefn \outcome_i - \treateff(\State_i,
\Action_i)$ and the empirical covariance operator
\begin{align*}
\widehat{\bigcov}_{\numobs}^\sigma \mydefn \frac{1}{\numobs} \sum_{i =
  1}^\numobs \sighat_\numobs^{-2} (\State_i, \Action_i) \phi(\State_i,
\Action_i) \phi(\State_i, \Action_i)^\top,
\end{align*}
the error vector admits the representation
\begin{align*}
\projparamhat_\numobs - \projparam_* =
\big(\widehat{\bigcov}_{\numobs}^\sigma + \regustagethree \EigenMat^{-1}
\big)^{-1} \frac{1}{\numobs} \sum_{i = 1}^\numobs \Big\{
\sighat_\numobs^{-2} \noise_i \phi(\State_i, \Action_i) -
\regustagethree \EigenMat^{-1} \projparam_* \Big\}
\end{align*}
We can bound such an error conditionally on the state-action pairs
$(\State_i, \Action_i)_{i = 1}^\numobs$ and the estimated conditional
covariance function $\sighat_\numobs$, as stated in the following
lemma.
\begin{lemma}
\label{lemma:weighted-krr-conditional-bound}
Under the set-up above, conditionally on the function
$\sighat_\numobs$ and the state-action pairs $(\State_i, \Action_i)_{i
  = 1}^\numobs$ such that the event $\dsthree{\Event}$ happens, with
probability $1 - \delta$, we have the upper bound
\begin{subequations}
\begin{multline}
  \abss{\plainvec^\top (\projparamhat_\numobs - \projparam_*)} \leq c
  \vecnorm{\big(\widehat{\bigcov}_{\numobs}^\sigma + \regu_{\numobs,
      *} \EigenMat^{-1} \big)^{-1/2} \plainvec}{\ell^2}\\ \times
  \left\{ 2 \sqrt{\frac{\log (1 / \delta)}{\numobs}} + \frac{\log
    \numobs \log (1 / \delta) \subgaussian }{\numobs \lowersigma^2 }
  \sup_{(\state, \action)}
  \vecnorm{\big(\widehat{\bigcov}_{\numobs}^\sigma + \regu_{\numobs,
      *} \EigenMat^{-1} \big)^{-1 / 2} \phi(\state, \action)}{\ell^2}
  \right\}
\end{multline}
\end{subequations}
\end{lemma}
\noindent See~\Cref{app:subsec-weighted-krr-conditional} for the
proof.

\medskip

Taking this lemma as given, we proceed with the proof
of~\Cref{lemma:weighted-krr-directional-bound}.  Define the truncated
variance function and the corresponding reweighted operator.
\begin{align*}
\widetilde{\sigma}_\numobs (\state, \action) & \mydefn \begin{cases}
  \sighat_\numobs^2(\state, \action) & \quad \mbox{if
    $\frac{\sighat_\numobs^2(\state, \action)}{\sigma^2(\state,
      \action)} \in \big(1/2, 2 \big)$}, \\
\sigma^2(\state, \action) & \mbox{otherwise}
\end{cases} \quad \mbox{and} \\
\widetilde{\bigcov}_{\numobs}^\sigma & \mydefn \frac{1}{\numobs}
\sum_{i = 1}^\numobs \frac{1}{\widetilde{\sigma}_\numobs^{2}
  (\State_i, \Action_i)} \phi(\State_i, \Action_i) \phi(\State_i,
\Action_i)^\top.
\end{align*}
Conditioned on the event $\dsthree{\Event}$, we have
$\widetilde{\bigcov}_\numobs^\sigma =
\widehat{\bigcov}_\numobs^\sigma$. On the other hand, we invoke
\Cref{lemma:emp-cov-concentration} with the weight function $q =
\widetilde{\sigma}_\numobs^{-2}$ and $\offpar = 1/2$. Note that the
condition~\eqref{eq:regu-requirement-in-emp-cov-concentration-general}
becomes
\begin{align*}
     (\varbound^2 / \lowersigma^2)\log \big( \frac{\kappa^2
  }{\regustagethree \delta} \big) \cdot \frac{\effdim
    \big(\regustagethree \varbound^2 \big)}{\numobs} \leq
  \frac{1}{32},
\end{align*}
which is satisfied under the sample size
requirement~\eqref{eq:hetero-sample-size-req-first-stage} and
regularization parameter
choice~\eqref{eq:optimal-lambda-four-stage}. Therefore, on the event
$\dsthree{\Event}$, with probability $1 - \delta$, we have
\begin{align}
\label{eq:bigcov-sigma-conc-in-weighted-krr-directional}  
\frac{1}{2} \big(\bigcovsig + \regustagethree \EigenMat^{-1}\big) \preceq
\widetilde{\bigcov}_\numobs^\sigma + \regustagethree \EigenMat^{-1} =
\widehat{\bigcov}_\numobs^\sigma + \regustagethree \EigenMat^{-1} \preceq
2 \big(\bigcovsig + \regustagethree \EigenMat^{-1}\big).
\end{align}
Substituting
equation~\eqref{eq:bigcov-sigma-conc-in-weighted-krr-directional} into
the guarantee from
~\Cref{lemma:conditional-high-prob-bound-for-krr-directional}, we
find that
\begin{multline*}
    \abss{\inprod{\plainvec}{\projparamhat_\numobs - \projparam_*}}
    \leq 4 c \vecnorm{\big(\bigcovsig + \regustagethree
      \EigenMat^{-1} \big)^{-1/2} \plainvec}{\ell^2} \cdot \Biggr \{
    \sqrt{\frac{\log (1 / \delta)}{\numobs}} \\
    + \frac{\log \numobs \log
      (1 / \delta) \subgaussian }{\numobs \lowersigma^2 }
    \sup_{(\state, \action)} \vecnorm{\big(\bigcovsig +
      \regustagethree \EigenMat^{-1} \big)^{-1 / 2} \phi(\state,
      \action)}{\ell^2} \Biggr\}.
\end{multline*}
By the definition~\eqref{eq:defn-effdim} of effective dimension, we
have
\begin{align*}
  \sup_{ (\state, \action)} \vecnorm{\big(\bigcovsig + \regustagethree
    \EigenMat^{-1} \big)^{-1 / 2} \phi(\state, \action)}{\ell^2} \leq
  \varbound \sup_{(\state, \action)} \vecnorm{\big( \IdMat +
    \varbound^2 \regustagethree \EigenMat^{-1} \big)^{-1 / 2}
    \phi(\state, \action)}{\ell^2} = \varbound
  \sqrt{\effdim(\regustageone)}.
\end{align*}
Moreover, given the 
condition~\eqref{eq:hetero-sample-size-req-first-stage} on the
sample size, it follows that
\begin{align*}
\sqrt{\frac{\log (1 / \delta)}{\numobs}} \geq \frac{\log \numobs \log
  (1 / \delta) \subgaussian^2}{\numobs \lowersigma^2}
\sqrt{\effdim(\regustageone)} \geq \frac{\log \numobs \log (1 /
  \delta) \subgaussian }{\numobs \lowersigma^2 } \sup_{(\state,
  \action)} \vecnorm{\big(\bigcovsig + \regustagethree
  \EigenMat^{-1} \big)^{-1 / 2} \phi(\state, \action)}{\ell^2}.
\end{align*}
Thus, we conclude that
\begin{align*}
\abss{ \inprod{\plainvec}{\projparamhat_\numobs - \projparam_*}} \leq
8 c \vecnorm{\big(\bigcovsig + \regustagethree \EigenMat^{-1}
  \big)^{-1/2} \plainvec}{\ell^2} \sqrt{\frac{\log (1 /
    \delta)}{\numobs}}
\end{align*}
with probability $1 - \delta$, with establishes the claim
in~\Cref{lemma:weighted-krr-directional-bound}.


\subsubsection{Proof of equation~\eqref{EqnCATEclaim}}\label{subsubsec:proof-cor-cate}

Recall the definition~\eqref{eq:defn-featureweightbar-general} of the
infinite-dimensional vector $\featweighbar (\delta_\xzero)$
\begin{align*}
	\featweighbar (\delta_\xzero) = \int_\ActionSpace \phi (\xzero,
    \action) d \wenweight{\action}{\xzero}.
\end{align*}
Using this definition, the estimation error admits a basis-function
representation
\begin{align*}
\tauhat_{ \numobs} (\xzero) - \taustar (\xzero) =
\inprod{\projparamhat_\numobs - \projparam_*}{ \widetilde{\featweigh}
  (\xzero) },
\end{align*}
where the vectors $\projparamhat_\numobs$ and $\projparam_*$ are
defined in~\Cref{subsubsec:setup-for-aux-results}.

Applying~\Cref{lemma:weighted-krr-directional-bound} with $\plainvec =
\widetilde{\featweigh} (\xzero)$ yields
\begin{align}
\label{eq:apply-weighted-krr-direc-lemma-to-cate}  
\abss{ \inprod{\widetilde{\featweigh} (\xzero)}{\projparamhat_\numobs
    - \projparam_*}} \leq c \vecnorm{\big(\bigcovsig + \regustagethree
  \EigenMat^{-1} \big)^{-1/2} \widetilde{\featweigh} (\xzero)}{\ell^2}
\sqrt{\frac{\log (1 / \delta)}{\numobs}}.
\end{align}
By~\Cref{PropExplicit}, we have
\begin{align*}
	\vecnorm{\big(\bigcovsig + \regustagethree \EigenMat^{-1}
          \big)^{-1/2} \widetilde{\featweigh} (\xzero)}{\ell^2}^2 =
        \widetilde{\featweigh}^\top (\xzero) \big(\bigcovsig +
        \tfrac{1}{\Rbound^2 \numobs} \EigenMat^{-1} \big)^{-1}
        \widetilde{\featweigh} \leq 4 \vsupersq{\delta_{\xzero}}.
\end{align*}

Substituting back completes the proof of the
claim~\eqref{EqnCATEclaim}.


\section{Discussion}
\label{SecDisussion}

In this paper, we studied the problem of estimating linear functionals
based on observational data.  Our main focus was the challenging
setting in which the importance ratio is poorly behaved.  In such
settings, the classical semi-parametric efficiency bound---based on a
presumptive $\sqrt{\numobs}$-rate of convergence---can be infinite,
and so fail to characterize the problem.  So as to remedy this
deficiency, the main contribution of this paper was to propose a
modified risk functional, defined as the optimal value of a
variational problem that respects the geometry of the function
class. The resulting minimax risks interpolate between the classical
regimes of semi-parametric efficiency with the $\sqrt{\numobs}$-rate,
and nonparametric rates for functional estimation. Focusing on the
case of RKHS, we analyze an outcome-based regression estimator, and
showed that it achieves our instance-dependent lower bound (up to a
universal constant pre-factor).  This estimator is attractive in not
requiring any knowledge of the behavioral policy.  Nonetheless,
despite its agnostic nature, it matches our lower bound that applies
even to oracle estimators that have full knowledge of the policy. When
applied to various off-policy estimation problems with singularities
in the importance ratio, our results uncover a novel class of minimax
rates, as well as instance-dependent optimality, adaptively achieved
by our estimators.

While this paper takes an initial step in characterizing
instance-dependent optimality for off-policy estimation beyond
semi-parametric efficiency, there are many open directions.

\bcar
\item Our optimality results impose assumptions on the conditional
  variance function $\sigma^2$.  We either require it to be uniformly
  bounded (for achieving the worst-case variance bound
  $\vclassical^2(\smallimpratio)$), or require additional structure that allows
  for consistent estimation (for optimal adaptation to the conditional
  variance structure). It is not clear if such requirements are
  necessary. In the classical $\sqrt{n}$-regime of semi-parametric
  efficiency, regime, AIPW estimators adapt to the conditional
  variance structure without knowledge of $\sigma^2$; for instance,
  see the paper~\cite{chernozhukov2018double}. An important open
  question, therefore, is whether such adaptivity is possible in the
  more challenging regime considered by our paper without additional
  assumptions on the conditional variance.
\item In this paper, we established achievability of our lower bounds
  only for reproducing kernel Hilbert spaces.  Thus, an important
  question is to what extent our results can be extended to more
  general function classes.  In the special case of homoskedastic
  models, we conjecture that the minimax linear estimation
  strategy~\cite{donoho1994statistical,hirshberg2021augmented} could
  yield an optimal estimator---in the same sense as the results
  presented in this paper---for any function class $\funcClass$
  satisfying the Donsker property. For non-Donsker classes, it is
  known from past work~\cite{robins2009semiparametric} that knowledge
  of the behavior policy plays a role. An important direction of
  future research, therefore, is to identify the optimal risk for
  estimation, jointly determined by the structural assumptions on the
  treatment effect function, the behavior policy function, and
  singularities in the importance ratio function.
\item Our results focus on the classical off-policy contextual bandit
  setup, where the data $(\State_i, \Action_i, \outcome_i)_{i =
    1}^\numobs$ are independent and identically distributed. However,
  many decision-making problems involve collecting data in an adaptive
  manner (e.g., by running a bandit algorithm), or following a Markov
  chain (e.g., in reinforcement learning). Such settings can lead to
  poorly controlled importance ratios, which in turn induces practical
  challenges~\cite{yin2020asymptotically,khamaru2021near}. Our
  estimation framework and risk functional, being optimally agnostic
  to the singularity in the importance ratio, could well be helpful
  for problems of this type.
  \ecar

\section*{Acknowledgement}
The authors thank Fangzhou Su for helpful discussion.
We gratefully acknowledge the support of the NSF through NSF-CCF grant 1955450 and NSF-DMS grant 2015454 to MJW, NSF-IIS grant 1909365 and NSF grant DMS-2023505 to MJW and PLB, and of the ONR through MURI award N000142112431 to PLB.

\AtNextBibliography{\small} \printbibliography


\appendix

\section{Properties of effective dimension}
\label{AppSupNorm}
In this section, we develop various bounds on the effective
dimension under decay rates on the eigenvalues, along with some
regularity conditions on the eigenfunctions.

\subsection{Regularity conditions on eigenfunctions}

The most straightforward assumption on the eigenfunctions is
the \emph{uniform boundedness condition}
\begin{align}
  \label{EqnUniformBound}
  \|\phi_j\|_\infty = \sup_{(\state, \action)} |\phi_j(\state,
  \action)| < \infty \qquad \mbox{for all $j = 1, 2, \ldots.$}
\end{align}
This condition appears frequently in the 
literature~\cite{williamson2001generalization,mendelson2010regularization,duan2021optimal,nie2021quasi}, but as noted, it is not satisfied
by all kernels. See
paper~\cite{zhou2002covering} and Appendix~\ref{app:subsec-proof-missing-data-rates} for some natural counterexamples.

In this paper, we consider the following relaxed growth condition:
there exists a scalar $\smallexponent \in [0, 1/2)$, such that the
  sup-norm of eigenfunctions satisfy the bound
  \begin{align}
\tag{Eig$(\smallexponent)$}
  \label{EqnEigfunGrowth}
  \phimax \mydefn \sup_{j \geq 1} \sup_{(\state, \action)}
  \eigen_j^{\smallexponent} |\phi_j(\state, \action)| < \infty.
\end{align}

\noindent We note that the requirement $\smallexponent \in [0, 1/2)$
  is natural, since the kernel boundedness condition~\eqref{EqnKerBou}
  implies that condition~\eqref{EqnEigfunGrowth} holds with
  $\smallexponent = 1/2$ and $\phimax = \kappa$. An exponent
  $\smallexponent$ strictly less than $1/2$ guarantees slightly more
  regularity.  The growth condition~\eqref{EqnEigfunGrowth} with
  $\smallexponent = 0$ is equivalent to the uniform boundedness
  condition~\eqref{EqnUniformBound}.  However, when $\smallexponent >
  0$, the relaxed condition, on the other hand, is much weaker. For
  example, it is shown by Mendelson and
  Neeman~\cite{mendelson2010regularization} that the counterexample in
  the paper~\cite{zhou2002covering} satisfies
  equation~\eqref{EqnEigfunGrowth} for any $\smallexponent > 0$.

Under Assumption~\eqref{EqnEigfunGrowth}, we have the upper bound
\begin{align}
\label{eq:effdim-upper-bound-under-eig}  
\effdim(\regu) \leq \sum_{j = 1}^{\infty} \sup_{(\state, \action)}
\frac{\eigen_j \phi_j^2 (\state, \action)}{\eigen_j + \regu}\leq
\phimax^2 \sum_{j = 1}^{\infty} \frac{\eigen_j^{1 - 2
    \smallexponent}}{\eigen_j + \regu}.
\end{align}
This upper bound, when combined with decay conditions on the
eigenvalue sequence $\{\eigen_j \}_{j=1}^\infty$, allows us to derive
explicit bounds on the effective dimension.  Two natural classes of
eigenvalue decay are the polynomial condition
\begin{subequations}
\begin{align}
\label{EqnPolyDecay}  
\eigen_j \leq \eigen_0 j^{- \alpha} \qquad \mbox{for some $\alpha > 1$,}
\end{align}
and the \emph{exponential decay}
\begin{align}
  \label{EqnExpDecay}
\eigen_j \leq \eigen_0 \exp (- c_0 j) \qquad \mbox{for some $c_0 >
    0$.}
\end{align}
\end{subequations}
\begin{proposition}
\label{prop:effective-dim}
  Under Assumptions~\eqref{EqnKerBou} and~\eqref{EqnEigfunGrowth}, we
  have
\begin{enumerate}
  \item[(a)] For eigenvalues with
    $\alpha$-polynomial-decay~\eqref{EqnPolyDecay} for some $\alpha >
    \frac{1}{1 - 2 \smallexponent}$, we have
    \begin{subequations}
    \begin{align}
      \effdim(\regu) & \leq c \regu^{- \frac{1}{\alpha} -
        2\smallexponent}.
    \end{align}
  \item[(b)] For eigenvalues with exponential
    decay~\eqref{EqnExpDecay}, we have
    \begin{align}
      \effdim(\regu) \leq c \log \big( \tfrac{\eigen_0}{\regu} \big).
    \end{align}
    \end{subequations}
\end{enumerate}
\end{proposition}
\noindent In these bounds, the constant $c$ can depend on problem
parameters $(\eigen_0, \alpha, c_0, \smallexponent)$ but is
independent of $\regu$.  See~\Cref{SecProofEffDim} for the proof.

In our main theorems, the bounds on the effective dimension is used to
establish the sample size
requirement~\eqref{eq:sample-size-req-based-on-eigendecay}
and~\eqref{eq:hetero-sample-size-req-first-stage}. In order for them
to be true, up to logarithmic factors of $(\numobs, \subgaussian,
\lowersigma^{-1}, \Rbound)$, we need sample sizes
\begin{align*}
\numobs & \gtrsim \big( \subgaussian / \lowersigma \big)^{\frac{2
    \alpha}{\alpha - 1 - 2 \smallexponent \alpha}} \cdot \big(
\Rbound/\varbound \big)^{\frac{1 + 2 \alpha \smallexponent}{\alpha - 1
    - 2 \smallexponent \alpha}}, \quad
\mbox{under~\Cref{prop:effective-dim}(a)},\\
\numobs &\gtrsim \big( \subgaussian/\lowersigma \big)^2, \quad
\mbox{under~\Cref{prop:effective-dim}(b)}.
\end{align*}
In words, the sample size requirement depends on two important
objects: the tail conditions of the observation noise $W = \outcome -
\treateff (\State, \Action)$, measured by the ratio between its
largest Orlicz norm and smallest variance; and the richness of the
kernel class, measured by the eigenvalue decay rates and the radius
of the RKHS ball.


\subsection{Proof of~\Cref{prop:effective-dim}}
\label{SecProofEffDim}

Since the eigenvalue sequence converges to zero, the cut-off integer
$J \mydefn \sup \{ j \geq 1 \mid \eigen_j > \regu\}$ is guaranteed to
be finite.  By the definition of the effective dimension, we have
\begin{align*}
\effdim(\regu) = \sum_{j = 1}^{\infty} \frac{\eigen_j^{1 - 2
    \smallexponent}}{\eigen_j + \regu} \leq \sum_{j \leq J} \eigen^{-
  2\smallexponent} + \regu^{-1} \sum_{j > J} \eigen_j^{1 - 2
  \smallexponent} \leq \regu^{- 2 \smallexponent} J + \regu^{-1}
\sum_{j > J} \eigen_j^{1 - 2 \smallexponent}.
\end{align*}
We prove the results for two cases separately.

For the polynomially-decaying eigenvalues, we have $ J \leq \big(
\frac{\eigen_0}{\regu} \big)^{1/\alpha}$, and
\begin{align*}
\sum_{j > J} \eigen_j^{1 - 2 \smallexponent} \leq \eigen_0 \sum_{j >
  J} j^{- \alpha (1 - 2 \smallexponent)} < \frac{1}{\alpha (1 - 2
  \smallexponent) - 1} J^{1 - \alpha (1 - 2 \smallexponent)}.
\end{align*}
Combining these bounds yields
\begin{align*}
D(\regu) \leq \eigen_0^{1/\alpha} \regu^{ - 2\smallexponent -
  1/\alpha} + \frac{1}{\alpha (1 - 2 \smallexponent) - 1}
\eigen_0^{1/\alpha - 1 + 2 \smallexponent} \regu^{ - 1/\alpha - 2
  \smallexponent } \leq c (\eigen_0, \smallexponent, \alpha) \regu^{-
  1 / \alpha - 2 \smallexponent}.
\end{align*}

For exponentially-decaying eigenvalues, we have $J \leq c_0^{-1} \log
\big( \frac{\eigen_0}{\regu} \big)$, and
\begin{align*}
 \sum_{j > J} \eigen_j^{1 - 2 \smallexponent} \leq \eigen_0 \sum_{j >
   J} \exp \big( - c_0 (1 - 2 \smallexponent) j \big) \leq
 \frac{\eigen_0}{(1 - 2 \smallexponent)c_0} \exp (- c_0 (1 - 2
 \smallexponent) J) \leq \frac{\regu}{c_0},
\end{align*}
which leads to the effective dimension bound
\begin{align*}
\effdim(\regu) \leq c(c_0, \eigen_0, \smallexponent) \log
\big(\frac{\eigen_0}{\regu} \big),
\end{align*}
completing the proof of~\Cref{prop:effective-dim}.


\section{Relaxing the effective dimension condition}
\label{subsubsec:near-optimal-rates}

Recall that~\Cref{ThmHomo} requires certain growth conditions on the
effective dimension.  In this section, we discuss how these conditions
can be relaxed, thereby obtaining a bound that remains
instance-optimal up to logarithmic factors.

\subsection{Near-optimal rates}

Our result involves the modified regularization parameter
\begin{align}
\label{eq:regu-choice-in-weak-cor}  
\regu_\numobs = \frac{\varbound^2}{\Rbound^2 \numobs} \vee \frac{32
  \kappa^2 \log (\numobs / \delta)}{\numobs},
\end{align}
along with the modified higher-order term
$\highorder_\numobs' \mydefn \big( \subgaussian + \kappa \Rbound \big)
\frac{\log (\numobs/ \delta)}{\sqrt{\numobs}}$.

\begin{corollary}
\label{cor:optimal-rate-weak-assumptions}
Suppose that Assumptions~\ref{EqnKerBou} and~\ref{EqnCleanSubGauss}
are in force, and we that implement the method with regularization
parameter~\eqref{eq:regu-choice-in-weak-cor}.  Then for for any sample
size $\numobs \geq 2$ and any $\delta \in (0, 1)$, we have
\begin{align}
\abss{\tauhat_\numobs - \taustar} \leq c v_\probx (\treateff)
\sqrt{\frac{\log (1 / \delta)}{\numobs} } + \frac{c (\subgaussian +
  \kappa \Rbound)}{\varbound} \vbar_\sigma \big( \ball_{\rkhs}
(\Rbound); \numobs \big) \frac{\log (\numobs / \delta)}{\sqrt{\numobs}
} + c \highorder_\numobs'
\end{align}
with probability at least $1 - \delta$.
\end{corollary}
\noindent See~\Cref{subsec:proof-corollary-optimal-rate} for the proof.

\medskip

A few remarks are in order. First,
~\Cref{cor:optimal-rate-weak-assumptions} holds for \emph{any} sample
size, and is completely agnostic to conditions on the effective
dimension $\effdim$. Compared to the optimal instance-dependent bounds
in~\Cref{ThmHomo}, ~\Cref{cor:optimal-rate-weak-assumptions} exhibits
two differences:
\begin{itemize}
\item The variance functional $\vwenlonghilbar$ is multiplied with a
  problem-dependent factor $\frac{\subgaussian + \kappa
    \Rbound}{\varbound}$, as well as logarithmic factors in the ratio
  $\numobs / \delta$.
\item The high-order term $\highorder_\numobs'$ is of order $\Theta
  (\numobs^{-1/2})$, with additional logarithmic factors. Such a
  convergence rate is slower than the high-order term bound
  $\highorder_\numobs$ established in~\Cref{ThmHomo}, which decays at
  a rate $o (\numobs^{-1/2})$
\end{itemize}
Due to these two major differences, the bound
in~\Cref{cor:optimal-rate-weak-assumptions} may not be always
instance-optimal. However, we remark that the near-optimal rate of
convergence (as a function of sample size $\numobs$) is still
preserved.  In the high-noise regime where the quantities
$(\subgaussian, \varbound, \kappa \Rbound)$ are of the same order, the
leading-order terms in~\Cref{ThmHomo}
and~\Cref{cor:optimal-rate-weak-assumptions} differ only by
logarithmic factors. The term $\highorder'$ is dominated by the
leading-order one, up to logarithmic factors. In combination, results
in~\Cref{cor:optimal-rate-weak-assumptions} under the weak assumptions
can be worse than~\Cref{ThmHomo} only by logarithmic factors and
problem dependent constants. Note that a variety of convergence rates
can be established beyond the classical $\sqrt{\numobs}$-regime (see
\Cref{SecExamples} for concrete examples). These convergence rates,
though depending on the intricate properties of the policy
$\propscore$, are automatically achieved without the effective
dimension condition.


\subsection{Proof of~\Cref{cor:optimal-rate-weak-assumptions}}
\label{subsec:proof-corollary-optimal-rate}

We use the same notation $(\featweighhat_\numobs, \projparam_\numobs,
\projparam_*)$ as in the proof of~\Cref{ThmHomo}. Recall from the
decomposition~\eqref{eq:linear-regression-error-decomp} that
$\tauhat_{ \numobs} - \taustar = \inprod{\featweighhat_{\numobs} -
  \widebar{\featweigh}}{\projparam_*} +
\inprod{\widebar{\featweigh}}{\projparamhat_\numobs - \projparam_*} +
\inprod{\featweighhat_\numobs - \widebar{\featweigh}
}{\projparamhat_\numobs - \projparam_*}$.  Since the
bound~\eqref{eq:ubar-beta-noise-bound} does not rely on the effective
dimension, it still holds under our current assumptions---that is, we
have
\begin{align}
\label{eq:ubar-noise-bound-copied}
\abss{ \inprod{\featweighhat_\numobs -
    \widebar{\featweigh}}{\projparam_*}} \leq 2 \vprobx (\treateff)
\sqrt{\tfrac{\log (1 / \delta)}{\numobs}} + 6 \kappa \Rbound \tfrac{\log
  (1 / \delta)}{\numobs},
\end{align}
with probability $1 - \delta$.

The rest of this section is devoted to the control of the other two
terms in the decomposition. We use the following lemma, which is
analogous to~\Cref{lemma:krr-directional-bound}.

\begin{lemma}
\label{lemma:krr-directional-bound-weak-assumption}
Uner the assumptions of ~\Cref{cor:optimal-rate-weak-assumptions}, for
any fixed $z \in \ell^2$, we have
\begin{align*}
\abss{\inprod{\plainvec}{ \projparamhat_\numobs - \projparam_*}} &
\leq c \vecnorm{\big( \IdMat + \regu_\numobs \EigenMat^{-1} \big)^{-1
    / 2} \plainvec}{\ell^2} \cdot (\subgaussian + \kappa \Rbound
)\sqrt{\frac{\log (\numobs / \delta) \log (1 / \delta)}{\numobs}}.
\end{align*}
with probability at least $1 - \delta$.
\end{lemma}
\noindent
See~\Cref{app:subsec-proof-krr-directional-bound-weak-assumption} for
the proof.

\medskip

Taking this lemma as given, we proceed with the proof
of~\Cref{cor:optimal-rate-weak-assumptions}.
Applying~\Cref{lemma:krr-directional-bound-weak-assumption} with
$\plainvec = \featweighbar$, we have
\begin{align}
  \abss{\inprod{\featweighbar}{ \projparamhat_\numobs - \projparam_*}}
  & \leq c \vecnorm{\big( \IdMat + \regu_\numobs \EigenMat^{-1} \big)^{-1 /
      2} \featweighbar}{\ell^2} \cdot (\subgaussian + \kappa \Rbound )
  \frac{\log (\numobs / \delta)}{\sqrt{\numobs}} \nonumber\\
& \leq c \vecnorm{\big( \IdMat + \frac{\varbound^2}{\Rbound^2 \numobs}
    \EigenMat^{-1} \big)^{-1 / 2} \plainvec}{\ell^2} \cdot (\subgaussian +
  \kappa \Rbound )\sqrt{\frac{\log (\numobs / \delta) \log (1 /
      \delta)}{\numobs}} \nonumber\\
\label{eq:betastar-noise-bound-weak}  
& \leq 4 c \vbar_\sigma \big( \ball_\rkhs (\Rbound); \numobs \big)
\cdot \frac{\subgaussian + \kappa \Rbound}{\varbound} \cdot \frac{\log
  (\numobs / \delta)}{\sqrt{\numobs}},
\end{align}
with probability at least $1 - \delta$.

In our next step, we apply
equation~\eqref{eq:avg-preconditioned-bound}
from~\Cref{lemma:avg-in-rkhs-directional}, as well as
condition~\eqref{EqnKerBou}. Doing so yields
\begin{align*}
\vecnorm{\big( \IdMat + \regu_\numobs \EigenMat^{-1} \big)^{-1 / 2}
  (\featweighhat_\numobs - \featweighbar)}{\ell^2} \leq c
\sqrt{\frac{\effdim (\regu_\numobs)}{\numobs}\log (1 / \delta)} \leq c \kappa
\sqrt{\frac{\log (1 / \delta)}{\regu_\numobs \numobs}} \leq c,
\end{align*}
with probability $1 - \delta$.

Combining with~\Cref{lemma:krr-directional-bound-weak-assumption}
yields
\begin{align}
\abss{\inprod{\featweighhat_\numobs - \featweighbar}{
    \projparamhat_\numobs - \projparam_* }} & \leq c \vecnorm{\big(
  \IdMat + \regu_\numobs \EigenMat^{-1} \big)^{-1 / 2}
  (\featweighhat_\numobs - \featweighbar)}{\ell^2} \cdot (\subgaussian
+ \kappa \Rbound )\frac{\log (\numobs / \delta)}{\sqrt{\numobs}}
\nonumber \\
\label{eq:cross-noise-bound-weak}
& \leq c' (\subgaussian + \kappa \Rbound ) \frac{\log (\numobs /
  \delta)}{\sqrt{\numobs}}.
\end{align}
Combining equations~\eqref{eq:ubar-noise-bound-copied},
~\eqref{eq:betastar-noise-bound-weak}
and~\eqref{eq:cross-noise-bound-weak} completes the proof
of~\Cref{cor:optimal-rate-weak-assumptions}.


\section{Proof of technical lemmas}
\label{SecTechnicalHomo}

We collect the proofs of auxiliary lemmas in the proof
of~\Cref{ThmHomo} in this section.

\subsection{Proof of~\Cref{lemma:conditional-high-prob-bound-for-krr-directional}}
\label{app:subsec-proof-conditional-high-prob-bound-for-krr-directional}

We start with the decomposition
\begin{align*}
  \plainvec^\top \big( \projparamhat_\numobs - \projparam_* \big) & =
  \myunder{\frac{1}{\numobs} \sum_{i=1}^\numobs W_i}{noise part} -
  \myunder{\reguhomo \plainvec^\top \EigenMat^{-1}
    \projparam_*}{bias part},
\end{align*}
where $W_i \defn \noise_i \, \plainvec^\top\big(
\widehat{\bigcov}_\numobs + \reguhomo \EigenMat^{-1} \big)^{-1}
\phi(\State_i, \Action_i)$.

Beginning with the bias term, we note that
\begin{align}
  \abss{ \plainvec^\top \big( \widehat{\bigcov}_\numobs +
    \reguhomo \EigenMat^{-1} \big)^{-1} \cdot \reguhomo
    \EigenMat^{-1} \projparam_*} & = \reguhomo \cdot
  \abss{\inprod{\EigenMat^{-1/2} \big( \widehat{\bigcov}_\numobs +
      \reguhomo \EigenMat^{-1} \big)^{-1} \plainvec
    }{\EigenMat^{-1/2} \projparam_*}} \notag \\
& \leq \reguhomo \rkhsnorm{\treateff} \cdot
  \vecnorm{\EigenMat^{-1/2} \big( \widehat{\bigcov}_\numobs +
    \reguhomo \EigenMat^{-1} \big)^{-1} \plainvec}{\ell^2}
  \nonumber \\
\label{eq:bias-part-in-krr-directional-bound-conditional}
& \leq \sqrt{\reguhomo} \rkhsnorm{\treateff} \cdot \vecnorm{\big(
  \widehat{\bigcov}_\numobs + \reguhomo \EigenMat^{-1} \big)^{-1 /
    2} \plainvec}{\ell^2}.
\end{align}
Here the final step is based on the fact that $\vecnorm{\newAmat
  x}{\ell^2} \leq \vecnorm{\Bmat x}{\ell^2}$ for any pair $(\newAmat,
\Bmat)$ of operators such that $\newAmat \preceq \Bmat$.

For the stochastic part, we note that the noise variables
$\{\noise_i\}_{i = 1}^\numobs$ are independent conditioned on
$(\State_i, \Action_i)_{i = 1}^\numobs$. For each $i \in [\numobs]$,
the conditional variance  takes the form
\begin{align*}
    \var \Big(W_i \mid (\State_i, \Action_i)_{i = 1}^\numobs \Big) =
    \sigma^2(\State_j, \Action_j) \cdot \abss{ \phi(\State_i,
      \Action_i)^\top (\widehat{\bigcov}_\numobs + \reguhomo
      \EigenMat^{-1} \big)^{-1} \plainvec}^2.
\end{align*}
Summing up these relations yields
\begin{align*}
\var \Big ( \frac{1}{\numobs} \sum_{i = 1}^\numobs W_i \mid (\State_i,
\Action_i)_{i = 1}^\numobs \Big) & = \plainvec^\top
(\widehat{\bigcov}_\numobs + \reguhomo \EigenMat^{-1} \big)^{-1}
\sum_{i = 1}^\numobs \sigma^2(\State_j, \Action_j) \phi(\State_i,
\Action_i) \phi(\State_i, \Action_i)^\top (\widehat{\bigcov}_\numobs +
\reguhomo \EigenMat^{-1} \big)^{-1} \plainvec \\
& \leq \varbound^2 \vecnorm{\widehat{\bigcov}_\numobs^{1/2} \big(
    \widehat{\bigcov}_\numobs + \reguhomo \EigenMat^{-1} \big)^{-1}
    \plainvec}{2}^2.
\end{align*}

Introducing the shorthand $\Mmat = \widehat{\bigcov}_\numobs +
\reguhomo \EigenMat^{-1}$, by the noise tail
assumption~\eqref{EqnCleanSubGauss} and Adamczak's concentration
inequality~\cite{adamczak2008tail}, conditionally on $(\State_i,
\Action_i)_{i = 1}^\numobs$, we have
\begin{multline*}
\abss{ \plainvec^\top \Mmat^{-1} \cdot \tfrac{1}{\numobs} \sum_{i =
    1}^\numobs \big\{ \noise_i \phi(\State_i, \Action_i) \big\} } \\
\leq c \varbound \vecnorm{\widehat{\bigcov}_\numobs^{1/2} \Mmat^{-1}
  \plainvec}{\ell^2} \sqrt{\tfrac{\log (1 / \delta)}{\numobs}} +
\max_{i \in [\numobs]} \abss{\plainvec^\top \Mmat^{-1} \phi(\State_i,
  \Action_i)} \tfrac{\subgaussian \log \numobs \log (1 /
  \delta)}{\numobs},
\end{multline*}
with probability $1 - \delta$.

In order to control the max term on the RHS, we invoke the
Cauchy--Schwarz inequality, thereby finding that
\begin{align*}
\abss{\plainvec^\top \Mmat^{-1} \phi(\State_i, \Action_i)} &\leq
\vecnorm{\Mmat^{-1/2} \plainvec}{\ell^2} \cdot \sup_{(\state, \action)
  \in \statespace \times \actionspace} \vecnorm{ \Mmat^{-1/2}
  \phi(\state, \action)}{\ell^2}.
\end{align*}
Combining above two bounds yields
\begin{multline}
\abss{ \plainvec^\top \Mmat^{-1} \cdot \tfrac{1}{\numobs} \sum_{i =
    1}^\numobs \Big\{ \noise_i \phi(\State_i, \Action_i) \Big\}} \\
\leq c \vecnorm{\big(\Mmat^{-1 / 2} \plainvec}{\ell^2} \Big\{
\varbound \sqrt{\tfrac{\log (1 / \delta)}{\numobs}} + \subgaussian
\sup_{(\state, \action) \in \statespace \times \actionspace}
\vecnorm{\Mmat^{-1/2} \phi(\state, \action)}{\ell^2} \cdot \tfrac{\log
  (1 / \delta) \log \numobs}{\numobs}
\Big\},\label{eq:noise-part-in-krr-directional-bound-conditional}
\end{multline}
with probability $1 - \delta$.

Combining
equations~\eqref{eq:bias-part-in-krr-directional-bound-conditional}
and~\eqref{eq:noise-part-in-krr-directional-bound-conditional}, we
conclude that
\begin{multline*}
  \abss{ \inprod{\plainvec}{\projparamhat_\numobs - \projparam_*}}
  \leq c \vecnorm{\Mmat^{-1/2} \plainvec}{\ell^2} \Big\{
  \sqrt{\reguhomo} \rkhsnorm{\treateff} + \varbound
  \sqrt{\tfrac{\log (1/\delta)}{\numobs}} + \subgaussian
  \sup_{(\state, \action) \in \statespace \times \actionspace}
  \vecnorm{\Mmat^{-1/2} \phi(\state, \action)}{\ell^2} \cdot
  \tfrac{\log (1 / \delta) \log \numobs}{\numobs} \Big\},
\end{multline*}
which completes the proof
of~\Cref{lemma:conditional-high-prob-bound-for-krr-directional}.


\subsection{Proof of~\Cref{lemma:emp-cov-concentration}}
\label{subsubsec:proof-lemma-emp-cov-concen}

For use in this proof, we note that the population-level covariance
operator $\bigcov_{*, q}$ satisfies the sandwich relation
\begin{align}
\label{eq:sigma-q-sandwich}  
\underline{q} \IdMat \preceq \bigcov_{*, q} \preceq \overline{q}
\IdMat.
\end{align}
Our argument adopts the approach used in the
paper~\cite{ma2022optimally}, but involves more refined arguments so
as to obtain sharper bounds that allow small value of
$\regu_\numobs$. By multiplying with the operator $(\bigcov_{*, q} +
\regu_\numobs \EigenMat^{-1})^{-1/2}$ from both the left and the right
of equation~\eqref{eq:matrix-dominance-in-emp-cov-concentration}, we
find that it suffices to bound the operator norm of the following
pre-conditioned error operator:
\begin{align*}
\Delhat_\numobs \mydefn (\bigcov_{*, q} + \regu_\numobs
\EigenMat^{-1})^{-1/2} \big( \widehat{\bigcov}_{\numobs, q} -
\bigcov_{*, q} \big) (\bigcov_{*, q} + \regu_\numobs
\EigenMat^{-1})^{-1/2}.
\end{align*}
Note that $\Delhat_\numobs$ is sum of $\mathrm{i.i.d.}$ random
operators. In order to bound its operator norm, we invoke a known
Bernstein inequality in Hilbert spaces.  It applies to an $\myiid$
sequence $\{X_i\}_{i=1}^\numobs$ of self-adjoint zero-mean operators
on a separable Hilbert space $\vecspace$.
\begin{proposition}[Minsker~\cite{minsker2017some}]
\label{prop:minsker}
Consider a sequence such that
\begin{align*}
\opnorm{\Exs [X_i^2]} \leq \sigma^2, \quad \trace (\Exs [X_i^2]) \leq
V < \infty, \quad \mbox{and} \quad \opnorm{X_i} \leq U, ~
\mbox{almost surely}.
\end{align*}
Then we have the concentration inequality
\begin{align*}
  \Prob \Big( \opnorm{\sum_{j = 1}^n X_i} \geq t \Big) \leq \frac{14
    V}{\sigma^2} \exp \left( - \frac{t^2 / 2}{\numobs \sigma^2 + t U /
    3} \right), \quad \mbox{for any $t > 0$}.
\end{align*}
\end{proposition}
\noindent A form of this result is stated as as Theorem 3.1 in the
paper~\cite{minsker2017some}; see also \S 3.1 of the same paper for
the extension to the infinite-dimensional case.

\medskip

Using this auxiliary result, let us now prove~\Cref{lemma:emp-cov-concentration}. In doing so, we make use the
shorthand notation $\phi^i \mydefn \phi(\State_i, \Action_i)$ and $q^i
\mydefn q(\State_i, \Action_i)$, along with the sequence of random
linear 
operators
\begin{align*}
\Delta_i \mydefn (\bigcov_{*, q} + \regu_\numobs
\EigenMat^{-1})^{-1/2} \big( q^i \phi^i (\phi^i)^\top - \bigcov_{*, q}
\big) (\bigcov_{*, q} + \regu_\numobs \EigenMat^{-1})^{-1/2} \qquad
\mbox{for each $i \in [\numobs]$.}
\end{align*}
We need to bound the relevant quantities required to
apply~\Cref{prop:minsker}.  Beginning with the variance, we have
\begin{align}
\Exs \big[ \Delta_i^2 \big] &\preceq \Exs \Big[ \Big\{ (\bigcov_{*, q}
  + \regu_\numobs \EigenMat^{-1})^{-1/2} \big( q^i \phi^i
  (\phi^i)^\top - \bigcov_{*, q} \big) (\bigcov_{*, q} + \regu_\numobs
  \EigenMat^{-1})^{-1/2} \Big\}^{2}\Big] \nonumber \\ & = \Exs \Big[
  \Big\{q^i (\phi^i)^\top (\bigcov_{*, q} + \regu_\numobs
  \EigenMat^{-1})^{-1} \phi^i \Big\} \cdot \Big\{ (\bigcov_{*, q} +
  \regu_\numobs \EigenMat^{-1})^{-1/2} q^i \phi^i (\phi^i)^\top
  (\bigcov_{*, q} + \regu_\numobs \EigenMat^{-1})^{-1/2}\Big\}
  \Big]. \label{eq:delta-sq-bound-in-cov-conc-proof}
\end{align}
Define the quantity $\Phi_{\max} \mydefn \sup_{(\state, \action)}
\abss{q(\state, \action) \phi(\state, \action)^\top (\bigcov_{*, q} +
  \regu_\numobs \EigenMat^{-1})^{-1} \phi(\state, \action) }$.  We can
use this uniform bound to control the right-hand-side of the
relation~\eqref{eq:delta-sq-bound-in-cov-conc-proof}, and obtain
\begin{align*}
\Exs \big[ \Delta_i^2 \big] & \preceq \Phi_{\max}\cdot \Exs \Big[
  (\bigcov_{*, q} + \regu_\numobs \EigenMat^{-1})^{-1/2} q^i \phi^i
  (\phi^i)^\top (\bigcov_{*, q} + \regu_\numobs \EigenMat^{-1})^{-1/2} \Big]
\\
& = \Phi_{\max} \cdot (\bigcov_{*, q} + \regu_\numobs \EigenMat^{-1})^{-1/2}
\bigcov_{*, q} (\bigcov_{*, q} + \regu_\numobs \EigenMat^{-1})^{-1/2}
\end{align*}
We can then bound the operator norm and trace of $\Exs [\Delta_i^2]$
as
\begin{subequations}
\label{eq:ingredients-for-hilbert-space-operator-conc}
\begin{align}
\opnorm{\Exs \big[ \Delta_i^2 \big]} \leq \Phi_{\max} \cdot
\opnorm{(\bigcov_{*, q} + \regu_\numobs \EigenMat^{-1})^{-1/2} \bigcov_{*,
    q} (\bigcov_{*, q} + \regu_\numobs \EigenMat^{-1})^{-1/2}} \leq
\Phi_{\max},
\end{align}
and by equation~\eqref{eq:sigma-q-sandwich}, we have
\begin{align}
\trace \Big( \Exs \big[ \Delta_i^2 \big] \Big) & \leq \Phi_{\max}
\cdot \trace \big( (\bigcov_{*, q} + \regu_\numobs \EigenMat^{-1})^{-1}
\bigcov_{*, q} \big) \: \leq \overline{q} \Phi_{\max} \cdot \trace
\big( \underline{q} \IdMat + \regu_\numobs \EigenMat^{-1})^{-1}\big).
\end{align}
Finally, we note that
\begin{align}
\opnorm{\Delta_i} &\leq \trace \Big( (\bigcov_{*, q} + \regu_\numobs
\EigenMat^{-1})^{-1/2} q^i \phi^i (\phi^i)^\top (\bigcov_{*, q} +
\regu_\numobs \EigenMat^{-1})^{-1/2} \Big) \nonumber \\
& = q^i (\phi^i)^\top (\bigcov_{*, q} + \regu_\numobs \EigenMat^{-1})^{-1}
\phi^i \leq \Phi_{\max}
\end{align}
\end{subequations}
almost surely.

Combining the different parts of
equation~\eqref{eq:ingredients-for-hilbert-space-operator-conc}
with~\Cref{prop:minsker} yields the tail bound
\begin{align*}
   \Prob\left( \opnorm{\Delhat_\numobs} \leq t \right) \leq 14
   \overline{q} \trace \big( (\underline{q} \IdMat + \regu_\numobs
   \EigenMat^{-1})^{-1} \big) \cdot \exp \left\{ \frac{- \numobs t^2 / 2}{(1
     + t / 3) \Phi_{\max}} \right\}, \qquad \mbox{valid for any $t >
     0$.}
\end{align*}
Noting that $ \trace \big((\underline{q} \IdMat + \regu_\numobs
\EigenMat^{-1})^{-1} \big) \leq \regu_\numobs^{-1} \trace (\EigenMat)
\leq \tfrac{\kappa^2}{\regu_\numobs}$, for the event $\Event_\delta$
defined as
\begin{align*}
\Event_\delta \mydefn \left\{ \opnorm{\Delhat_\numobs} \leq
\sqrt{\frac{ 2 \Phi_{\max}}{\numobs} \log \big( \frac{\kappa^2
  }{\regu_\numobs \delta} \big)} + \frac{6 \Phi_{\max}}{\numobs} \log
\big( \frac{\kappa^2 }{\regu_\numobs \delta} \big) \right\},
\end{align*}
we have $\Prob (\Event_\delta) \geq 1 - \delta$.

Note that the quantity $\Phi_{\max}$ admits the bound
\begin{align*}
\Phi_{\max} \leq \overline{q} \sup_{\state \in \Xspace, \action
  \in \actionspace} \abss{ \phi(\state, \action)^\top (\underline{q}
  \IdMat + \regu_\numobs \EigenMat^{-1})^{-1} \phi(\state, \action) } \leq
\overline{q} / \underline{q} \cdot \effdim \big(\regu_\numobs /
\underline{q} \big).
\end{align*}

On the event $\Event_\delta$, the
conditions~\eqref{eq:regu-requirement-in-emp-cov-concentration-general}
imply that $\opnorm{\Delhat_\numobs} \leq \offpar$. Therefore, we
conclude that the following bound holds true with probability $1 -
\delta$:
\begin{align*}
    (1 - \offpar) \IdMat \preceq (\bigcov_{*, q} + \regu_\numobs
  \EigenMat^{-1})^{-1/2} \big( \widehat{\bigcov}_{\numobs, q} +
  \regu_\numobs \EigenMat^{-1} \big) (\bigcov_{*, q} + \regu_\numobs
  \EigenMat^{-1})^{-1/2} \preceq (1 + \offpar) \IdMat,
\end{align*}
which completes the proof of~\Cref{lemma:emp-cov-concentration}.


\subsection{Proof of~\Cref{lemma:krr-directional-bound-weak-assumption}}
\label{app:subsec-proof-krr-directional-bound-weak-assumption}

Recall the error decomposition in the proof
of~\Cref{lemma:krr-directional-bound}:
\begin{align*}
    \projparamhat_\numobs - \projparam_* = \big(
    \widehat{\bigcov}_\numobs + \regu_\numobs \EigenMat^{-1}
    \big)^{-1} \cdot \frac{1}{\numobs} \sum_{i = 1}^\numobs \Big\{
    \noise_i \phi(\State_i, \Action_i) - \regu_\numobs \EigenMat^{-1}
    \projparam_* \Big\},
\end{align*}
where we define the noise function $\noise_i \mydefn \outcome_i -
\treateff(\State_i, \Action_i)$.

Since~\Cref{lemma:conditional-high-prob-bound-for-krr-directional}
does not depend on the
condition~\eqref{eq:sample-size-req-based-on-eigendecay} on the
effective dimension, conditionally on the state-action pairs
$(\State_i, \Action_i)_{i = 1}^\numobs$, with probability $1 -
\delta$, we have
\begin{multline}
  \abss{\plainvec^\top (\projparamhat_\numobs - \projparam_*)} \leq c
  \vecnorm{\big( \widehat{\bigcov}_{\numobs} + \regu_\numobs
    \EigenMat^{-1} \big)^{-1 / 2} \plainvec}{\ell^2}\\ \times \Big\{
  \sqrt{\regu_\numobs} \rkhsnorm{\treateff} + \varbound
  \sqrt{\tfrac{\log (1 / \delta)}{\numobs}} + \subgaussian
  \sup_{\state \in \statespace, \action \in \actionspace}
  \vecnorm{\big( \widehat{\bigcov}_\numobs + \regu_\numobs
    \EigenMat^{-1} \big)^{-1/2} \phi(\state, \action)}{\ell^2} \cdot
  \tfrac{\log (1 / \delta) \log \numobs}{\numobs}
  \Big\}.\label{eq:conditional-high-prob-bound-for-krr-copied}
\end{multline}
On the other hand, note that under Assumption~\eqref{EqnKerBou}, given
the regularization parameter
choice~\eqref{eq:regu-choice-in-weak-cor}, we have
\begin{align*}
    \log \Big( \frac{\kappa^2}{\regu_\numobs \delta} \Big)
    \frac{\effdim (\regu_\numobs)}{\numobs} \leq \log \Big(
    \frac{\kappa^2}{\regu_\numobs \delta} \Big)
    \frac{\kappa^2}{\numobs \regu_\numobs} \leq \frac{1}{32},
\end{align*}
which verifies the
condition~\eqref{eq:regu-requirement-in-emp-cov-concentration-general}
with $\overline{q} = \underline{q} = 1$ and $\offpar = 1/2$. Invoking
the empirical covariance concentration
lemma~\ref{lemma:emp-cov-concentration} with $q \equiv 1$ yields
\begin{align*}
  \frac{1}{2} (\IdMat + \regu_\numobs \EigenMat^{-1}) \preceq
  \widehat{\bigcov}_\numobs + \regu_\numobs \EigenMat^{-1} \preceq 2
  (\IdMat + \regu_\numobs \EigenMat^{-1}), \quad \mbox{with
    probability $1 - \delta$}.
\end{align*}
We can therefore control the the relevant terms in
equation~\eqref{eq:conditional-high-prob-bound-for-krr-copied},
leading to the following inequalities with probability $1 - \delta$.
\begin{align*}
\vecnorm{\big( \widehat{\bigcov}_{\numobs} + \regu_\numobs
  \EigenMat^{-1} \big)^{-1 / 2} \plainvec}{\ell^2} &\leq 2
\vecnorm{\big( \IdMat + \regu_\numobs \EigenMat^{-1} \big)^{-1 / 2}
  \plainvec}{\ell^2}, \quad \mbox{and}, \\
\sup_{(\state, \action) \in \statespace \times \actionspace}
\vecnorm{\big( \widehat{\bigcov}_\numobs + \regu_\numobs
  \EigenMat^{-1} \big)^{-1/2} \phi(\state, \action)}{\ell^2} &\leq
\regu_\numobs^{-1/2} \cdot \sup_{(\state, \action) \in \statespace
  \times \actionspace} \vecnorm{\EigenMat^{1/2} \phi(\state,
  \action)}{\ell^2} \leq \kappa / \sqrt{\regu_\numobs},
\end{align*}
where the last step follows from the uniform upper
bound~\eqref{EqnKerBou}.

Substituting these results back into
equation~\eqref{eq:conditional-high-prob-bound-for-krr-copied}, and
taking the regularization parameter according to
equation~\eqref{eq:regu-choice-in-weak-cor}, we conclude that
\begin{align*}
\abss{\plainvec^\top (\projparamhat_\numobs - \projparam_*)} & \leq 4
c \vecnorm{\big( \IdMat + \regu_\numobs \EigenMat^{-1} \big)^{-1 / 2}
  \plainvec}{\ell^2} \Big\{ \sqrt{\regu_\numobs} \rkhsnorm{\treateff}
+ \varbound \sqrt{\tfrac{\log (1 / \delta)}{\numobs}} +
\frac{\subgaussian \kappa \log (1 / \delta) \log \numobs}{\numobs
  \sqrt{\regu_\numobs}} \Big\}\\
& \leq c' \vecnorm{\big( \IdMat + \regu_\numobs \EigenMat^{-1}
  \big)^{-1 / 2} \plainvec}{\ell^2} \cdot (\subgaussian + \kappa
\Rbound )\sqrt{\frac{\log (\numobs / \delta) \log (1 /
    \delta)}{\numobs}},
\end{align*}
with probability $1 - \delta$, which
proves~\Cref{lemma:krr-directional-bound-weak-assumption}.


\subsection{Proof of~\Cref{lemma:weighted-krr-conditional-bound}}
\label{app:subsec-weighted-krr-conditional}

As with the proof
of~\Cref{lemma:conditional-high-prob-bound-for-krr-directional}, we
decompose the error into a noise and bias term---namely
\begin{align*}
  \inprod{\plainvec}{\projparamhat_\numobs - \projparam_*} & =
  \myunder{\frac{1}{\numobs} \sum_{i=1}^\numobs W_i}{noise part} -
  \myunder{\regustagethree \plainvec^\top
    \big(\widehat{\bigcov}_{\numobs}^\sigma + \regustagethree
    \EigenMat^{-1} \big)^{-1} \EigenMat^{-1} \projparam_*}{bias part},
\end{align*}
where $W_i \defn \plainvec^\top
\big(\widehat{\bigcov}_{\numobs}^\sigma + \regustagethree
\EigenMat^{-1} \big)^{-1} \noise_i \sighat_\numobs^{-2} \phi(\State_i,
\Action_i)$

For the bias part, applying the Cauchy--Schwarz inequality yields
\begin{align}
\regustagethree \abss{\plainvec^\top
  \big(\widehat{\bigcov}_{\numobs}^\sigma + \regustagethree
  \EigenMat^{-1} \big)^{-1} \EigenMat^{-1} \projparam_*} & \leq
\regustagethree \vecnorm{\EigenMat^{-1/2} \projparam_*}{\ell^2} \cdot
\vecnorm{\EigenMat^{-1/2} \big(\widehat{\bigcov}_{\numobs}^\sigma +
  \regustagethree \EigenMat^{-1} \big)^{-1} \plainvec}{\ell^2}
\nonumber \\
    & \overset{(i)}{\leq} \sqrt{\regustagethree} \rkhsnorm{\treateff}
\cdot \vecnorm{\big(\widehat{\bigcov}_{\numobs}^\sigma +
  \regustagethree \EigenMat^{-1} \big)^{-1/2} \plainvec}{\ell^2}
\nonumber \\
\label{eq:weighted-krr-conditional-bias-part-bound}    
& \overset{(ii)}{\leq} \frac{1}{\sqrt{\numobs}}
\vecnorm{\big(\widehat{\bigcov}_{\numobs}^\sigma + \regustagethree
  \EigenMat^{-1} \big)^{-1/2} \plainvec}{\ell^2},
\end{align}
where in step (i), we use the fact $\opnorm{\EigenMat^{-1/2}
  \big(\widehat{\bigcov}_{\numobs}^\sigma + \regustagethree
  \EigenMat^{-1} \big)^{-1/2}} \leq (\regustagethree)^{-1/2}$, and in
step (ii), we substitute with the regularization parameter choice
$\regustagethree = \frac{1}{\Rbound \numobs}$.

For the noise part, we use Adamczak's concentration inequality to
establish high-probability bounds. We start with the expression for
the conditional variance
\begin{align*}
    \Exs \Big[ W_i^2 \mid \State_i, \Action_i, \sighat_\numobs \Big] =
    \frac{\sigma^2(\State_i, \Action_i)}{\sighat_\numobs^4 (\State_i,
      \Action_i)} \Big(\plainvec^\top
    \big(\widehat{\bigcov}_{\numobs}^\sigma + \regustagethree
    \EigenMat^{-1} \big)^{-1} \phi(\State_i, \Action_i) \Big)^2,
\end{align*}
which leads to the bound
\begin{align}
\frac{1}{\numobs}\sum_{i = 1}^\numobs \Exs \big[ W_i^2 \mid \State_i,
  \Action_i, \sighat_\numobs \big] & \leq \max_{i \in [\numobs]}
\frac{\sigma^4 (\State_i, \Action_i)}{\sighat_\numobs^4 (\State_i,
  \Action_i)} \cdot z^\top \big(\widehat{\bigcov}_{\numobs}^\sigma +
\regustagethree \EigenMat^{-1} \big)^{-1}
\widehat{\bigcov}_\numobs^\sigma
\big(\widehat{\bigcov}_{\numobs}^\sigma + \regustagethree
\EigenMat^{-1} \big)^{-1} \plainvec \nonumber \\
\label{eq:second-moment-in-weighted-krr-conditional-ineq}
& \leq 4 z^\top \big(\widehat{\bigcov}_{\numobs}^\sigma +
\regustagethree \EigenMat^{-1} \big)^{-1} \plainvec.
\end{align}
In the last step, we use the fact $\sighat_\numobs^2(\State_i,
\Action_i) \geq \frac{1}{2} \sigma^2(\State_i, \Action_i)$ for any $i
\in [\numobs]$.

On the other hand, for any $p > 0$, we have the conditional moment
bound
\begin{align}
 \left \{ \Exs \Big[ |W_i|^p \mid \State_i, \Action_i, \sighat_\numobs
   \Big] \right\}^{1/p} & \leq \sqrt{p} \subgaussian \cdot
 \abss{\plainvec^\top \big(\widehat{\bigcov}_{\numobs}^\sigma +
   \regustagethree \EigenMat^{-1} \big)^{-1} \sighat_\numobs^{-2}
   \phi(\State_i, \Action_i)} \nonumber \\
\label{eq:orlicz-bound-in-weighted-krr-conditional-ineq}  
& \leq \frac{2 \sqrt{p} \subgaussian }{\lowersigma^2}
\vecnorm{\big(\widehat{\bigcov}_{\numobs}^\sigma + \regustagethree
  \EigenMat^{-1} \big)^{-1/2} \plainvec}{\ell^2} \cdot \sup_{( \state,
  \action)} \vecnorm{\big(\widehat{\bigcov}_{\numobs}^\sigma +
  \regustagethree \EigenMat^{-1} \big)^{-1 / 2} \phi(\state,
  \action)}{\ell^2}.
\end{align}
Combining
equations~\eqref{eq:second-moment-in-weighted-krr-conditional-ineq}
and~\eqref{eq:orlicz-bound-in-weighted-krr-conditional-ineq} with
Adamczak's inequality, we conclude that
\begin{multline}
\label{eq:weighted-krr-conditional-noise-part-bound}
|\frac{1}{\numobs} \sum_{i=1}^\numobs W_i| \leq c
\vecnorm{\big(\widehat{\bigcov}_{\numobs}^\sigma + \regustagethree
  \EigenMat^{-1} \big)^{-1/2} \plainvec}{\ell^2} \Biggr \{ 2
\sqrt{\frac{\log (1 / \delta)}{\numobs}} \\
+ \frac{\log \numobs \log (1 / \delta) \subgaussian }{\numobs
  \lowersigma^2 } \sup_{(\state, \action)}
\vecnorm{\big(\widehat{\bigcov}_{\numobs}^\sigma + \regustagethree
  \EigenMat^{-1} \big)^{-1 / 2} \phi(\state, \action)}{\ell^2} \Biggr
\}.
\end{multline}
Finally, putting together
equations~\eqref{eq:weighted-krr-conditional-bias-part-bound}
and~\eqref{eq:weighted-krr-conditional-noise-part-bound} completes the
proof of this lemma.


\section{Conditional variance estimation and robust risk property}
\label{subsubsec:cond-var-est-example}

In this section, we discuss the problem of estimating the conditional
variance function \mbox{$(\state, \action) \mapsto \sigma^2(\state,
  \action)$.}


\subsection{Some conditional variance estimators}

In this section, we construct concrete estimators for the conditional
variance $\sigma^2$ that satisfy the robust pointwise risk
property. In combination with the four-stage
framework~\eqref{eq:four-stage-framework}, these results immediately
lead to instance-optimal results in~\Cref{ThmHetero}.

\paragraph{Kernel ridge regression:}
Consider a positive semi-definite kernel function $\kernelfuncprime:
(\Xspace \times \actionspace) \times (\Xspace \times \actionspace)
\rightarrow \real$ that defines an RKHS $\rkhsprime$ with the Mercer
decomposition
\begin{align}
\label{eq:mercer-decomposition-prime}  
\kernelfuncprime \big( (\state_1, \action_1), (\state_2, \action_2)
\big) = \sum_{j = 1}^{\infty} \eigen_j \basisfunc_j (\state_1,
\action_1) \basisfunc_j (\state_2, \action_2).
\end{align}
We assume that the RKHS $\rkhsprime$ satisfies the regularity
assumption~\eqref{EqnKerBou}, and that the true conditional variance
function lies in this RKHS, i.e.,
\begin{align}
\rkhsprimenorm{\sigma^2} \leq \Rbound^\sigma.
\end{align}
Following the definition~\eqref{eq:defn-effdim}, for any $\regu > 0$
we define $\effdim_\sigma (\regu)$ as the effective dimension
associated to the regularization parameter $\regu > 0$ for the RKHS
$\rkhsprime$.

We consider the penalized least-square estimator
\begin{align}
\label{eq:krr-for-cond-var}  
\sighat_\numobs^2 \mydefn \arg\min_{h} \Big\{ \frac{1}{\numobs}
\sum_{i = 1}^\numobs \big(Z_i - h (\State_i, \Action_i) \big)^2 +
\regu \rkhsprimenorm{h}^2 \Big\}.
\end{align}
\begin{proposition}
\label{prop:krr-robust-ptwise-risk}
Let $\regu_0 (\varepsilon)$ be the smallest value of $\regu$ such that
$\regu \effdim_\sigma (\regu) \leq \big(\varepsilon /
\Rbound^\sigma\big)^{2} $, the estimator~\eqref{eq:krr-for-cond-var}
with parameter choice $\regu_\numobs =\regu_0 (\varepsilon) $
satisfies the robust pointwise risk property with
\begin{align}
\label{eq:robust-pt-wise-risk-for-krr}  
\maux(\varepsilon, \delta) \mydefn c \frac{\subgaussian^4 \log (1 /
  \delta)}{\varepsilon^2} \effdim_\sigma (\regu_0 (\varepsilon)) + c
\frac{(\Rbound^\sigma)^2}{\subgaussian^4 \regu_0 (\varepsilon)}, \quad
\mbox{and} \quad \overline{\smallbias} (\varepsilon, \delta) \mydefn
\frac{\varepsilon}{c \sqrt{\effdim_\sigma (\regu_0 (\varepsilon))}}.
   \end{align}
\end{proposition}
\noindent See~\Cref{app:subsec-proof-krr-robust-ptwise-risk} for the
proof.

\medskip

A few remarks are in order. First, \Cref{prop:krr-robust-ptwise-risk}
requires that the effective dimension of the RKHS $\rkhs_\sigma$ to
satisfy that $\regu \effdim_\sigma (\regu) \rightarrow 0$ for $\regu
\rightarrow 0^+$. A similar condition is also imposed on the RKHS
$\rkhs$ used to estimate the treatment effect function, which can be
verified under certain conditions on the eigenfunctions. (see
equation~\eqref{eq:effdim-upper-bound-under-eig} and
\Cref{prop:effective-dim} in the appendix for the statement of such
results.) In particular, suppose that the effective dimension
satisfies a decay condition $\effdim (\regu) \leq \effdim_0
\regu^{\offpar - 1}$ for some scalar $\offpar \in (0, 1]$, by seeing
the scalars $(\effdim_0, \Rbound^\sigma, \subgaussian)$ as constants,
we choose $\regu_0 (\varepsilon) = \varepsilon^{\frac{2}{\offpar}}$
the condition~\eqref{eq:robust-pt-wise-risk-for-krr} becomes
\begin{align*}
 \maux(\varepsilon, \delta) \asymp \varepsilon^{- 2 / \offpar} \log (1
 / \delta) \quad \mbox{and} \quad \overline{\smallbias} (\varepsilon,
 \delta) \asymp \varepsilon^{1 / \offpar},
\end{align*}
Such a requirement on the sample size $\maux$ and the bias upper bound
$\overline{\smallbias}$ may not always achieve the optimal rate for
estimating the function $\sigma^2$. However, since we only need the
estimation error to be smaller than a constant $\lowersigma^2/2$, as
required in equation~\eqref{eq:hetero-sample-size-req-second-stage}, a
polynomial dependency on the accuracy level $\varepsilon$ and
poly-logarithmic dependency on the failure probability $ \delta$
suffices our purposes.

\paragraph{Local average estimator:}
Let the statespace $\Xspace$ be a compact subset of $\real^\usedim$
and let the action space $\actionspace$ be discrete.  Define the class
of $L$-Lipschitz functions as
\begin{align*}
\funcClass_L = \Big\{ f: \Xspace \rightarrow \real, \abss{f(\state) -
  f (y)} \leq L \vecnorm{x - y}{2}~ \mbox{for any $x, y
  \in \Xspace$}\Big\}
\end{align*}
We assume that the conditional variances are smooth enough.
\begin{align}
  \sigma^2(\cdot, \action) \in \funcClass_L, \quad \mbox{for each
    $\action \in \actionspace$.}\label{eq:sigmafunc-belong-to-holder}
\end{align}
To make estimation possible with random design, we need an additional
regularity assumption on the density.
\begin{align}
\label{eq:density-assumption-for-local-avg}  
\inf_{\state \in \Xspace, \action \in \actionspace} \big(
\Joint{\probxstar}{\propscore} \big) \Big( \ball (\state, r) \times
\{\action\} \Big) \geq p_0 r^{d_0}, \quad \mbox{for any $r \in (0,
  \localradius)$.}
\end{align}
Given a tuning parameter $r_\numobs > 0$, we consider the local
averaging estimator
\begin{align}
\label{eq:local-avg-estimator}  
\sighat_\numobs^2(\xzero, \actionzero) \mydefn |\avgset|^{-1} \sum_{i
  \in \avgset} Z_i, \quad \mbox{where } \avgset \mydefn \Big\{ i \in
       [\numobs] : \State_i \in \ball (\xzero, r_\numobs), \Action_i =
       \actionzero \Big\}.
\end{align}
\begin{proposition}
\label{prop:local-avg-robust-ptwise-risk}
For any $\varepsilon > 0$ and $\delta \in (0, 1)$, there exists a
universal constant $c > 0$, such that the
estimator~\eqref{eq:local-avg-estimator} satisfies the robust
pointwise risk property with
\begin{align*}
 \maux(\varepsilon, \delta) = \frac{L^{2 \usedim_0} \log (1 /
   \delta)}{p_0} \Big( \frac{c }{\varepsilon} \Big)^{\usedim_0 + 2} +
 \frac{ \log (1 / \delta)}{\Ltwospace p_0 \localradius^{\usedim_0 + 1}}
 L^{\usedim_0} \frac{\log (1 / \delta)}{p_0} \log^{\usedim_0 + 2} (1 /
 \varepsilon), \quad \mbox{and} \quad\overline{\smallbias}
 (\varepsilon, \delta) = \varepsilon / 2.
 \end{align*}
\end{proposition}
\noindent See~\Cref{app:subsec-proof-local-avg-robust-ptwise-risk} for
the proof.

\medskip

A few remarks are in order. Compared to
\Cref{prop:krr-robust-ptwise-risk}, the local average estimator only
requires the target function $\sigma^2 (\cdot, a)$ to be Lipschitz,
for any $a \in \actionspace$. In dimension larger than $1$, this
usually requires less order of smoothness than the RKHS case
in~\Cref{prop:krr-robust-ptwise-risk}, while being less flexible with
the structure of the function class. The regularity
condition~\eqref{eq:density-assumption-for-local-avg} ensures that any
small ball in $\Xspace$ and any action $\action$ get sufficiently
large probability of being sampled. For example, when the probability
distribution $\probx$ has a density function uniformly bounded by
$\probx_{\min} > 0$, and when the probability of choosing any action
$\action$ is at least $\propscore_{\min}$, the
condition~\eqref{eq:density-assumption-for-local-avg} is satisfied
with $p_0 = c_d \probx_{\min} \propscore_{\min}$ and $\usedim_0 =
\usedim$, for a constant $c_d > 0$ depending only on $d$. More
generally, even if the function $(\state, \action) \mapsto
\probxstar(\state) \propscore(\state, \action)$ can attain $0$ at some
points, as long as appropriate growth conditions are imposed around
these points, the
condition~\eqref{eq:density-assumption-for-local-avg} will still be
satisfied. Finally, though we only study the Lipschitz case, in
literature optimal results for general H\"{o}lder classes have been
established for the estimation problems of conditional
variance~\cite{shen2020optimal}. In combination with their results, we
can also obtain optimal instance-dependent guarantees
in~\Cref{ThmHetero}.


\subsection{Proofs of robust pointwise risk properties}
\label{app:subsec-proof-robust-ptwise}

In this appendix, we establish the robust pointwise risk properties
for various estimators discussed
in~\Cref{subsubsec:cond-var-est-example}.


\subsubsection{Proof of~\Cref{prop:krr-robust-ptwise-risk}}
\label{app:subsec-proof-krr-robust-ptwise-risk}

The proof is similar to that of \Cref{lemma:krr-directional-bound},
with specific treatment given to the deterministic bias part. Define
the infinite-dimensional vectors
\begin{align*}
\projparam_* \mydefn \Psi (\sigma^2) \quad \mbox{and} \quad
\projparamhat_\numobs \mydefn \Psi (\sighat_\numobs^2).
\end{align*}
We can represent the error using basis functions.
\begin{align*}
\sigma^2(\state_0, \actionzero) - \sighat^2_\numobs (\state_0,
\actionzero) = \inprod{\projparamhat_\numobs -
  \projparam_*}{\phi(\state_0, \actionzero)}.
\end{align*}
Defining the noise and bias parts
\begin{align*}
    \noise_i \mydefn Z_i - \Exs [Z_i \mid \State_i, \Action_i], \quad
    \mbox{and} \quad \smallbias(\State_i, \Action_i) \mydefn \Exs[Z_i
      \mid \State_i, \Action_i] - \sigma^2(\State_i, \Action_i).
\end{align*}
We also define the empirical covariance operator
\begin{align*}
    \widehat{\bigcov}_\numobs \mydefn \frac{1}{\numobs} \sum_{i =
      1}^\numobs \phi(\State_i, \Action_i) \phi(\State_i,
    \Action_i)^\top,
\end{align*}
the error vector $\projparamhat_\numobs - \projparam_*$ admits a
representation
\begin{align}
\projparamhat_\numobs - \projparam_* = \big( \widehat{\bigcov}_\numobs
+ \regu \EigenMat^{-1} \big)^{-1} \frac{1}{\numobs} \sum_{i =
  1}^\numobs \Big\{ \noise_i (\State_i, \Action_i) \phi(\State_i,
\Action_i) + \smallbias (\State_i, \Action_i) \phi(\State_i,
\Action_i) - \regu \EigenMat^{-1} \projparam_* \Big\}.
\end{align}
Define the event
\begin{align*}
\Event_{\varepsilon, \delta} \mydefn \Big\{ \max_{1 \leq i \leq
  \numobs} \abss{\smallbias (\State_i, \Action_i)} \leq
\overline{\smallbias} (\varepsilon, \delta) \Big\}.
\end{align*}

Clearly, the error consists of three parts: a part induced by
stochastic (unbiased) noise $\noise_i$; a part involving the
observation bias $\smallbias (\State_i, \Action_i)$; and the bias
introduced by the regularization $\regu$. We claim that the following
bounds hold true with probability $1 - \delta$ on the event
$\Event_{\varepsilon, \delta}$.
\begin{subequations}
\begin{align}
    \abss{\phi(\state_0, \actionzero)^\top (\widehat{\bigcov}_\numobs
      + \regu \EigenMat^{-1} \big)^{-1} \frac{1}{\numobs} \sum_{i =
        1}^\numobs \noise_i \phi(\State_i, \Action_i) } &\leq c
    \big(\subgaussian^2 + \overline{\smallbias} (\varepsilon, \delta)
    \big) \sqrt{\frac{\effdim_\sigma(\regu) \log (1 /
        \delta)}{\numobs}
    }. \label{eq:noise-part-in-krr-one-point-bound}\\ \abss{\phi(\state_0,
      \actionzero)^\top (\widehat{\bigcov}_\numobs + \regu
      \EigenMat^{-1} \big)^{-1} \frac{1}{\numobs} \sum_{i = 1}^\numobs
      \smallbias (\State_i, \Action_i) \phi(\State_i, \Action_i) }
    &\leq c \overline{\smallbias} (\varepsilon, \delta)
    \sqrt{\effdim_\sigma(\regu)}
    ,\label{eq:bias-part-in-krr-one-point-bound}\\ \abss{\phi(\state_0,
      \actionzero)^\top (\widehat{\bigcov}_\numobs + \regu
      \EigenMat^{-1} \big)^{-1} \regu \EigenMat^{-1} \projparam_* }
    &\leq c \Rbound^\sigma \sqrt{\regu
      \effdim_\sigma(\regu)}. \label{eq:regu-part-in-krr-one-point-bound}
\end{align}
\end{subequations}
Taking these three bounds as given, for any $\varepsilon > 0$ and
$\delta \in (0, 1)$, we take $\regu_0 (\varepsilon)$ be the smallest
value of $\regu$ such that $\regu \effdim_\sigma (\regu) \leq \big(
\frac{\varepsilon}{\Rbound^\sigma} \big)^2$ (which is guaranteed to
exist for if $\numobs \cdot \effdim_\sigma (1 / \numobs) \rightarrow
0$), the robust pointwise risk condition is satisfied with
\begin{align*}
    \maux (\varepsilon, \delta) = c \frac{ \subgaussian^4 \log (1 /
      \delta)}{\varepsilon^2} \effdim_\sigma (\regu_0 (\varepsilon)) +
    c \frac{(\Rbound^\sigma)^2}{\subgaussian^4 \regu_0 (\varepsilon)},
    \quad \mbox{and} \quad \overline{\smallbias} =
    \frac{\varepsilon}{c \sqrt{\effdim_\sigma (\regu_0 (\varepsilon))}
    },
\end{align*}
completing the proof of~\Cref{prop:krr-robust-ptwise-risk}.

The rest of this section is devoted to the proofs of
equations~\eqref{eq:noise-part-in-krr-one-point-bound}--~\eqref{eq:regu-part-in-krr-one-point-bound}.

\paragraph{Proof of equation~\eqref{eq:noise-part-in-krr-one-point-bound}:}
By definition, note that the noise satisfies the conditional
$\psi_1$-norm bound
\begin{align*}
    \vecnorm{\noise_i \mid \State_i, \Action_i}{\psi_1} \leq 4
    (\subgaussian^2 + |\smallbias (\State_i, \Action_i)|).
\end{align*}
Invoking Adamczak's concentration inequality, conditionally on
$(\State_i, \Action_i)_{i = 1}^\numobs$, with probability $1 -
\delta$, we have
\begin{multline*}
  \abss{\phi(\state_0, \actionzero)^\top (\widehat{\bigcov}_\numobs +
    \regu \EigenMat^{-1} \big)^{-1} \frac{1}{\numobs} \sum_{i =
      1}^\numobs \noise_i (\State_i, \Action_i) \phi(\State_i,
    \Action_i)} \\
  \leq c \big( \subgaussian^2 + \max_{i \in [\numobs]} |\smallbias
  (\State_i, \Action_i) | \big) \times \Bigg\{
  \vecnorm{\widehat{\bigcov}_\numobs^{1/2} (\widehat{\bigcov}_\numobs
    + \regu \EigenMat^{-1} \big)^{-1} \phi(\state_0,
    \actionzero)}{\ell^2} \sqrt{\frac{\log (1 / \delta)}{\numobs} }
  \\ + \sup_{\state', \action'} \abss{\phi(\state', \action')^\top
    (\widehat{\bigcov}_\numobs + \regu \EigenMat^{-1} \big)^{-1}
    \phi(\state_0, \actionzero)} \frac{\log \numobs \log
    (1/\delta)}{\numobs} \Bigg\}.
\end{multline*}
On the event $\Event_{\varepsilon, \delta}$, we have $\max_{i \in
  [\numobs]} |\smallbias (\State_i, \Action_i) | \leq
\overline{\smallbias} (\varepsilon,
\delta)$. By~\Cref{lemma:emp-cov-concentration}, with probability $1
-\delta$, we have
\begin{align*}
  \sup_{\state, \action} \vecnorm{ (\widehat{\bigcov}_\numobs + \regu
    \EigenMat^{-1} \big)^{-1 / 2} \phi(\state, \action)}{\ell^2} \leq
  2 \sup_{\state, \action} \vecnorm{ (\IdMat + \regu \EigenMat^{-1}
    \big)^{-1 / 2} \phi(\state, \action)}{\ell^2} \leq c
  \sqrt{\effdim_\sigma (\regu)}.
\end{align*}
Putting together the pieces completes the proof of
equation~\eqref{eq:noise-part-in-krr-one-point-bound}.


\paragraph{Proof of equation~\eqref{eq:bias-part-in-krr-one-point-bound}:}
Applying the Cauchy--Schwarz inequality to the finite summation yields
\begin{align*}
     &\abss{\phi(\state_0, \actionzero)^\top
    (\widehat{\bigcov}_\numobs + \regu \EigenMat^{-1} \big)^{-1}
    \frac{1}{\numobs} \sum_{i = 1}^\numobs \smallbias (\State_i,
    \Action_i) \phi(\State_i, \Action_i) }^2 \\
& \leq \frac{1}{\numobs} \sum_{i = 1}^\numobs \abss{\phi(\state_0,
    \actionzero)^\top (\widehat{\bigcov}_\numobs + \regu
    \EigenMat^{-1} \big)^{-1} \smallbias (\State_i, \Action_i)
    \phi(\State_i, \Action_i)}^2 \\ &= \max_{i \in [\numobs]}
  \smallbias^2(\State_i, \Action_i) \cdot \frac{1}{\numobs} \sum_{i =
    1}^\numobs \phi(\state_0, \actionzero)^\top
  (\widehat{\bigcov}_\numobs + \regu \EigenMat^{-1} \big)^{-1}
  \phi(\State_i, \Action_i) \phi(\State_i, \Action_i)^\top
  (\widehat{\bigcov}_\numobs + \regu \EigenMat^{-1} \big)^{-1}
  \phi(\state_0, \actionzero) \\
& \leq \max_{i \in [\numobs]} \smallbias^2(\State_i, \Action_i) \cdot
  \phi(\state_0, \actionzero)^\top (\widehat{\bigcov}_\numobs + \regu
  \EigenMat^{-1} \big)^{-1} \phi(\state_0, \actionzero)
\end{align*}
On the other hand, we can apply~\Cref{lemma:emp-cov-concentration} to
the empirical covariance operator $\widehat{\bigcov}_\numobs$ in the
RKHS $\rkhsprime$. Given the regularization parameter $\regu \geq c
\kappa^2 \frac{\log (\numobs/ \delta)}{\numobs}$, with probability $1
- \delta$, we have
\begin{align*}
\frac{1}{2} (\IdMat + \regu \EigenMat^{-1}) \preceq
\widehat{\bigcov}_\numobs + \regu \EigenMat^{-1} \preceq 2 (\IdMat +
\regu \EigenMat^{-1}).
\end{align*}
Consequently, on the event $\Event_{\varepsilon, \delta}$, we have the
upper bound
\begin{multline*}
\abss{\phi(\state_0, \actionzero)^\top (\widehat{\bigcov}_\numobs +
  \regu \EigenMat^{-1} \big)^{-1} \frac{1}{\numobs} \sum_{i =
    1}^\numobs \smallbias (\State_i, \Action_i) \phi(\State_i,
  \Action_i) } \leq c \overline{\smallbias} (\varepsilon, \delta)
\cdot \sup_{\state, \action} \vecnorm{(\widehat{\bigcov}_\numobs +
  \regu \EigenMat^{-1} \big)^{-1/2} \phi (\state,
  \action)}{\ell^2}^2\\ \leq 2c \overline{\smallbias} (\varepsilon,
\delta) \cdot \sup_{\state, \action} \vecnorm{(\IdMat + \regu
  \EigenMat^{-1} \big)^{-1/2} \phi (\state, \action)}{\ell^2}^2 = 2 c
\overline{\smallbias} (\varepsilon, \delta) \effdim_\sigma (\regu).
\end{multline*}
with probability $1 - \delta$.


\paragraph{Proof of equation~\eqref{eq:regu-part-in-krr-one-point-bound}:}

By Cauchy--Schwarz inequality, we note that
\begin{align*}
& \abss{\phi(\state_0, \actionzero)^\top (\widehat{\bigcov}_\numobs +
    \regu \EigenMat^{-1} \big)^{-1} \regu \EigenMat^{-1} \projparam_*
  } \\
& \leq \sqrt{\regu} \cdot \vecnorm{(\widehat{\bigcov}_\numobs + \regu
    \EigenMat^{-1} \big)^{-1/2} \phi(\state_0, \actionzero)}{\ell^2}
  \cdot \opnorm{(\widehat{\bigcov}_\numobs + \regu \EigenMat^{-1}
    \big)^{-1/2} (\regu \EigenMat^{-1})^{1/2}} \cdot
  \vecnorm{\EigenMat^{-1/2} \projparam_* }{\ell^2} \\
& \leq 2 \sqrt{\regu} \Rbound^\sigma
  \vecnorm{(\widehat{\bigcov}_\numobs + \regu \EigenMat^{-1}
    \big)^{-1/2} \phi(\state_0, \actionzero)}{\ell^2}.
\end{align*}
By~\Cref{lemma:emp-cov-concentration}, with probability $1 - \delta$,
we have
\begin{align*}
\vecnorm{(\widehat{\bigcov}_\numobs + \regu \EigenMat^{-1}
  \big)^{-1/2} \phi (\state_0, \actionzero)}{\ell^2} \leq 2
\vecnorm{(\IdMat + \regu \EigenMat^{-1} \big)^{-1/2} \phi(\state_0,
  \actionzero)}{\ell^2} \leq 2 \sqrt{\effdim_\sigma (\regu)},
\end{align*}
which proves equation~\eqref{eq:regu-part-in-krr-one-point-bound}.


\subsubsection{Proof of~\Cref{prop:local-avg-robust-ptwise-risk}}
\label{app:subsec-proof-local-avg-robust-ptwise-risk}

We start with a decomposition of the error
\begin{align}
\label{eq:err-decomp-local-avg}  
  \sighat_\numobs^2(\xzero, \actionzero) - \sigma^2 (\xzero,
  \actionzero) = \abss{\avgset}^{-1} \sum_{i \in \avgset} \Big\{
  \noise_i + \smallbias(\State_i, \Action_i) + \big(\sigma^2
  (\State_i, \Action_i) - \sigma^2 (\xzero, \actionzero)\big) \Big\},
\end{align}
where the noise $\noise_i$ is defined as $\noise_i \mydefn Z_i - \Exs
[Z_i | \State_i, \Action_i]$ for each $i \in [\numobs]$.

Recall from the definition of the set $\avgset$ that for each $i \in
\avgset$, we have $\Action_i = \actionzero$ and $\State_i \in \ball
(\xzero, \radiusnumobs)$. Applying the Lipschitz
condition~\ref{eq:sigmafunc-belong-to-holder} then leads to the bound
\begin{align}
\label{eq:local-avg-lip-bias-part}  
\abss{\avgset}^{-1} \sum_{i \in \avgset} \abss{\sigma^2 (\State_i,
  \Action_i) - \sigma^2 (\xzero, \actionzero)} \leq L \radiusnumobs.
\end{align}
Defining the event
\begin{align*}
  \Event_{\varepsilon, \delta} \mydefn \Big\{ \max_{i \in [\numobs]}
  |\smallbias (\State_i, \Action_i)| \leq \overline{\smallbias}
  (\varepsilon, \delta) \Big\},
\end{align*}
on this event, we can control the additional bias in the observations
\begin{align}
\label{eq:local-avg-additional-bias-part}  
  \abss{\avgset}^{-1} \sum_{i \in \avgset} \abss{\smallbias (\State_i,
    \Action_i)} \leq \overline{\smallbias} (\varepsilon, \delta).
\end{align}
For the stochastic noise, we claim the following bound holds true
whenever the tuning parameter $\radiusnumobs$ satisfies $\radiusnumobs
\leq \localradius$ and $\frac{p_0\numobs
  \radiusnumobs^{\usedim_0}}{\log^2 \numobs} \geq \log (1 / \delta)$.
\begin{align}
\label{eq:concentration-in-local-avg-proof}  
  \Big| \sum_{i \in \avgset} \noise_i \Big| \leq c \sqrt{\frac{\log (1
      / \delta)}{\numobs p_0 \radiusnumobs^{\usedim_0} } }, \quad
  \mbox{with probability $1 - \delta$, on the event
    $\Event_{\varepsilon, \delta}$}
\end{align}
We prove this inequality at the end of this section.

Combining
equations~\eqref{eq:local-avg-lip-bias-part},~\eqref{eq:local-avg-additional-bias-part},
and~\eqref{eq:concentration-in-local-avg-proof}, we choose the local
radius as
\begin{align*}
  \radiusnumobs \mydefn \Big\{ \frac{\log (1 / \delta)}{\Ltwospace p_0
    \numobs} \Big\}^{\frac{1}{\usedim_0 + 2}}.
\end{align*}
Whenever the sample size $\numobs$ satisfies
\begin{align*}
  \numobs \geq \frac{ \log (1 / \delta)}{\Ltwospace p_0
    \localradius^{\usedim_0 + 1}}, \quad \mbox{and} \quad
  \frac{\numobs}{\log^{\usedim_0 + 2} \numobs} \geq L^{\usedim_0}
  \frac{\log (1 / \delta)}{p_0},
\end{align*}
on the event $\Event_{\varepsilon, \delta}$, we have the upper bound
with probability $1 - \delta$,
\begin{align*}
  \abss{\sighat_\numobs^2 (\xzero, \actionzero) - \sigma^2 (\xzero,
    \actionzero)} \leq \overline{\smallbias} (\varepsilon, \delta) + c
  \cdot L^{\frac{2 \usedim_0}{\usedim_0 + 2}} \Big\{ \frac{\log (1 /
    \delta)}{p_0 \numobs} \Big\}^{\frac{1}{\usedim_0 + 2}}.
\end{align*}


\paragraph{Proof of equation~\eqref{eq:concentration-in-local-avg-proof}:}

We start by exhibiting a lower bound on the cardinality of the set
$\avgset$.  When the averaging radius $\radiusnumobs$ satisfies
$\radiusnumobs \leq \localradius$, by the density
condition~\eqref{eq:density-assumption-for-local-avg}, we have
\begin{align*}
p_* \mydefn \Prob \left( i \in \avgset \right) \geq p_0
\radiusnumobs^{\usedim_0}, \quad \mbox{for each $i \in [\numobs]$}.
\end{align*}
The indicators $\bm{1}_{i \in \avgset}$ are independent for each $i
\in [\numobs]$. By Chernoff bound in the entropy form, we have
\begin{align*}
    \Prob \left( |\avgset| \leq \frac{\numobs p_*}{2} \right) \leq
    \exp \Big\{ - \numobs \kull{p_* / 2}{p_*} \Big\} \leq \exp \left(
    - c p_* \numobs \right),
\end{align*}
for a universal constant $c > 0$.

Defining the event
\begin{align*}
 \Event' \mydefn \Big\{ |\avgset| \geq \frac{\numobs p_*}{2} \Big\},
\end{align*}
the concentration inequality above implies that $\Prob (\Event') \geq
1 - \delta$ whenever $\numobs p_* \geq c' \log (1 / \delta)$.

Let us condition on the state-action pairs $(\State_i, \Action_i)_{i =
  1}^\numobs$ such that the event $\Event_{\varepsilon, \delta} \cap
\Event'$ holds true, applying Adamczak's concentration inequality to
its summation, with probability $1 - \delta$ under the conditional
law, we have
\begin{align*}
  \Big| \sum_{i \in \avgset} \noise_i \Big| \leq c \big(\subgaussian^2
  + \overline{\smallbias} (\varepsilon, \delta) \big) \cdot \Big\{
  \sqrt{\frac{\log (1 / \delta)}{|\avgset|}} + \frac{\log (1 / \delta)
    \log \numobs}{|\avgset|} \Big\},
\end{align*}
for a universal constant $c > 0$.

Taking into account the random design points $(\State_i, \Action_i)_{i
  = 1}^\numobs$, as long as the sample size and the radius satisfies
\begin{align*}
  \numobs p_* \geq \numobs p_0 \radiusnumobs^{\usedim_0} \geq \log (1
  / \delta) \cdot \log^2 \numobs,
\end{align*}
with probability $1 - \delta$, we have
\begin{align*}
  \Big| \sum_{i \in \avgset} \noise_i \Big| \leq c \sqrt{\frac{\log (1
      / \delta)}{\numobs p_0 \radiusnumobs^{\usedim_0} } },
\end{align*}
which proves equation~\eqref{eq:concentration-in-local-avg-proof}.


\section{Proofs for the examples}

We collect the proofs for the examples in this section.


\subsection{Proof of~\Cref{cor:missing-data-rates}}
\label{app:subsec-proof-missing-data-rates}

We first establish the effective dimension
condition~\eqref{EqnEffdimBound} by verifying the sup-norm growth
bound~\eqref{EqnEigfunGrowth}. Doing so ensures that the optimal risk
is determined (up to universal constant factors) by the risk
functionals $\vprobx^2(\treateff) + \vwenlongsq$ and
$\vsupersq{\delta_{\xzero}}$.  We then use Theorems~\ref{ThmLower}
and~\ref{ThmHetero} to prove the
bounds~\eqref{eq:missing-data-example-ate-rate}
and~\eqref{eq:missing-data-example-cate-rate}, respectively.


\subsubsection{Establishing the effective dimension condition}

We start by establishing tight bounds on the sup-norm growth
condition~\eqref{EqnEigfunGrowth}, which comes with a non-trivial (yet
well-controlled) exponent $\smallexponent$. This result is of
independent interest, illustrating the growth of eigenfunctions as a
natural phenomenon for RKHS applied to data whose densities have
singularities.
\begin{lemma}
\label{prop:supnorm-growth-missing-data-example}
Under the set-up above, there exists a pair of positive constants
$c_1, c_2$ that depends only on $\alpha$, such that for each $j \geq
1$, the eigenfunctions $\phi_j$ (normalized with
$\vecnorm{\phi_j}{\Ltwospace (\probx)} = 1$) associated to eigenvalue
$\eigen_j$ satisfy
\begin{align*}
c_1 \eigen_j^{- \frac{\alpha}{4(\alpha + 2)}} \leq
\vecnorm{\phi_j}{\infty} \leq c_2 \eigen_j^{- \frac{\alpha}{4(\alpha +
    2)}}.
\end{align*}
\end{lemma}
\noindent See~\Cref{sec:app-growth-condition} for the proof of this
lemma.

Consequently, the condition~\eqref{EqnEigfunGrowth} is satisfied with
exponent $\smallexponent = \frac{\alpha}{4(\alpha + 2)}$ and constant
$\phimax$ depending only on $\alpha$. With eigenvalue decay $\eigen_j
\asymp j^{-2}$ of the first-order Sobolev space
(see~\cite{widom1963asymptotic}),~\Cref{prop:effective-dim} yields
\begin{align*}
\effdim(\regu) \asymp \regu^{- \frac{\alpha + 1}{\alpha + 2}}, \quad
\mbox{and} \quad \effdim(\regu_\numobs) / \numobs \asymp \numobs^{-
  \frac{1}{\alpha + 2}}.
\end{align*}
which ensures the regularity
condition~\eqref{eq:sample-size-req-based-on-eigendecay} for sample
size $\numobs$ larger than a threshold depending only on $\alpha$.


\subsubsection{Proof of equation~\eqref{eq:missing-data-example-ate-rate}}

Note that for any function $f \in \ball_\rkhs (1)$ and $x \in [0, 1]$,
the Cauchy--Schwarz inequality yields
\begin{align*}
|f(x)| = \abss{\int_0^x f'(t) dt} \leq \sqrt{x \cdot \int_0^x
  (f'(t))^2 dt} \leq \sqrt{x}.
\end{align*}

So we have $\vecnorm{\treateff}{\infty} \leq 1$ and consequently
$\vprobx (\treateff) \leq 1$.

By the definition~\eqref{EqnDefnFunctional}, the variance function
$\vwenlong$ is defined (up to universal constant factors) as the
optimum value of the following variational problem:

\begin{subequations}
  \label{eq:variational-problem-var-functional-in-missing-data}
\begin{align}
\label{eq:variational-problem-objective}
\sup_{q} \left\{ \int_0^1 q(x) dx \right\}, &\quad \mbox{such that} \\
\label{eq:variational-problem-constraint}    
q(0) = 0, ~ \int_0^1 (q' (x))^2 dx &\leq \numobs, \quad \mbox{and}
\quad \int_0^1 (1 - x)^\alpha q^2(x) dx \leq 1.
\end{align}
\end{subequations}

It suffices to establish upper and lower bounds on the variance
functional $\vwenlong$ under different regimes.

Our proof relies on a technical lemma regarding the
constraint~\eqref{eq:variational-problem-constraint}, stated as
\begin{lemma}
\label{lemma:perturb-func-bound-under-cons-for-missing-data-example}
Under the constraint~\eqref{eq:variational-problem-constraint}, we
have
\begin{align*}
\sup_{\state \in [0, 1]} |q(x)| \leq c_\alpha \numobs^{\frac{1 +
    \alpha}{2 (2 + \alpha)}}, \quad \mbox{and} \quad \sup_{\state \in
  [0, 1 - \varepsilon]} |q(x)| \leq c_{\alpha, \varepsilon}
\numobs^{1/4},
\end{align*}
for a constant $c_\alpha$ depending on $\alpha > 0$, and a constant
$c_{\alpha, \varepsilon}$ depending on $\alpha > 0$ and $\varepsilon
\in (0, 1)$.
\end{lemma}
\noindent
We prove this lemma
in~\Cref{subsubsec:proof-lemma-perturb-func-bound-under-cons-for-missing-data-example}.

\medskip 

Taking it as given, we now proceed the proof of
equation~\eqref{eq:missing-data-example-ate-rate}.  It suffices to
establish upper and lower bounds on the variance functional
$\vwenlong$ under different regimes.

\paragraph{Upper bounds on the variance functional:}

Given a function $q$ satisfying the constraint~\eqref{eq:variational-problem-constraint}, for any $\varepsilon \in (0, 1)$,
we decompose the integral $\int_0^1 q(x) dx$ into parts $\int_0^{1 -
  \varepsilon}$ and $\int_{1 - \varepsilon}^1$, and bound them in
different ways.

By the Cauchy--Schwarz inequality, we note that
\begin{align}
\Big(\int_0^{1 - \varepsilon} q(x) dx\Big)^2 \leq \int_0^{1 -
  \varepsilon} (1 - x)^\alpha q^2(x) dx \cdot \int_0^{1 - \varepsilon}
\frac{dx}{(1 - x)^\alpha} \leq \begin{cases} \frac{1}{1 - \alpha} &
  \alpha < 1,\\ \log (1 / \varepsilon) & \alpha = 1,\\ \frac{1}{\alpha
    - 1} \varepsilon^{1 - \alpha} & \alpha > 1.
    \end{cases} \label{eq:wide-part-bound-in-variational-problem-upper}
\end{align}
For the second part, integration-by-parts yields
\begin{align*}
    \int_{1 - \varepsilon}^1 q(x) dx = q(1) - (1 - \varepsilon) q(1 -
    \varepsilon) - \int_{1 - \varepsilon}^1 x q' (x) dx = \varepsilon
    q(1 - \varepsilon) + \int_{1 - \varepsilon}^1 (1 - x) q' (x) dx.
\end{align*}
For the integral term, applying the Cauchy--Schwarz inequality yields
\begin{align}
\label{eq:integral-bound-in-variational-problem-upper}  
|\int_{1 - \varepsilon}^1 (1 - x) q' (x) dx| \leq \rkhsnorm{q}
\cdot \sqrt{\int_{1 - \varepsilon}^1 (1 - x)^2 dx } \leq
\sqrt{\varepsilon^3 \numobs}.
\end{align}
By~\Cref{lemma:perturb-func-bound-under-cons-for-missing-data-example},
we have
\begin{align}
\label{eq:qstar-bound-in-variational-problem-upper}  
|\varepsilon q(1 - \varepsilon)  | \leq \varepsilon \numobs^{\frac{1 +
    \alpha}{2 (2 + \alpha)}}.
\end{align}
Combining
equations~\eqref{eq:wide-part-bound-in-variational-problem-upper},~\eqref{eq:integral-bound-in-variational-problem-upper},~\eqref{eq:qstar-bound-in-variational-problem-upper}
yields
\begin{align*}
    \abss{\int_0^1 q(x) dx} \leq \sqrt{\varepsilon^3 \numobs} +
    c_\alpha \varepsilon \numobs^{\frac{1 + \alpha}{2 (2 + \alpha)}} +
    c'_\alpha \times \begin{cases} 1 & \alpha < 1,\\ \sqrt{\log (1 /
        \varepsilon)} & \alpha = 1,\\ \varepsilon^{- (\alpha - 1) / 2}
      & \alpha > 1.
    \end{cases} 
\end{align*}
We consider three cases:
\begin{itemize}
\item When $\alpha < 1$, we take $\varepsilon = 0$, and obtain
      that $\vwenlong \lesssim_\alpha 1$.
\item When $\alpha = 1$, we take $\varepsilon = \numobs^{- 1}$, and
  obtain that $\vwenlong \lesssim \sqrt{\log \numobs}$.
\item When $\alpha > 1$, we take $\varepsilon =
  \numobs^{\frac{-1}{\alpha + 2}}$, and obtain that $\vwenlong
  \lesssim_\alpha \numobs^{\frac{\alpha - 1}{2(\alpha + 2)}}$.
\end{itemize}


\paragraph{Lower bounds on the variance functional:}
On the lower bound side, for positive scalars $\varepsilon \in (0, 1)$
and $h > 0$, we construct the function
\begin{align*}
q_{\varepsilon, h} (x) \mydefn \begin{cases} 0 & x \leq 1 -
  \varepsilon,\\ \frac{h}{\varepsilon} (x - 1 + \varepsilon) & x > 1 -
  \varepsilon.
\end{cases}
\end{align*}
Clearly, we have $q_{\varepsilon, h} (0) = 0$, and straightforward
calculation yields
\begin{align*}
    \int_0^1 q_{\varepsilon, h} (x) dx = \frac{\varepsilon h}{2},
    \quad \mbox{and}\quad \int_0^1 (q_{\varepsilon, h}' (x))^2 dx =
    \frac{h^2}{\varepsilon}, \quad \mbox{and} \quad \int_0^1 (1 -
    x)^\alpha q_{\varepsilon, h}^2(x) dx \leq \varepsilon^{\alpha + 1}
    h^2.
\end{align*}
For $\alpha > 1$, under the choice of parameters
\begin{align*}
    \varepsilon = \numobs^{\frac{- 1}{2 + \alpha}}, \quad
    \mbox{and}\quad h = \numobs^{\frac{1 + \alpha}{2 (2 + \alpha)}},
\end{align*}
we have that $\int_0^1 q_{\varepsilon, h} (x) dx \asymp
\numobs^{\frac{ \alpha - 1}{2 (2 + \alpha)}}$.

For $\alpha < 1$, taking $\varepsilon = 1$ and $h = 1$, we have that
$\int_0^1 q_{\varepsilon, h} (x) dx \asymp 1$.

Consequently, for $\alpha \neq 1$, we have the lower bounds
\begin{align}
 \vwenlonghil \gtrsim_\alpha \begin{cases} \numobs^{\frac{\alpha -
       1}{2(\alpha + 2)}}, & \alpha > 1,\\ 1 & \alpha < 1.
    \end{cases}  \label{eq:variational-problem-lower-bound-alpha-not-1}
\end{align}

For the case of $\alpha = 1$, we use a different construction. Define
the function
\begin{align*}
  q_\numobs (x) \mydefn \begin{cases} \frac{(\log \numobs)^{-1/2}}{1 -
      x}, & x \in [0, 1 - \numobs^{-1/3}],\\ \numobs^{1/3}(\log
    \numobs)^{-1/2}, & x \in [1 - \numobs^{-1/3}, 1].
    \end{cases}
\end{align*}
Straightforward calculation yields
\begin{align*}
    \int_0^1 (q_\numobs' (x))^2 dx = \frac{\numobs}{12 \log \numobs} <
    \numobs, \quad \mbox{and} \quad \int_0^1 q_\numobs^2(x) (1 - x)dx
    = \frac{1}{3} + \frac{1}{2 \log \numobs} < 1,
\end{align*}
which verifies the
constraint~\eqref{eq:variational-problem-constraint}.

Therefore, we can lower bound the quantity $\vwenlonghil$ using the
value of the variational problem at $q_\numobs$, leading to the result
\begin{align}
\vwenlonghil \geq \int_0^1 q_\numobs (x) dx \geq \frac{1}{3}
\sqrt{\log
  \numobs}. \label{eq:variational-problem-lower-bound-alpha-1}
\end{align}
Combining
equations~\eqref{eq:variational-problem-lower-bound-alpha-not-1}
and~\eqref{eq:variational-problem-lower-bound-alpha-1} completes the
proof of the lower bound on $\vwenlonghil$.


\subsubsection{Proof of equation~\eqref{eq:missing-data-example-cate-rate}}

By~\Cref{ThmLower} and the claim~\eqref{EqnCATEclaim}, the minimax
risk is determined by the variance functional
$\vsuper{\delta_{\xzero}}$, defined as the optimum value of the
variational problem
\begin{align}
\label{eq:variational-problem-for-cate-rate-missing-data}  
\sup_{q} |q(\xzero)|, \quad \mbox{such that } q(0) = 0, ~ \int_0^1 (q'
(x))^2 dx &\leq \numobs, \quad \mbox{and} \quad \int_0^1 (1 -
x)^\alpha q^2(x) dx \leq 1.
\end{align}
If $\xzero = 0$, we have the trivial solution
$\vsuper{\delta_{\xzero}} = 0$. The rest of this section deals with
the case of $\xzero = 1$ and $\xzero \in (0, 1)$, respectively.

\paragraph{Case I: $\xzero = 1$.}

By~\Cref{lemma:perturb-func-bound-under-cons-for-missing-data-example},
we have $|q(1)| \leq c_\alpha \numobs^{\frac{1 + \alpha}{2 (1 +
    \alpha)}}$ for any function $q$ satisfying the constraints in the
variational
problem~\eqref{eq:variational-problem-for-cate-rate-missing-data}. On
the other hand, consider the function
\begin{align*}
q(x) \mydefn \begin{cases} 0 & x \leq 1 - \numobs^{\frac{-1}{2 +
      \alpha}},\\ \numobs^{\frac{3 + \alpha}{2 (2 + \alpha)}} \big(x -
  1 + \numobs^{\frac{-1}{2 + \alpha}} \big) & x > 1 -
  \numobs^{\frac{-1}{2 + \alpha}}.
    \end{cases}
\end{align*}
Straightforward calculation verifies that the function $q$ satisfies
the constraint in the variational
problem~\eqref{eq:variational-problem-for-cate-rate-missing-data},
with $q(1) = \numobs^{\frac{1 + \alpha}{2 (2 + \alpha)}}$. Combining
with the upper bound establishes that
\begin{align*}
\vsuper{\delta_{\xzero}} \asymp \numobs^{\frac{1 + \alpha}{2 (2
    + \alpha)}}.
\end{align*}

\paragraph{Case II: $\xzero \in (0, 1)$.}

For any function $q$ satisfying the constraint in the variational
problem~\eqref{eq:variational-problem-for-cate-rate-missing-data},
~\Cref{lemma:perturb-func-bound-under-cons-for-missing-data-example}
yields
\begin{align*}
|q(\xzero)| \leq c_{\alpha, \xzero} \numobs^{-1/4}, \quad \mbox{for
  constant $c_{\alpha, \xzero}$ depending on $\alpha$ and $\xzero$.}
\end{align*}
On the other hand, given $\numobs \geq \max \big(\xzero^{-2}, (1 -
\xzero)^{-1}\big)$, we construct the function
\begin{align*}
    q(x) = \max \Big\{ \frac{\numobs^{1/4}}{2} -
    \frac{\numobs^{3/4}}{2} |x - \xzero| , 0 \Big\}.
\end{align*}
Straightforward calculation verifies that the function $q$ satisfies
the constraint in the variational
problem~\eqref{eq:variational-problem-for-cate-rate-missing-data},
with $q(\xzero) = \numobs^{1/4} / 2$. Combining with the upper bound
establishes that
\begin{align*}
\vwenlongxzero \asymp \numobs^{1/4}.
\end{align*}

Putting together the results under two cases completes the proof of
equation~\eqref{eq:missing-data-example-cate-rate}.

\subsubsection{Proof of~\Cref{lemma:perturb-func-bound-under-cons-for-missing-data-example}}
\label{subsubsec:proof-lemma-perturb-func-bound-under-cons-for-missing-data-example}

For any $x \in [0, 1]$ and $y \in [0, x]$, by applying the
Cauchy--Schwarz inequality, we find that
\begin{align*}
|q(y)| \geq |q(x)| - |q(y) - q(x)| \geq |q(x)| - \sqrt{ (x - y)
  \int_y^x q' (t)^2 dt} \geq |q(x)| - \sqrt{(x - y) \numobs}.
\end{align*}
Substituting into the second constraint in
equation~\eqref{eq:variational-problem-constraint} yields
\begin{align*}
   1 &\geq \int_0^1 (1 - y)^\alpha q^2(y) dy \geq \int_0^x (1 -
   y)^\alpha \Big[|q(x)| - \sqrt{(x- y) \numobs} \Big]_+^2 dy\\ &\geq
   \int_0^{\frac{q^2(x)}{\numobs}} (z + 1 - x)^\alpha \big(|q(x)| -
   \sqrt{\numobs z} \big)^2 dz \geq \int_0^{\frac{q^2(x)}{2 \numobs}}
   z^\alpha \big(|q(x)| - \sqrt{\numobs z} \big)^2 dz\\ &\geq
   \frac{1}{4 (1 + \alpha)} \cdot \frac{|q(x)|^{4 + 2
       \alpha}}{\numobs^{1 + \alpha}}.
\end{align*}
Since the choice of $x \in [0, 1]$ is arbitrary,
it follows that
\begin{align}
  \sup_{x \in [0, 1]} |q(x)| \leq c_\alpha \numobs^{\frac{1 +
      \alpha}{2 (2 + \alpha)}}.
\end{align}
On the other hand, when $x$ is bounded away from $1$, following the
same derivation, we have
\begin{align*}
   1 \geq \int_0^{\frac{q^2(x)}{\numobs}} (z + 1 - x)^\alpha
   \big(|q(x)| - \sqrt{\numobs z} \big)^2 dz \geq
   \int_0^{\frac{q^2(x)}{2 \numobs}} (1 - x)^\alpha \big(|q(x)| -
   \sqrt{\numobs z} \big)^2 dz \geq \frac{(1 - x)^\alpha}{4 \numobs}
   q^4 (x),
\end{align*}
which implies that
\begin{align}
\sup_{x \in [0, 1 - \varepsilon]} |q(x)| \leq c_{\alpha, \varepsilon}
\numobs^{1/4}
\end{align}
Putting together the pieces completes the proof
of~\Cref{lemma:perturb-func-bound-under-cons-for-missing-data-example}.


\subsubsection{Proof of~\Cref{prop:supnorm-growth-missing-data-example}}
\label{sec:app-growth-condition}

Since we focus on the ratio between $\Ltwospace$-norm and sup-norm of
the eigenfunction $\phi_j$, we slightly abuse the notation, and use
$\phi_j$ to denote a constant multiple of an eigenfunction of
$\kernelfunc$ under $\Ltwospace$ associated to the eigenvalue
$\eigen_j$. The orthogonality condition gives
\begin{align}
\label{eq:kernel-integral-eq-for-growth-cond}  
\eigen_j \phi_j(x) = \int_0^1 \kernelfunc (x, y) \phi_j(y) \propscore
(y) dy.
\end{align}

Substituting the kernel function $\kernelfunc$ the integral
equation~\eqref{eq:kernel-integral-eq-for-growth-cond}, we have
\begin{align*}
\eigen_j \phi_j(x) = \int_0^x y \phi_j(y) \propscore(y) dy + x \cdot
\int_x^1 \phi_j(y) \propscore(y) dy.
\end{align*}
Taking the derivative twice yields the ordinary differential equation
\begin{align*}
    \eigen_j \phi_j''(x) + (1 - x)^\alpha \phi_j(x) = 0 \quad \mbox{on
      $[0, 1]$.}
\end{align*}
Define the auxiliary function $\psi_j (z) \mydefn \phi_j \big(
\eigen_j ^{\frac{1}{2 + \alpha}} (1 - z) \big)$, the differential
equation can be converted into a standard form
\begin{align*}
  \psi_j''(z) + z^\alpha \psi_j = 0, \quad \mbox{for $z \in \big( 0,
    \eigen_j^{\frac{- 1}{2 + \alpha}} \big)$}
\end{align*}
Using $J_\alpha: \real_+ \rightarrow \real$ to denote the Bessel
function of first kind (see~\cite{watson1922treatise}), the ODE above
admits the closed-form solution
\begin{align}
\label{eq:closed-form-eq-for-growth-cond}  
 \psi_j(z) = \sqrt{z} \left\{ \gamma_1 (j) J_{\frac{- 1}{\alpha + 2}}
 \Big( \frac{2 }{\alpha + 2} z^{\alpha / 2 + 1} \Big) + \gamma_2 (j)
 J_{\frac{1}{\alpha + 2}} \Big( \frac{2}{\alpha + 2} z^{\alpha / 2 +
   1} \Big) \right\}, \quad \mbox{for $z \in \big( 0,
   \eigen_j^{\frac{- 1}{2 + \alpha}} \big)$.}
\end{align}
for a pair of constants $\gamma_1 (j)$ and $\gamma_2 (j)$ that may
depend on $j$.

Since we focus on the ratio between $\Ltwospace$-norm and sup-norm of
$\psi_j$, we can assume $\gamma_1^2(j) + \gamma_2^2(j) = 1$ without
loss of generality. Let $\phi_j$ be induced by such function
$\psi_j$. Under this setup, we claim the following relations for any
$j \geq 1$
\begin{subequations}
\begin{align}
    \underline{c} &\leq \vecnorm{\phi_j}{\infty} \leq
    \overline{c}, \label{eq:sup-norm-in-missing-data-example}\\ \underline{c}
    \eigen_j^{\frac{\alpha}{4\alpha + 8}} &\leq
    \vecnorm{\phi_j}{\Ltwospace (\propscore)} \leq \overline{c}
    \eigen_j^{\frac{\alpha}{4\alpha +
        8}}, \label{eq:l2-norm-in-missing-data-example}
\end{align}
\end{subequations}
for a pair $(\underline{c}, \overline{c})$ of constants depending only
on $\alpha$.

Renormalizing the eigenfunction $\phi_j$ to the quantity $\phi_j /
\vecnorm{\phi_j}{\Ltwospace (\probx)}$, we conclude that
\begin{align*}
    \underline{c} / \overline{c} \cdot \eigen_j^{- \frac{\alpha}{4
        (\alpha + 2)}} \leq \vecnorm{\frac{\phi_j}{
        \vecnorm{\phi_j}{\Ltwospace (\probx)}} }{\infty} \leq
    \overline{c} / \underline{c} \cdot \eigen_j^{- \frac{\alpha}{4
        (\alpha + 2)}},
\end{align*}
which proves~\Cref{prop:supnorm-growth-missing-data-example}.

The rest of this section is devoted to the proofs of
equations~\eqref{eq:sup-norm-in-missing-data-example}
and~\eqref{eq:l2-norm-in-missing-data-example}.

\paragraph{Proof of equation~\eqref{eq:sup-norm-in-missing-data-example}:}
For $\alpha > 0$ fixed, by definition, we note have
\begin{align*}
    \inf_{j \geq 1} \vecnorm{\phi_j}{\infty} & \geq \inf_{\gamma_1^2 +
      \gamma_2^2 = 1} \sup_{z \in [1, 2]} \abss{\gamma_1 J_{\frac{-
          1}{\alpha + 2}} \Big( \frac{2 }{\alpha + 2} z^{\alpha / 2 +
        1} \Big) + \gamma_2 J_{\frac{1}{\alpha + 2}} \Big(
      \frac{2}{\alpha + 2} z^{\alpha / 2 + 1} \Big) }\\ &= \sup_{z \in
      [1, 2]} \abss{\gamma_1^* J_{\frac{- 1}{\alpha + 2}} \Big(
      \frac{2 }{\alpha + 2} z^{\alpha / 2 + 1} \Big) + \gamma_2^*
      J_{\frac{1}{\alpha + 2}} \Big( \frac{2}{\alpha + 2} z^{\alpha /
        2 + 1} \Big) },
\end{align*}
where the constants $(\gamma_1^*, \gamma_2^*)$ minimizes the
expression above (the expression is uniformly continuous in
$(\gamma_1, \gamma_2, z)$, which implies continuity of the supremum in
$(\gamma_1, \gamma_2)$, and guarantees existence of a minimizer on a
compact domain). Since Bessel functions $\dudley_{\frac{- 1}{1 +
    \alpha}}$ and $\dudley_{\frac{1}{1 + \alpha}}$ are linearly
independent on any open interval~\cite{watson1922treatise}, there
exists a constant $\underline{c} > 0$ depending only on $\alpha$, such
that
\begin{align*}
  \inf_{j \geq 1} \vecnorm{\phi_j}{\infty} \geq \underline{c}.
\end{align*}
On the other hand, using the asymptotic formulae for Bessel functions,
we note that
\begin{align*}
  \abss{J_{\frac{\pm 1}{\alpha + 2}} \Big( \frac{2 }{\alpha + 2}
    z^{\alpha / 2 + 1} \Big)} = \begin{cases} O \big(1 / \sqrt{z}
    \big) & z \rightarrow 0,\\ O \big( z^{- \frac{\alpha}{4} -
      \frac{1}{2}} \big) & |z| \rightarrow \infty.
  \end{cases}
\end{align*}
Combining with the
expression~\eqref{eq:closed-form-eq-for-growth-cond} implies that the
class of functions $\{\phi_j\}_{j \geq 1}$ admits a uniform upper
bound $\overline{c}$, which is independent of $j$.

\paragraph{Proof of equation~\eqref{eq:l2-norm-in-missing-data-example}:}

Define the auxiliary functions
\begin{align*}
\widetilde{\psi}_j (z) = \bm{1}_{z > 1} \sqrt{\frac{\alpha + 2}{\pi}}
z^{ - \frac{\alpha}{4}} \left\{ \gamma_1 (j) \cos \Big( \frac{2
}{\alpha + 2} z^{\alpha / 2 + 1} - \frac{\alpha \pi}{4 (\alpha + 2)}
\Big) + \gamma_2 (j) \cos \Big( \frac{2 }{\alpha + 2} z^{\alpha / 2 +
  1} + \frac{\alpha \pi}{4 (\alpha + 2)} \Big) \right\}.
\end{align*}
By the asymptotic approximation properties for Bessel
functions~\cite{watson1922treatise}, we have
\begin{align*}
  \abss{\widetilde{\psi}_j (z) - \psi_j (z) } \leq c' (1 + z)^{- 1 -
    \frac{3}{4} \alpha}, \quad \mbox{for any $z \in \real$.}
\end{align*}
for a constant $c' > 0$ depending only on $\alpha$.

Let $\widetilde{\phi}_j \mydefn \psi \big( \eigen_j^{\frac{-1}{2 +
    \alpha}} (1 - x) \big)$, we have
\begin{align}
    \vecnorm{\phi_j - \widetilde{\phi}_j}{\Ltwospace (\propscore)}^2
    &\leq \int_{0}^1 \abss{ \psi_j \big( (1/ \eigen_j)^{\frac{1}{2 +
          \alpha}} (1 - x) \big) - \widetilde{\psi}_j \big( (1 /
      \eigen_j)^{\frac{1}{2 + \alpha}} (1 - x) \big) }^2 \propscore
    (x) dx \nonumber\\ &\leq c' \int_0^{1} \frac{x^\alpha}{1 + \big(
      (1 / \eigen_j)^{\frac{1}{2 + \alpha}} x \big)^{2 + \frac{3}{2}
        \alpha}} dx \nonumber\\ &\leq c_2 (\alpha)
    \eigen_j^{\frac{\alpha + 1}{\alpha +
        2}},\label{eq:approx-error-in-ltwo-norm-for-missing-data-example}
\end{align}
where the constant $c_2 (\alpha)$ depends only on $\alpha$.

For the function $\widetilde{\phi}_j$, we can compute its $\Ltwospace
(\propscore)$-norm.
\begin{align}
    \vecnorm{\widetilde{\phi}_j}{\Ltwospace (\propscore)}^2 &=
    \eigen_j^{\frac{1 + \alpha}{2 + \alpha}} \int_1^{ (\frac{1 }{
        \eigen_j})^{\frac{1}{2 + \alpha}}} \widetilde{\psi}_j^2(z)
    z^\alpha dz \nonumber\\ &= \frac{2}{\alpha + 2} \cdot
    {\eigen_j}^{\frac{1 + \alpha}{2 + \alpha}} \int_1^{1/
      \sqrt{\eigen_j}} \widetilde{\psi}_j^2 \big( \theta^{\frac{2}{2 +
        \alpha}} \big) \theta^{\frac{ \alpha}{\alpha + 2}} d\theta
    \nonumber\\ &= \frac{2}{\pi} \cdot \eigen_j^{\frac{1 + \alpha}{2 +
        \alpha}} \int_1^{1 / \sqrt{\eigen_j}} \left\{ \gamma_1 (j)
    \cos \Big( \frac{2 \theta - \alpha \pi / 4}{\alpha + 2 } \Big) +
    \gamma_2 (j) \cos \Big( \frac{2 \theta + \alpha \pi / 4}{\alpha +
      2}\Big) \right\}^2
    d\theta.\label{eq:expression-for-ltwo-norm-inm-missing-data-example}
\end{align}
Note that the integral is with respect to a periodic function, we have
the upper bound
\begin{align*}
    \int_1^{ 1 / \sqrt{\eigen_j}} \left\{ \gamma_1 (j) \cos \Big(
    \frac{2 \theta - \alpha \pi / 4}{\alpha + 2 } \Big) + \gamma_2 (j)
    \cos \Big( \frac{2 \theta + \alpha \pi / 4}{\alpha + 2}\Big)
    \right\}^2 d\theta \leq c_3 (\alpha) \cdot \eigen_j^{ -
      \frac{1}{2}},
\end{align*}
and the lower bound
\begin{align*}
    &\int_1^{1 / \sqrt{\eigen_j}} \left\{ \gamma_1 (j) \cos \Big(
  \frac{2 \theta - \alpha \pi / 4}{\alpha + 2 } \Big) + \gamma_2 (j)
  \cos \Big( \frac{2 \theta + \alpha \pi / 4}{\alpha + 2}\Big)
  \right\}^2 d\theta \\ &\geq \Big\{\frac{1}{\pi (\alpha +
    2)\sqrt{\eigen_j}} - 2 \Big\} \cdot \int_0^{\pi (\alpha + 2)}
  \left\{ \gamma_1 (j) \cos \Big( \frac{2 \theta - \alpha \pi /
    4}{\alpha + 2 } \Big) + \gamma_2 (j) \cos \Big( \frac{2 \theta +
    \alpha \pi / 4}{\alpha + 2}\Big) \right\}^2 d\theta
\end{align*}
For any pair $(\gamma_1, \gamma_2)$ such that $\gamma_1^2 + \gamma_2^2
= 1$, we have
\begin{align*}
    &\int_0^{2 \pi} \left\{ \gamma_1 \cos \Big( \theta - \frac{\alpha
    \pi / 4}{\alpha + 2 } \Big) + \gamma_2 \cos \Big(\theta +
  \frac{\alpha \pi / 4}{\alpha + 2}\Big) \right\}^2 d\theta\\ &= \pi -
  2 \gamma_1 \gamma_2 \int_0^{2 \pi} \cos \Big( \theta - \frac{\alpha
    \pi / 4}{\alpha + 2 } \Big) \cos \Big(\theta + \frac{\alpha \pi /
    4}{\alpha + 2}\Big) d \theta\\ &= \pi - 2 \pi \gamma_1 \gamma_2
  \cos \big( \frac{\alpha \pi}{2 \alpha + 4} \big) \\ &\geq \pi \Big\{
  1 - \cos \big( \frac{\alpha \pi}{2 \alpha + 4} \big) \Big\},
\end{align*}
which is a positive constant depending only on $\alpha$, and
independent of $\gamma_1$ and $\gamma_2$.

Substituting these bounds back to
equation~\eqref{eq:expression-for-ltwo-norm-inm-missing-data-example}
yields
\begin{align*}
c_4'(\alpha) \eigen_j^{\frac{\alpha}{4\alpha + 8}} \leq
\vecnorm{\widetilde{\phi}_j}{\Ltwospace (\propscore)} \leq c_3'
(\alpha) \eigen_j^{\frac{\alpha}{4\alpha + 8}}
\end{align*}
Combining with
equation~\eqref{eq:approx-error-in-ltwo-norm-for-missing-data-example}
completes the proof of
equation~\eqref{eq:l2-norm-in-missing-data-example}.

\subsection{Proof of~\Cref{cor:continuum-bandit-rates}}
\label{app:subsec-proof-continuum-rates}

We prove the claims about the averaged functional $\taustar$ and the
one-point functional $\LinFuncZero$ separately in the following two
subsections.


\subsubsection{Bounds on the minimax risk for the averaged functional}

By~\Cref{ThmLower,ThmHetero}, we have $\minimaxRisk \big( \ball_\rkhs
(1) \big) \asymp \numobs^{-1} \big\{ \vprobx^2(\treateff) +
\vwenlongsq(\funcClass; \numobs) \big\}$. Taking the measure
$\wenweight{\cdot}{\state} = \delta_\state$ for any $\state
\in \Xspace$, it can be seen that
\begin{align*}
\vprobx^2(\treateff) = \var_{\State \sim \probx} \Big(
\treateff(\State, \pitarget (\State))\Big).
\end{align*}
By the generalized Morrey's embedding theorem (e.g., see \S 5.6.3 in
Evans~\cite{evans2010partial}), for any smoothness
index $\smoothorder > (\dimx + \dima) / 2$, we have
\begin{align*}
\sup_{(\state, \action) \in \statespace \times \actionspace}
|\treateff(\state, \action)| \leq c' \rkhsnorm{\treateff} \leq c',
\end{align*}
for a constant $c'$ depending on the triple $(\dimx, \dima,
\smoothorder)$.  So we have $\vprobx^2(\treateff) \lesssim 1$ in the
worst case.

Since the conditional variance function is a constant, we have
\begin{align*}
  \vwenlongsq = \sum_{j, k \geq 1} \frac{\featweighbar_{j, k}^2}{1 +
    \tfrac{1}{\numobs} \tfrac{1}{\eigen_{j, k}}},
\end{align*}
where we define the projection coefficients
\begin{align*}
\featweighbar_{j, k} = \int_{\Xspace} \phi_j(\state) \psi_k
(\pitarget(\state)) d \state = \inprod{\phi_j}{\psi_k \circ
  \pitarget}_{\Ltwospace (\Xspace)}.
\end{align*}
By the eigenvalue decay condition~\eqref{eq:sobo-eigendecay-general},
we have $\eigen_{j, k}^{-1} \asymp j^{2 \smoothorder / \dimx} + k^{2
  \smoothorder / \dima}$, which implies that
\begin{align}
\label{eq:instance-dependent-bound-for-continuum-bandit-ate}  
\numobs^{-1} \vwenlongsq(\funcClass; \numobs) \asymp
\vwenlonghilbar = \sum_{j, k \geq 1}
\frac{|\inprod{\phi_j}{\psi_k \circ \pitarget}_{\Ltwospace
    (\Xspace)}|^2 }{\numobs + j^{2 \smoothorder / \dimx} + k^{2
    \smoothorder / \dima}},
\end{align}
which proves the instance-dependent bound.

By Parseval's identity, for each $k \geq 1$, we have
\begin{align*}
\sum_{j = 1}^{\infty} \inprod{\phi_j}{\psi_k \circ
  \pitarget}_{\Ltwospace (\Xspace)}^2 = \vecnorm{\psi_k \circ
  \pitarget}{\Ltwospace (\Xspace)}^2 \leq \vecnorm{\psi_k}{\infty}^2 =
1.
\end{align*}
Substituting into the instance-dependent
bound~\eqref{eq:instance-dependent-bound-for-continuum-bandit-ate}, we
have the worst-case instantiation
\begin{align*}
\frac{1}{\numobs} \vwenlongsq(\funcClass; \numobs) & \leq \sum_{k =
  1}^{\numobs^{\frac{\dima}{2 \smoothorder}}} \sum_{j = 1}^{\infty}
\frac{|\inprod{\phi_j}{\psi_k \circ \pitarget}_{\Ltwospace
    (\Xspace)}|^2}{\numobs} + \sum_{k = \numobs^{\frac{\dima}{2
      \smoothorder}}}^{\infty} \sum_{j = 1}^{\infty}
\frac{|\inprod{\phi_j}{\psi_k \circ \pitarget}_{\Ltwospace
    (\Xspace)}|^2}{k^{2 \smoothorder / \dima}} \\
& \leq \Big\{1 + \frac{2 \smoothorder}{2 \smoothorder - \dima} \Big\}
\numobs^{\frac{\dima}{2 \smoothorder} - 1}.
\end{align*}
Thus, we obtain the worst-case upper bound $\sup_{\pitarget}
\minimaxRisk \big( \ball_\rkhs(1) \big) \lesssim
\numobs^{\frac{\dima}{2 \smoothorder} - 1}$.

On the other hand, for any $\actionzero \in \actionspace$, taking the
target functional $\pitarget_{a_0}(\state) \equiv \actionzero$ for any
$\state \in \Xspace$, we have
\begin{align*}
     \numobs^{-1} \vwenlongsq = \sum_{j, k \geq 1} \frac{|\psi_k
       (\actionzero)|^2 |\inprod{\phi_j}{\bm{1}}_{\Ltwospace
         (\Xspace)}|^2 }{\numobs + j^{2 \smoothorder / \dimx} + k^{2
         \smoothorder / \dima}} = \sum_{k \geq 1} \frac{|\psi_k
       (\actionzero)|^2 }{\numobs + 1 + k^{2 \smoothorder / \dima}}
\end{align*}
For the Fourier basis $\psi_k$, we have $|\psi_k (\actionzero)|$ for
any $\actionzero \in \actionspace$, which leads to the lower bound
\begin{align*}
     \sup_{\pitarget} \minimaxRisk \big( \ball_\rkhs (1) \big) \gtrsim
     \sum_{k \geq 1} \frac{1}{\numobs + 1 + k^{2 \smoothorder /
         \dima}} \gtrsim \numobs^{\frac{\dima}{2 \smoothorder} - 1}.
\end{align*}


\subsubsection{Minimax bounds for the one-point functional}

For any $\xzero \in \Xspace$, Theorem~\ref{ThmLower}(b) and
equation~\eqref{EqnCATEclaim} imply that $\minimaxRisk \big( \xzero;
\funcClass \big) \asymp \numobs^{-1} \vwenlongsqxzero$.  By the
variational representation of $\vwenlongsqxzero$, it can be seen that
\begin{align*}
  \vwenlongsqxzero \asymp \sum_{j, k
    \geq 1} \frac{\abss{\phi_j(\xzero) \psi_k (\pitarget
      (\xzero))}^2}{1 + \numobs^{-1} \eigen_{j, k}} = \sum_{j, k
    \geq 1} \frac{1}{1 + \numobs^{-1} \eigen_{j, k}},
\end{align*}
where the last equation follows from the fact that the complex Fourier
bases $\phi_j$ and $\psi_k$ take value at unit circle.

Substituting with the eigenvalue decay
condition~\eqref{eq:sobo-eigendecay-general}, we obtain that
\begin{align*}
   \numobs^{-1} \vwenlongsqxzero \asymp
   \sum_{j, k \geq 1} \frac{1}{\numobs + j^{2 \smoothorder / \dimx} +
     k^{2 \smoothorder / \dima}} =: S_\numobs.
\end{align*}
It remains to study the summation $S_\numobs$. On the one hand, we
note that
\begin{align*}
S_\numobs \geq \sum_{j = 1}^{\numobs^{\frac{\dimx}{2 \smoothorder}}}
\sum_{k = 1}^{\numobs^{\frac{\dima}{2 \smoothorder}}} \frac{1}{\numobs
  + j^{2 \smoothorder / \dimx} + k^{2 \smoothorder / \dima}} \geq
\frac{1}{3} \numobs^{\frac{\dimx + \dima}{2 \smoothorder} - 1}.
\end{align*}
On the other hand, we have the upper bound
\begin{align*}
S_\numobs \leq c_{\dimx, \smoothorder} \sum_{j = 1}^{\infty} \big(
\numobs + j^{2 \smoothorder / \dimx} \big)^{\frac{\dima}{2
    \smoothorder} - 1} \leq c_{\dimx, \dima, \smoothorder}
\numobs^{\frac{\dimx + \dima}{2 \smoothorder} - 1}.
\end{align*}
Therefore, we conclude that $\minimaxRisk (\xzero; \funcClass)
\asymp \numobs^{\frac{\dimx + \dima}{2
    \smoothorder} - 1}$ for any $\xzero \in \Xspace$ and deterministic
policy $\pitarget$.


\end{document}